\documentclass{article}
\usepackage{amssymb}
\usepackage{amsmath}
\usepackage{arxiv}
\usepackage[utf8]{inputenc}
\usepackage[T1]{fontenc}
\usepackage{url}
\usepackage{booktabs}
\usepackage{amsfonts}
\usepackage{nicefrac}
\usepackage{microtype}
\usepackage{color}
\usepackage{lipsum}
\usepackage{amsthm}
\usepackage{indentfirst}
\usepackage{graphicx}
\usepackage{epstopdf}
\usepackage{fourier}
\usepackage{bm}
\usepackage{mathtools}
\usepackage{latexsym,enumerate}
\usepackage{multicol}
\usepackage[skip=2pt]{caption}
\usepackage[font=small,skip=0pt]{subcaption}
\usepackage{array}
\usepackage{ascii}
\usepackage{algorithm,algorithmic}
\usepackage{mathtools}
\usepackage{arydshln}

\usepackage[pdfstartview=FitH, CJKbookmarks=true,
bookmarksnumbered=true, bookmarksopen=true,
colorlinks,
linkcolor=green,
anchorcolor=blue, citecolor=blue
]{hyperref}

\newcommand{\comment}[1]{}

\newcommand{\BEA}{\begin{eqnarray}}
\newcommand{\EEA}{\end{eqnarray}}

\newtheorem{thm}{Theorem}[section]
\newtheorem{proposition}[thm]{Proposition}

\newtheorem{remark}[thm]{Remark}

\newcommand{\PreserveBackslash}[1]{\let\temp=\\#1\let\\=\temp}
\newcolumntype{C}[1]{>{\PreserveBackslash\centering}p{#1}}
\newcolumntype{R}[1]{>{\PreserveBackslash\raggedleft}p{#1}}
\newcolumntype{L}[1]{>{\PreserveBackslash\raggedright}p{#1}}

\newcommand{\stkout}[1]{\ifmmode\text{\sout{\ensuremath{#1}}}\else\sout{#1}\fi}

\begin{document}

\title{Generalized Moving Least-Squares Methods for Solving Vector-valued PDEs on Unknown Manifolds}
\author{
Rongji Li \\
School of Information Science and Technology, ShanghaiTech University, Shanghai
201210, China\\
\texttt{lirj2022@shanghaitech.edu.cn}
\And Qile Yan \\
School of Mathematics, University of Minnesota, 206 Church St SE, Minneapolis, MN 55455, USA\\
\texttt{yan00082@umn.edu}
\And Shixiao W. Jiang \\
Institute of Mathematical Sciences, ShanghaiTech University, Shanghai
201210, China\\
\texttt{jiangshx@shanghaitech.edu.cn}
}
\date{\today}
\maketitle

\begin{abstract}
In this paper, we extend the Generalized Moving Least-Squares (GMLS)
method in two different ways to solve the vector-valued PDEs
on unknown smooth manifolds
without boundaries, identified with randomly
sampled point cloud data. The two approaches are referred to as  the intrinsic method and the extrinsic method.  For the intrinsic method which  relies on local approximations of metric tensors, we simplify the formula of Laplacians and  covariant derivatives acting on vector fields at
the base point by calculating them in a local Monge coordinate system. On the other hand, the extrinsic method formulates tangential derivatives on a submanifold as the projection of the directional derivative in the ambient Euclidean space onto the tangent space of
the submanifold. One challenge of this method
is that the discretization of vector Laplacians yields a matrix whose
size relies on the ambient dimension. To overcome this issue, we reduce the
dimension of vector Laplacian matrices  by employing
a  coordinate transformation. The complexity of both methods scales well with the dimension of manifolds rather than the ambient dimension.
We also present supporting  numerical examples, including eigenvalue problems, linear Poisson equations, and nonlinear Burgers' equations, to examine the numerical accuracy of proposed methods on various smooth manifolds.
\end{abstract}

\lhead{}
\newpage

\section{Introduction}\label{sec:intro}
Solving vector-valued Partial Differential Equations (PDEs) on curved
surfaces or embedded manifolds has various applications in natural sciences
and practical engineering. Examples include computation of parallel transport of tangent vectors \cite{sharp2019vector}, anisotropic texture synthesis \cite{fisher2007design,vaxman2016directional}, simulation of the magnetohydrodynamics \cite{holme2007large},  surface hydrodynamics flows \cite{nitschke2019hydrodynamic,gross2018hydrodynamic}, surface dynamics of fluid vesicles and biomembranes \cite{barrett2015numerical}, fluid flows in computer graphics \cite{stam2003flows}.

For aforementioned applications, there have been several developing numerical
methods to approximate the solution of vector-valued PDEs on manifolds. Both Discrete
Exterior Calculus (DEC) \cite{hirani2003discrete,desbrun2006discrete} and
Finite Element Exterior Calculus (FEEC) \cite%
{arnold2006finite,holst2012geometric,gillette2017finite} are suitable for discretizing the Hodge Laplacian on manifolds. In contrast, the trace FEM \cite{gross2018trace} was developed to address the Bochner Laplacian and a vector Laplacian associated with the surface strain tensor. Given a sufficient high-quality simplicial complex or mesh, those methods can achieve the desirable convergence in PDE problems for surfaces embedded in $\mathbb{R}^3$. However, constructing meshes on surfaces with high quality can be difficult and complicated when the manifold is only identified by randomly sampled data.

To avoid this issue, various mesh-free approaches have been developed for solving PDEs in
recent decades. For unknown manifolds with only a point cloud given, the
first step of most mesh-free approaches is to characterize the manifolds
using a parametrization scheme, such as local SVD methods \cite%
{donoho2003hessian,zhang2004principal,tyagi2013tangent}, level set methods
\cite%
{bertalmio2001variational,greer2006improvement,xu2003eulerian,alvarez2021local}%
, closest point methods \cite{ruuth2008simple,petras2018rbf}, orthogonal
gradient methods \cite{piret2012orthogonal}, moving least squares \cite%
{liang2013solving}, just to name a few. The second step is to approximate
the differential operators acting on  functions in the strong form along
the tangent bundles using approaches including the Radial Basis Function
(RBF) methods \cite{flyer2009radial,Fuselier2009Stability}, Generalized
Moving Least-Squares (GMLS) \cite{liang2013solving,gross2020meshfree}, Generalized Finite
Difference Methods (GFDM) \cite{suchde2019fully,suchde2021meshfree}. In our previous work \cite{harlim2023radial}, the RBF method was extended to provide an accurate approximation of differential operators acting on vector fields of smooth manifolds. However, its effectiveness depends crucially on the choice of shape parameters in the radial basis function. Additionally, the computational cost of the method is high since the resulting discrete operators are dense matrices with their size dependent on the ambient dimension. On the other hand, a local mesh method was proposed in  \cite{lai2013local} which discretizes the weak formulation of differential operators using local connectivity between points. This approach was recently extended to handle vector Laplacians on manifolds in \cite{peoples2024higher}.

Alternatively, graph-based approaches, including Vector Diffusion Maps (VDM)
\cite{singer2012vector} and Spectral Exterior Calculus (SEC) \cite%
{berry2020spectral} can also be used as a mesh-free tool to
approximate vector Laplacians for processing raw data in manifold
learning applications. While these methods do not require a parametrization of the manifold, it is limited to estimating the Laplacian-type operators.

\comment{
So far, each of the aforementioned approaches has its own advantages
depending on the application scenario as well as having challenges.}

In this paper, we consider solving vector-valued PDEs on  $d$-dimensional unknown
manifolds $M$ without boundaries embedded in the Euclidean space $\mathbb{R}^{n}$, identified with randomly sampled point cloud data. Our goal here is to construct an approximation for differential operators acting on vector fields of smooth manifolds where the resulting matrices are sparse with size dependent of the intrinsic dimension $d$ and the number of data $N$. Basically, we translate the formulations from differential geometry into representations that can be concretely and efficiently  realized by vector and matrix algebra.
Based on two different formulations of vector differential operators, we extend the GMLS methods in two distinct ways to approximate these operators and subsequently solve vector-valued PDEs on manifolds.
These two approximations will be referred to as the intrinsic method and the extrinsic method. Our framework is quite flexible for approximating different differential operators such as covariant derivative and various vector Laplacians on simple smooth manifolds.

\comment{
Firstly, we employ the GMLS for the parametrization of manifolds by
locally fitting multivariate polynomial functions of the local tangent space
coordinates which are first approximated by local SVD methods \cite%
{liang2013solving,suchde2019meshfree,gross2020meshfree,jones2023generalized}%
, and subsequently apply the GMLS to approximate the tangential derivatives
and Laplacians of target functions or vector fields on the point
cloud. In order to ensure the stability and convergence of PDEs, we
choose the same specific weight functions as in \cite{liang2013solving} for
the local least-squares regression problems.}

\comment{The GMLS approaches for surface PDEs can be grouped into two classes. The
first class is referred to as intrinsic methods, which in general involves
with the local characterization of the Riemannian metric tensor and other
geomtric quantities of manifolds in local intrinsic coordinates \cite%
{liang2013solving,gross2020meshfree,jones2023generalized}.}

For the intrinsic GMLS method, it generally involves the local characterization of the Riemannian metric tensor and other
geometric quantities of manifolds in local intrinsic coordinates \cite%
{liang2013solving,gross2020meshfree,jones2023generalized}. In this paper, we derive the formulas of the covariant derivative and Laplacians acting on vector fields at
the base point in a local Monge coordinate system. It turns out that the
resulting formula can be much simplified by using the properties of the
Monge coordinate system so that it helps to reduce the computational cost
in constructing the vector Laplacian matrices. Basically we find appropriate local and global representations of vector/tensor fields, apply GMLS to each of their associated components, and then assemble these components into a $dN\times dN$ discrete vector Laplacian matrix. Here and throughout, $N$ denotes the size of the point cloud data
, $n$ is the ambient dimension and $d$ is the intrinsic dimension of the manifold.
Our derivations can be viewed
as a generalization from the Laplace-Beltrami of functions in \cite%
{liang2013solving} to the covariant derivative and Laplacians of
vector fields.

For the extrinsic GMLS method, \comment{the surface PDE and its associated variables and operators are extended to
be defined on an open neighborhood of the manifold so that}the surface
derivative in the local coordinate can be written as a projection of the
appropriate derivative in the ambient space \cite%
{suchde2019meshfree,suchde2021meshfree,jiang2024generalized}. Then the extrinsic method relies on the discretization in the ambient
space which can be computationally expensive for the vector field case. In this paper, we apply a coordinate transformation method to reduce the matrix dimensionality of discrete vector differential operators from  $nN \times nN$ to $dN \times dN$. This dimension reduction can be numerically realized since  all objects (e.g., functions and tensor fields) in this paper are confined to the
tangent bundle of manifolds and the
surfaces  are all static without movement in the normal
direction. In fact, this work improves the result in our previous work \cite{harlim2023radial}, saving the
memory cost and computational cost especially when the ambient dimension $n\sim 10$ is moderately large.


\comment{
 Here, we remark that the GMLS for both intrinsic and extrinsic
methods can be applied to the PDE problems with at most $d=1$ to $4$ and $%
n=10$ to $100 $ due to the curse of dimensionality of multivariate polynomial
basis functions.}


The paper is organized as follows. In Section \ref{sec:prel}, we provide a
short review of the GMLS approach and its application in the approximation of
tangent vectors on unknown manifolds identified by randomly sampled data
points. In Section \ref{sec:intf}, we present the intrinsic formulation for
approximating the Bochner Laplacian and the covariant derivative using GMLS. We distinguish our discussion between the cases of 2D surfaces in $\mathbb{R}^3$ and general $d$-dimensional manifolds in $\mathbb{R}^n$. In Section \ref{sec:extf}, using the extrinsic formulation, we apply a coordinate transformation to enhance the efficiency of approximation  for the Bochner
Laplacian and the covariant derivative. In Section \ref{sec:numr}, we provide
supporting numerical examples for solving vector-valued PDEs, particularly eigenvalue problems, screened Poisson equations, linear vector diffusion equations and nonlinear Burgers' equations. In Section \ref{sec:con}, we conclude the paper with a summary and some open problems.  We relegate the intrinsic approximation and extrinsic approximation of two other vector Laplacians
to the supplementary material Section~\ref{app:B} and Section~\ref{app:A}, respectively.


\section{Preliminaries\label{sec:prel}}

\subsection{Review of GMLS}\label{sec:gmlsm}

\begin{figure*}[tbph]
	\centering
	\includegraphics[width=3.8in, height=1.8in]{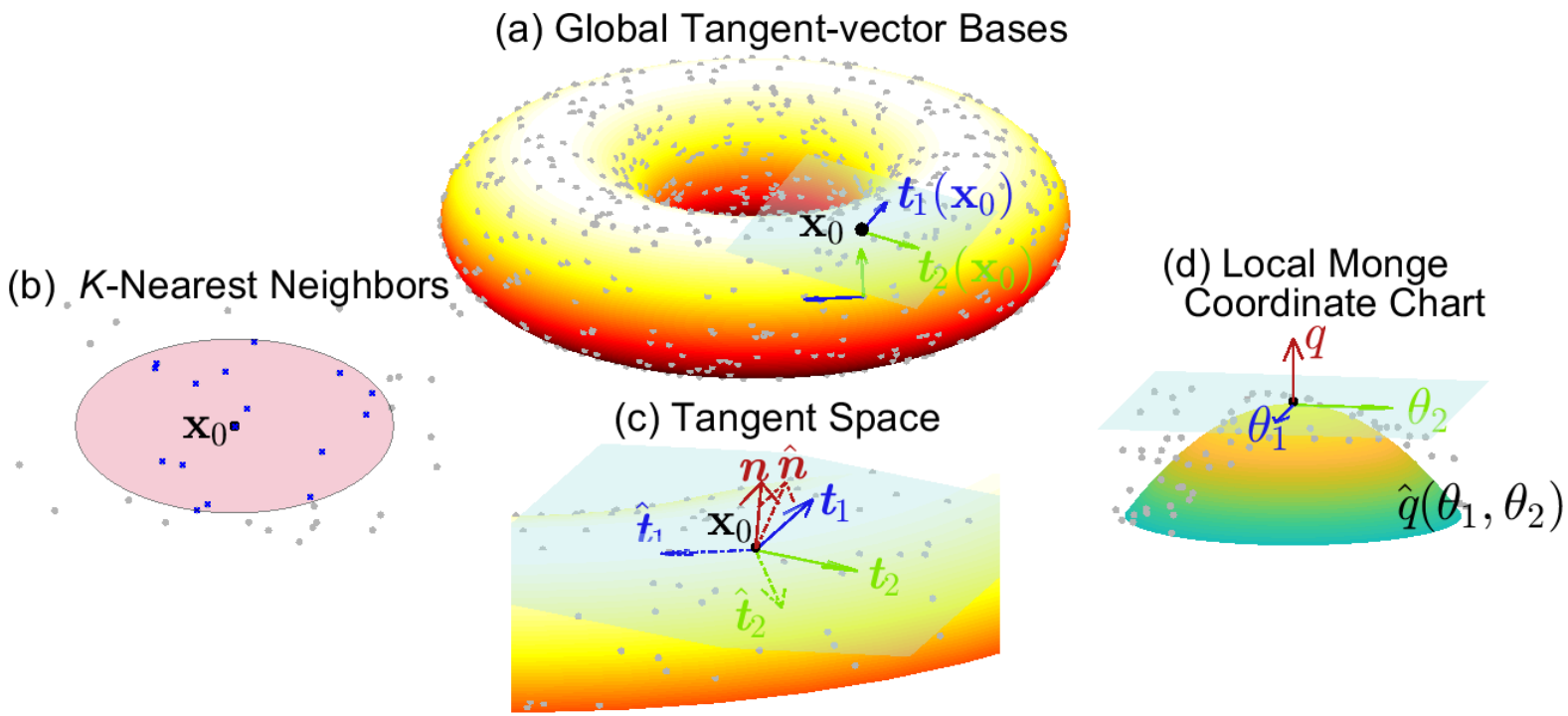}
\caption{Sketch of our setups. (a) Global tangent-vector bases $\{\boldsymbol{t}_1,\boldsymbol{t}_2\}$ on the point
cloud $\{\mathbf{x}_{i}\}_{i=1}^{N}$ (gray dots) without alignment.
 (b) $K$-nearest neighbors of the base
point $\mathbf{x}_{0}$ plotted as blue crosses inside the red circle. (c)
The tangent space $\mathrm{span}\{\boldsymbol{t}_1,\boldsymbol{t}_2\}$ and
the normal $\boldsymbol{n}$ at the base point $\mathbf{x}_{0}$. Also plotted are their estimations
$\{\boldsymbol{\hat{t}}_1,\boldsymbol{\hat{t}}_2,
\boldsymbol{\hat{n}}\}$ given in (\ref{eq:that}) when the manifold is unknown. (d) The
local Monge coordinate $(\protect\theta_1,\protect\theta_2)$
and the GMLS approximation of the manifold $\hat{q}(\theta_1,\theta_2)$ in the local coordinate system (equation (\ref{eqn:qhat})).  }
\label{fig1_sketch}
\end{figure*}

The  GMLS approach is a regression
technique for approximating functions by a linear combination of local basis
functions, such as local polynomials, based on a set of scattered data
samples (see e.g. \cite%
{nayroles1992generalizing,liu1997moving,levin1998approximation,Wendland2005Scat,mirzaei2012generalized}).
We now briefly
review the application of GMLS approach on a 2D smooth manifold embedded in $%
\mathbb{R}^{3}$. Given a set of (distinct) nodes $\mathbf{X}_{M}=\{\mathbf{x}%
_{i}\}_{i=1}^{N}\subset M$. For a base point ${\mathbf{x}_{0}}\in \mathbf{X}%
_{M}\subset M$, we denote its $K$-nearest neighbors in $\mathbf{X}_{M}$ by $%
S_{{\mathbf{x}_{0}}}=\{{\mathbf{x}_{0,k}}\}_{k=1}^{K}\subset \mathbf{X}_{M}$ (Fig.~\ref{fig1_sketch}(b)). In most literature, such a set $S_{{\mathbf{x}_{0}}}$ is called a
\textquotedblleft stencil\textquotedblright. By definition, we have ${%
\mathbf{x}_{0,1}}={\mathbf{x}_{0}}$ to be the base point. We first assume that the
true tangent space of the manifold $M$ is given and then
discuss the unknown manifold regime in the subsection \ref%
{sec:unkM}. Let $T_{\mathbf{x}}M$ be the tangent space at point $\mathbf{x}%
=(x^{1},x^{2},x^{3})\in M$ and denote its orthonormal basis by $\left\{
\boldsymbol{t}_{1}(\mathbf{x}),\boldsymbol{t}_{2}(\mathbf{x})\right\} $.

\comment{Notice that based on Hairy ball theorem, there is no nonvanishing continuous
tangent vector field\ on a 2D sphere and thus the orthonormal basis $\{\boldsymbol{%
t}_{1}(\mathbf{x}),\boldsymbol{t}_{2}(\mathbf{x})\}$ cannot be chosen as
continuous vector fields on it. To avoid this issue, we take a basis of tangent vector at each node ${\mathbf{x}}_{i}\in \mathbf{X}_{M}$,
denoted by $\left\{ \boldsymbol{t}_{1}(\mathbf{x}_{i}),\boldsymbol{t}_{2}(%
\mathbf{x}_{i})\right\} _{i=1}^{N},$ which will not be aligned (see
Fig. \ref{fig1_sketch}(a)). That is to say, the chosen basis $%
\left\{ \boldsymbol{t}_{1}(\mathbf{x}),\boldsymbol{t}_{2}(\mathbf{x}%
)\right\} $\ can take arbitrary orthonormal rotation in the tangnet space $%
T_{\mathbf{x}}M$ at each node.}


The GMLS approach in this paper is formulated similarly as in \cite%
{liang2013solving,gross2020meshfree,jones2023generalized,jiang2024generalized}%
. First, let $\mathbb{P}_{\mathbf{x}_{0}}^{l}$ be the space of
\textquotedblleft intrinsic\textquotedblright\ polynomials with degree up to
$l$ in two variables at the base ${\mathbf{x}_{0}}$, i.e.,
\begin{equation}
\mathbb{P}_{\mathbf{x}_{0}}^{l}=\mathrm{span}\{\theta _{1}^{\alpha
_{1}}\theta _{2}^{\alpha _{2}}|0\leq \left\vert \boldsymbol{\alpha }%
\right\vert\leq l\},  \label{eq:thet12}
\end{equation}%
where $\theta _{i}(\mathbf{x})=\boldsymbol{t}_{i}({\mathbf{x}_{0}})\cdot (\mathbf{x}-{%
\mathbf{x}_{0}})$ $(i=1,2)$ are the local coordinates
and $\boldsymbol{\alpha }=(\alpha _{1},\alpha _{2})\in \mathbb{N}^{2}$ is
the mult-index notation with $\left\vert \boldsymbol{\alpha }\right\vert =\alpha _{1}+\alpha _{2}$. By definition, the dimension of the space $\mathbb{P%
}_{\mathbf{x}_{0}}^{l}$ is $m=(l+2)(l+1)/2$. Given $K>m$ and fitting data $%
\mathbf{f}_{{\mathbf{x}_{0}}}=(f({\mathbf{x}_{0,1}}),...,f(\mathbf{x}%
_{0,K}))^{\top }$, we can define an operator $\mathcal{I}_{\mathbb{P}}:%
\mathbf{f}_{{\mathbf{x}_{0}}}\in \mathbb{R}^{K}\rightarrow \mathcal{I}_{%
\mathbb{P}}\mathbf{f}_{{\mathbf{x}_{0}}}\in \mathbb{P}_{\mathbf{x}_{0}}^{l}$
such that $\mathcal{I}_{\mathbb{P}}\mathbf{f}_{{\mathbf{x}_{0}}}$ is the
optimal solution of the following moving least-squares problem:%
\begin{equation}
\underset{\hat{f}\in \mathbb{P}_{\mathbf{x}_{0}}^{l}}{\min }%
\sum_{k=1}^{K}\lambda _{k}\left( f({\mathbf{x}_{0,k}})-\hat{f}({\mathbf{x}%
_{0,k}})\right) ^{2},  \label{eqn:int_LS}
\end{equation}%
where $\lambda _k$ is some positive weight function. There are various
choices of weighted $\ell ^{2}$-norms in the moving least-squares (\ref%
{eqn:int_LS}) as shown in {\cite%
{Wendland2005Scat,liang2013solving,gross2020meshfree}}.
As reported in previous works {\cite{liang2012geometric,liang2013solving,jiang2024generalized}}, the weight function
\begin{equation}
\lambda _{k}=\left\{ \begin{aligned} 1, \quad \mathrm{if} \quad k &= 1, \\
1/K, \quad \mathrm{if} \quad k & = 2,\ldots,K, \end{aligned}\right.
\label{kkernel}
\end{equation}%
can numerically provide a stable approximation to the Laplace-Beltrami operator across various data sets  \cite{liang2012geometric,liang2013solving,jiang2024generalized}. 
Here we also apply such a special weight function
for the numerically stable vector Laplacians.
For the covariant derivative, we choose the uniform weight function for a stable approximation.

The solution to the moving
least-squares problem (\ref{eqn:int_LS}) can be represented as $\hat{f}(%
\mathbf{x}):=(\mathcal{I}_{\mathbb{P}}\mathbf{f}_{{\mathbf{x}_{0}}})(\mathbf{%
x})=\sum_{0\leq \alpha _{1}+\alpha _{2}\leq l}b_{\alpha _{1},\alpha
_{2}}\theta _{1}^{\alpha _{1}}(\mathbf{x})\theta _{2}^{\alpha _{2}}(\mathbf{x%
})$, where the concatenated coefficients $\boldsymbol{b}=(b_{\boldsymbol{%
\alpha }(1)},\ldots ,b_{\boldsymbol{\alpha }(m)})^{\top }\in \mathbb{R}%
^{m\times 1}$ can be solved from a normal equation:
\begin{equation}
\boldsymbol{b}=(\boldsymbol{\Phi }^{\top }\boldsymbol{\Lambda \Phi })^{-1}%
\boldsymbol{\Phi }^{\top }\boldsymbol{\Lambda }\mathbf{f}_{{\mathbf{x}_{0}}%
}:=\boldsymbol{\Phi }^{\dag }\mathbf{f}_{{\mathbf{x}_{0}}},  \label{eq:phinv}
\end{equation}%
with $\boldsymbol{\Lambda }:=\mathrm{diag}(\lambda _{1},\cdots ,\lambda
_{K})\in \mathbb{R}^{K\times K}$ and $\boldsymbol{\Phi }^{\dag }:=(%
\boldsymbol{\Phi }^{\top }\boldsymbol{\Lambda \Phi })^{-1}\boldsymbol{\Phi }%
^{\top }\boldsymbol{\Lambda }\in \mathbb{R}^{m\times K}$. Here $\left\{
\boldsymbol{\alpha }(j)\right\} _{j=1}^{m}$ is the ordered multi-index $%
\boldsymbol{\alpha }(j)=(\alpha _{1}(j),\alpha _{2}(j))$ with $0\leq \alpha
_{1}(j)+\alpha _{2}(j)\leq l$ for all $j=1,\ldots ,m$. The above normal
equation can be uniquely identified if $\boldsymbol{\Phi }=(\boldsymbol{\Phi
}_{kj})\in \mathbb{R}^{K\times m}$ with
\begin{equation}
\boldsymbol{\Phi }_{kj}=\theta _{1}^{\alpha _{1}(j)}({\mathbf{x}_{0,k}}%
)\theta _{2}^{\alpha _{2}(j)}({\mathbf{x}_{0,k}}),\quad 1\leq k\leq K,\
1\leq j\leq m,  \label{eqn:matrix_phi}
\end{equation}%
is a full rank matrix.


%

We now present the theoretical error estimates \cite{mirzaei2012generalized} for such a GMLS
approximation for partial derivative operators. For a set of points $\mathbf{X}_M=\{\mathbf{x}%
_{i}\}_{i=1}^{N}\subset M$, define the fill distance $%
h_{\mathbf{X},M}$ and the separation distance $q_{\mathbf{X},M}$ in Euclidean space as
\begin{equation}
h_{\mathbf{X},M}=\sup_{\mathbf{x}\in M}\min_{1\leq j\leq N}\Vert \mathbf{x}-\mathbf{x}%
_{j}\Vert _{2},\text{ \ \ }q_{\mathbf{X},M}=\frac{1}{2}\min_{i\neq j}\Vert \mathbf{x}%
_{i}-\mathbf{x}_{j}\Vert _{2},  \notag
\end{equation}%
respectively, and assume that $q_{\mathbf{X},M}\leq h_{\mathbf{X},M}\leq c_{q}q_{\mathbf{X},M}$ for some
constant $c_{q}\geq 1$. Let $\mathbf{X}_{M}\subseteq M$ be such that $%
h_{\mathbf{X},M}\leq h_{0}$ for some constant $h_{0}>0$ and $M^{\ast
}=\bigcup_{\mathbf{x}\in M}B(\mathbf{x},C_{2}h_{0})$, a union of geodesic
ball over the length $C_{2}h_{0}$ with some constant $C_{2}>0$. Here, the basic idea is to guarantee both the local Taylor expansion and the local polynomial fitting to be valid in the local stencil $B(\mathbf{x},C_{2}h_{0})$ for each $\mathbf{x}\in M$. Then for any
$f\in C^{l+1}(M^{\ast })$ and $\boldsymbol{\delta }=(\delta _{1},\delta
_{2})$ that satisfies $|{\boldsymbol{\delta }}|=\delta _{1}+\delta
_{2}\leq l$,
\begin{equation}
\big|D^{{\boldsymbol{\delta }}}f(\mathbf{x})-\widehat{D^{{\boldsymbol{\delta
}}}f}(\mathbf{x})\big|\leq Ch_{\mathbf{X},M}^{l+1-|{\boldsymbol{\delta }}%
|}\max_{|\boldsymbol{\beta} |=l+1}\Vert D^{\boldsymbol{\beta} }f\Vert _{L^{\infty }(M^{\ast })},
\label{MLS_error}
\end{equation}%
for all $\mathbf{x}\in M$ and some $C>0$ (see Section 4 of \cite{mirzaei2012generalized}). In the error bound above, we used
the notation $\widehat{D^{{\boldsymbol{\delta }}}f}$ as the GMLS
approximation to $D^{{\boldsymbol{\delta }}}f$ using local polynomials up to
degree $l$, where $D^{{\boldsymbol{\delta }}}$ denotes a general
multi-dimensional derivative with multi-index ${\boldsymbol{\delta }}$. %
For uniformly i.i.d. random data {$\mathbf{X}_{M}\subseteq M$ }on manifolds,
with probability higher than $1-\frac{1}{N}$, one can show that $h_{\mathbf{X},M}=O((\log N/N)^{1/2})$ for a 2D manifold \cite{jiang2024generalized}.
Consequently, one has the error bound for a two dimensional manifold that
$ |D^{{\boldsymbol{\delta }}}f(\mathbf{x})-\widehat{D^{{\boldsymbol{\delta }}}f}(\mathbf{x})|=O\left(\left({\log N}/{N}\right)^{\frac{l+1-|{\boldsymbol{\delta }}|}{2}}\right) $.

\begin{remark}\label{rem:err2}
For a $d$-dimensional manifold, one has
$h_{\mathbf{X},M}=O((\log N/N)^{1/d})$ \cite{jiang2024generalized}. Utilizing the analysis presented in  \cite{mirzaei2012generalized}, the following error bound can be established:
\begin{equation}
|D^{{\boldsymbol{\delta }}}f(\mathbf{x})-\widehat{D^{{\boldsymbol{\delta }}}f}(\mathbf{x})|=O\left(\left(\frac{\log N}{N}\right)^{\frac{l+1-|{\boldsymbol{\delta }}|}{d}}\right). \label{eq:gmlrt}
\end{equation}
Although the
analysis in \cite{mirzaei2012generalized,jiang2024generalized} was not originally developed for vector fields in the manifold setting, our numerical results (Section~\ref{sec:numr}) show that the convergence rates are in good agreement with the error bound (\ref{eq:gmlrt}).
\end{remark}

\subsection{Approximation of tangent vectors for unknown manifolds\label%
{sec:unkM}}

For the unknown manifold setup, we assume that we are only given a point
cloud $\{\mathbf{x}_{i}\}_{i=1}^{N}\subset M$\ on the interior of the smooth
manifold. To construct approximations of differential operators, we need to approximate bases of tangent vector space $\left\{ \boldsymbol{t}%
_{1}(\mathbf{x}),\boldsymbol{t}_{2}(\mathbf{x})\right\} $ at each point $\{\mathbf{x}_{i}\}_{i=1}^{N}$. We first use the local singular
value decomposition (SVD) or local principal component analysis (PCA) for a
coarse approximation of the tangent vectors {\cite%
{donoho2003hessian,zhang2004principal,tyagi2013tangent,liang2013solving,harlim2023radial}%
}:

\begin{enumerate}
\item For each base point ${\mathbf{x}_{0}}\in \mathbb{R}^{3}$and its $K$%
-nearest neighbors $\{{\mathbf{x}_{0,k}}\}_{k=1}^{K}$, construct the
distance matrix $\mathbf{D}:=[\mathbf{D}_{1},\ldots ,\mathbf{D}_{K}]\in
\mathbb{R}^{3\times K},$ where $\mathbf{D}_{i}={\mathbf{x}_{0,i}-\mathbf{x}%
_{0}}$ ($i=1,\ldots ,K$).

\item Take a singular value decomposition of $\mathbf{D}=\mathbf{V}_{L}%
\boldsymbol{\Sigma }\mathbf{V}_{R}$, and then obtain the leading 2-columns
of $\mathbf{V}_{L}\in \mathbb{R}^{3\times 3}$, denoted by $\left\{
\boldsymbol{\tilde{t}}_{1}(\mathbf{x}_{0}),\boldsymbol{\tilde{t}}_{2}(%
\mathbf{x}_{0})\right\} $, which approximates a tangent basis of $T_{\mathbf{%
x}_0}M$, as well as the last column of $\mathbf{V}_{L}$, denoted by $%
\boldsymbol{\tilde{n}}({\mathbf{x}_{0}})$, which approximates a normal
vector (Fig. \ref{fig1_sketch}(c)).
\end{enumerate}

Once the coarse approximation to the tangent vectors is obtained, one can apply the GMLS approach to get a more accurate approximation of the
manifold as well as the tangent space  {\cite{liang2013solving,gross2020meshfree,jones2023generalized}}. In particular, one can construct a
local coordinate system with the origin $\mathbf{x}_{0}$\ and the
orthonormal basis $\{\boldsymbol{\tilde{t}}_{1},\boldsymbol{\tilde{t}}_{2},%
\boldsymbol{\tilde{n}}\}$. Then, a neighboring point $\mathbf{x}$\ around $%
\mathbf{x}_{0}$\ has the local coordinate $\left( \tilde{\theta}_{1},\tilde{%
\theta}_{2},\tilde{q}\right) $, where $\tilde{\theta}_{i}({\mathbf{x}})=%
\boldsymbol{\tilde{t}}_{i}({\mathbf{x}_{0}})\cdot (\mathbf{x}-{\mathbf{x}_{0}%
})$ $(i=1,2)$ and $\tilde{q}({\mathbf{x}})=\boldsymbol{\tilde{n}}({\mathbf{x}%
_{0}})\cdot (\mathbf{x}-{\mathbf{x}_{0}})$. An approximation of the manifold with higher accuracy could be calculated as follows:
\begin{enumerate}
\setcounter{enumi}{2}
\item Let $\bar{q}(\tilde{\theta}_{1},\tilde{\theta}_{2})=\sum_{1\leq \alpha
_{1}+\alpha _{2}\leq l}\tilde{a}_{\alpha _{1},\alpha _{2}}\tilde{\theta}%
_{1}^{\alpha _{1}}\tilde{\theta}_{2}^{\alpha _{2}}$ be the optimal solution
to the GMLS problem:%
\begin{equation*}
\underset{\bar{q}\in \mathbb{P}_{\mathbf{x}_{0}}^{l}}{\min }%
\sum_{k=1}^{K}\lambda _{k}\Big( \bar{q}\big( \tilde{\theta}_{1}({\mathbf{x}%
_{0,k}}),\tilde{\theta}_{2}({\mathbf{x}_{0,k}})\big) -\tilde{q}({\mathbf{x}%
_{0,k}})\Big) ^{2}.
\end{equation*}
\item Solve the normal equation $(\boldsymbol{\tilde{\Phi}}^{\top }%
\boldsymbol{\Lambda \tilde{\Phi}})\boldsymbol{\tilde{a}}=\boldsymbol{\tilde{%
\Phi}}^{\top }\boldsymbol{\Lambda }\mathbf{\tilde{q}}_{{\mathbf{x}_{0}}}$ to
obtain the concatenated regression coefficients $\boldsymbol{\tilde{a}}%
=\left( \tilde{a}_{\alpha _{1}(j),\alpha _{2}(j)}\right) _{0\leq \alpha
_{1}+\alpha _{2}\leq l}^{j=1,\ldots ,m}\in \mathbb{R}^{m\times 1},$ where $%
\mathbf{\tilde{q}}_{{\mathbf{x}_{0}}}=(\tilde{q}({\mathbf{x}_{0,1}}),...,%
\tilde{q}(\mathbf{x}_{0,K}))^{\top }$ and $\boldsymbol{\tilde{\Phi}}$ is in
the form of (\ref{eqn:matrix_phi}) but with $\{\theta _{1},\theta _{2}\}$
replaced by $\{\tilde{\theta}_{1},\tilde{\theta}_{2}\}$. Thus, one arrives
at the locally parametrized surface, $\bar{q}(\tilde{\theta}_{1},\tilde{%
\theta}_{2})$.

\item The two tangent vectors are $(1,0,\partial _{\tilde{\theta}_{1}}\bar{q}%
)^{\top }$\ and $(0,1,\partial _{\tilde{\theta}_{2}}\bar{q})^{\top }$\ in
the local coordinate system with the basis $\{\boldsymbol{\tilde{t}}_{1},%
\boldsymbol{\tilde{t}}_{2},\boldsymbol{\tilde{n}}\}.$
Then compute the two tangent vectors of $\bar{q}(\tilde{\theta}_{1},\tilde{%
\theta}_{2})$\ at the base $\mathbf{x}_{0}$\ in the Cartesian coordinate
system as: \
\begin{equation}
\boldsymbol{\hat{t}}_{k}({\mathbf{x}_{0}})=\boldsymbol{\tilde{t}}_{k}({%
\mathbf{x}_{0}})+\frac{\partial \bar{q}(0,0)}{\partial \tilde{\theta}_{k}}%
\boldsymbol{\tilde{n}}({\mathbf{x}_{0}})=\left\{
\begin{array}{cc}
\boldsymbol{\tilde{t}}_{k}({\mathbf{x}_{0}})+\tilde{a}_{1,0}\boldsymbol{%
\tilde{n}}({\mathbf{x}_{0}}), & k=1, \\
\boldsymbol{\tilde{t}}_{k}({\mathbf{x}_{0}})+\tilde{a}_{0,1}\boldsymbol{%
\tilde{n}}({\mathbf{x}_{0}}), & k=2.%
\end{array}%
\right.  \label{eq:that}
\end{equation}
Then $\{\hat{\boldsymbol{t}}_1(\mathbf{x}_0),\hat{\boldsymbol{t}}_2(\mathbf{x}_0)\}$ is orthonormalized using Gram-Schmidt process.
\item Compute the normal vector $\boldsymbol{\hat{n}}=\boldsymbol{\hat{t}}%
_{1}\times \boldsymbol{\hat{t}}_{2}$ at $\mathbf{x}_{0}$. Repeat the above
steps for each base point in $\{\mathbf{x}_{i}\}_{i=1}^{N}$.
\end{enumerate}

\begin{remark}
For a $d$ dimensional manifold $M$ embedded in $n$
dimensional ambient spaces $\mathbb{R}^{n}$ ($d\ll n$), a similar strategy can be used to approximate the $d$ tangent vectors $\{ \boldsymbol{\hat{t}}%
_{i}\}_{i=1}^d$ at a given point. Specifically, we first apply the local SVD as described in above Step 1 and Step 2 to obtain a coarse approximation of the tangent vectors  $\{\boldsymbol{\tilde{t}}%
_{i}\}_{i=1}^d$.
Next, we define the local coordinates $\tilde{\theta}_{i}({\mathbf{x}})=%
\boldsymbol{\tilde{t}}_{i}({\mathbf{x}_{0}})\cdot (\mathbf{x}-{\mathbf{x}_{0}%
}), 1\leq i \leq d$. Rather than applying the GMLS approximation to the local normal coordinate $\tilde{q}_i(\mathbf{x})=\tilde{\boldsymbol{n}}_i(\mathbf{x}_0)\cdot(\mathbf{x}-\mathbf{x}_0), 1\leq i\leq n-d,$ as done in above 2D case, we instead apply the GMLS to each ambient coordinate of $\mathbf{x}-\mathbf{x}_0$. This approach avoided the need for generating $n-d$ normal directions $\tilde{\boldsymbol{n}}_i(\mathbf{x}_0)$ and computing $n-d$ inner products   in each $\tilde{q}_i(\mathbf{x})$ at each point, significantly reducing computational cost when $d\ll n$. A GMLS approximation of $\{ \boldsymbol{\hat{t}}%
_{i}\}_{i=1}^d$ with higher accuracy can then be calculated as follows:
\begin{enumerate}	
	\setcounter{enumi}{2}
	\item Denote ${\mathbf{x}}=(x^{1},\ldots ,x^{n})$ and ${\mathbf{x}}%
	_{0}=(x_{0}^{1},\ldots ,x_{0}^{n})$. Let
	\begin{small}
	\begin{equation*}
		\bar{y}^{s}(\tilde{\theta}_{1},\ldots ,\tilde{\theta}_{d})=\bar{x}^{s}(%
		\tilde{\theta}_{1},\ldots ,\tilde{\theta}_{d})-x_{0}^{s}=\sum_{1\leq
			\left\vert \boldsymbol{\alpha }\right\vert \leq l}\tilde{b}_{\boldsymbol{%
				\alpha }}^{s}\boldsymbol{\tilde{\theta}}^{\boldsymbol{\alpha }}:=\sum_{1\leq
			\left\vert \boldsymbol{\alpha }\right\vert \leq l}\tilde{b}_{\alpha
			_{1},\ldots ,\alpha _{d}}^{s}\tilde{\theta}_{1}^{\alpha _{1}}\cdots \tilde{%
			\theta}_{d}^{\alpha _{d}},
	\end{equation*}%
	\end{small}
	be the optimal solution of the GMLS problem to approximate $%
	y^{s}=x^{s}-x_{0}^{s}$ for each dimension $s=1,\ldots ,n$. Here $\boldsymbol{%
		\alpha }=(\alpha _{1},\ldots ,\alpha _{d})$\ and $\left\vert \boldsymbol{%
		\alpha }\right\vert =\sum_{i=1}^{d}\alpha _{i}$.
	
	\item Solve the normal equation to obtain the regression coefficients $\{%
	\tilde{b}_{\boldsymbol{\alpha }}^{s}\}_{1\leq \left\vert
		\boldsymbol{\alpha }\right\vert \leq l}\ $for each $s=1,\ldots ,n$. The complexity is $O(nKm^{2}+nm^{3})$ at each ${\mathbf{x}}_{0}$, where $%
		m=\left(
		\begin{array}{c}
		l+d \\
		d%
		\end{array}%
		\right) $ is the number of monomial
		basis functions $\left\{ \boldsymbol{%
			\tilde{\theta} }^{\boldsymbol{\alpha }}|0\leq \left\vert \boldsymbol{\alpha }%
		\right\vert \leq l\right\} $.
	
	\item Compute the $d$ tangent vectors of the parametrized manifold $\bar{x}%
	^{s}(\tilde{\theta}_{1},\ldots ,\tilde{\theta}_{d})$\ at the base $\mathbf{x}%
	_{0}$\ in the Cartesian coordinate system as: \
	\begin{equation}
		\boldsymbol{\hat{t}}_{k}({\mathbf{x}_{0}})=\sum_{s=1}^{n}\frac{\partial \bar{%
				x}^{s}(0,\ldots ,0)}{\partial \tilde{\theta}_{k}}\mathbf{e}%
		_{s}=\sum_{s=1}^{n}\tilde{b}_{\boldsymbol{\alpha }(k)}^{s}\mathbf{e}_{s},
		\label{eq:that22}
	\end{equation}%
	where $\boldsymbol{\alpha }(k)=(0,\ldots ,1,0,\ldots ,0)$ with only the $k$%
	th entry being $1$ and $\mathbf{e}_{s}$ is the $s$th standard
	orthonormal basis of the Cartesian coordinate system. Then we orthonormalize $\{\hat{\boldsymbol{t}}_1(\mathbf{x}_0),\ldots,\hat{\boldsymbol{t}}_d(\mathbf{x}_0)\}$ using Gram-Schmidt process.
\end{enumerate}
Notably, we do not need to generate the $n-d$ normal vectors $\{\hat{\boldsymbol{n}}_i\}_{i=1}^{n-d}$, which would require a computational complexity of about $O(n^3)$ when $d\ll n$. Now, the total complexity for  Step 3 to Step 5 in the above algorithm is $O(n)$.
We also point out that the moving least-squares based approach proposed in \cite{sober2020manifold} can be used to handle  point cloud data in $\mathbb{R}^n$.

\end{remark}

\begin{figure*}[htbp]
	{\scriptsize \centering
		\setlength{\tabcolsep}{1pt}
	\begin{tabular}{ccc}
		{\normalsize (a) 2D Torus in $\mathbb{R}^3$} & {\normalsize (b) 2D Torus in $\mathbb{R}^9$} & {\normalsize (c) 3D Flat Torus in $\mathbb{R}^{12}$} \\
		\includegraphics[width=1.65in, height=1.35in]{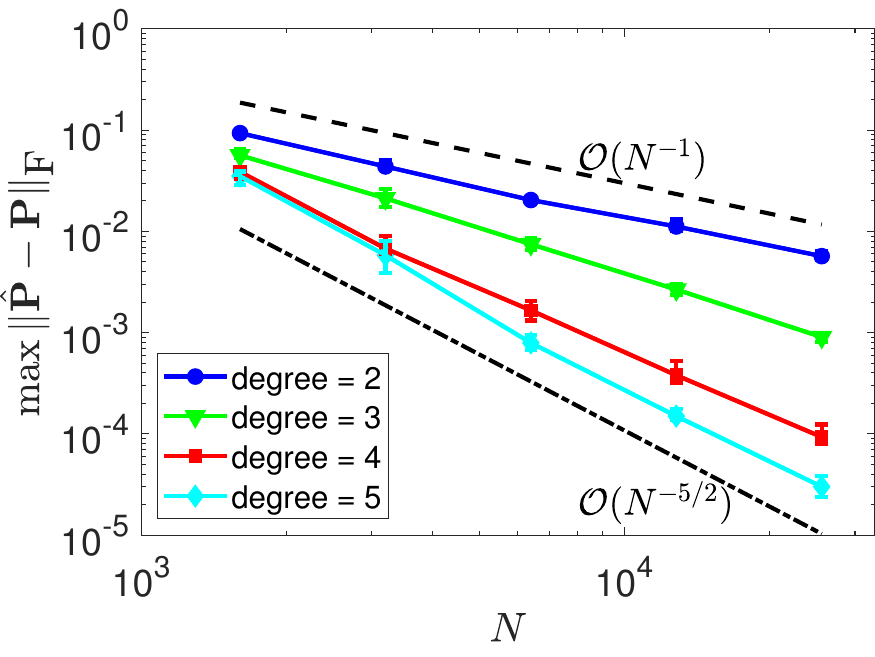} &
		\includegraphics[width=1.65in, height=1.35in]{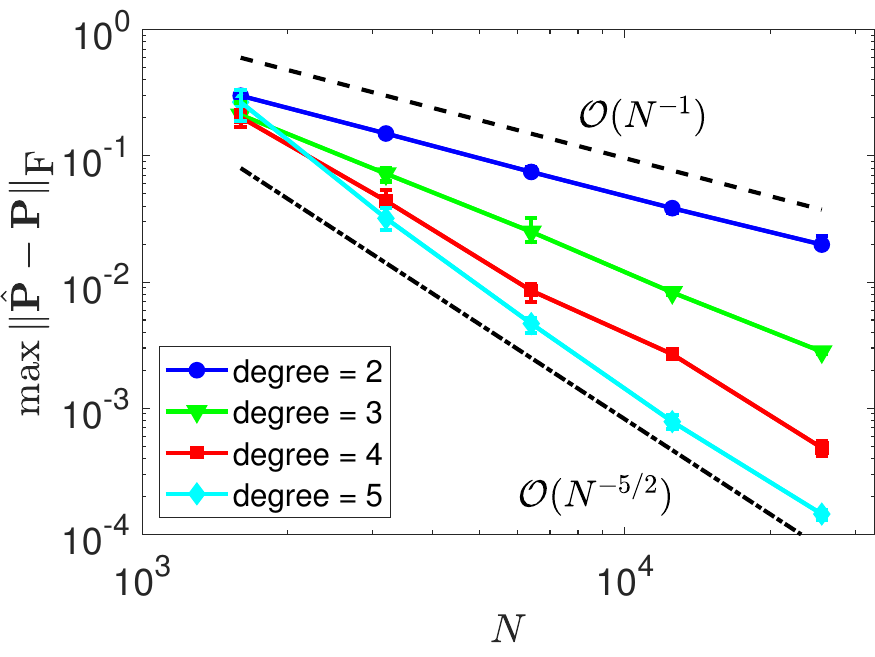} &
		\includegraphics[width=1.65in, height=1.35in]{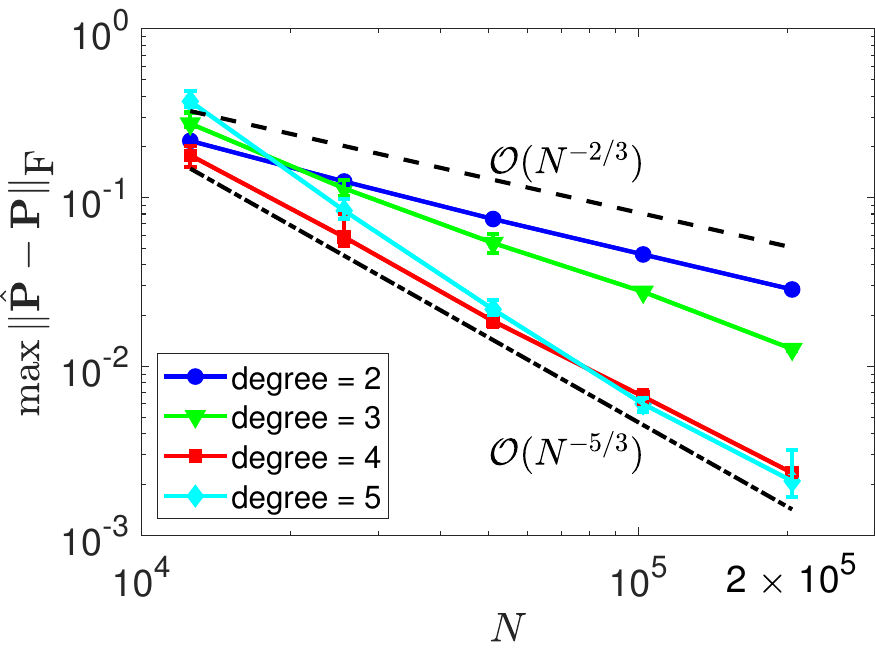}
	\end{tabular}
	}
	\caption{Convergence results for the estimations of projection matrices $\hat{\mathbf{P}}$ on various manifolds. (a) 2D torus in $\mathbb{R}^3$. (b) 2D torus in $\mathbb{R}^9$. (c) 3D flat torus in $\mathbb{R}^{12}$. Here, $\Vert \cdot \Vert_{\mathrm{F}}$ denotes the Frobenius matrix norm, the projection matrix is ${\mathbf{P}}=\sum_{i=1}^d\boldsymbol{t}_{i}\boldsymbol{t}^\top_{i}$ and its estimation is  $\hat{\mathbf{P}}=\sum_{i=1}^d\boldsymbol{\hat{t}}_{i}\boldsymbol{\hat{t}}^\top_{i}$ (see Section~\ref{sec:revext} or \cite{harlim2023radial} for more about the projection matrix).  }
	\label{fig_tang}
\end{figure*}

In Figs. \ref{fig_tang}(a)-(c), we numerically verify the
convergence rate for the estimation of the projection matrices $\hat{\mathbf{P}}=\sum_{i=1}^d\boldsymbol{\hat{t}}_{i}\boldsymbol{\hat{t}}^\top_{i}$  on various examples, including a 2D torus in $\mathbb{R}^3$, a 2D torus in $\mathbb{R}^9$ and a 3D torus in $\mathbb{R}^{12}$. Their parametrizations are defined in (\ref{eqn:torusr3}), (\ref{eqn:torus}), and (\ref{eqn:flattorus}), respectively.
As shown in Fig. \ref{fig_tang}, the convergence
rate is approximately $O(N^{-{l}/{d}})$, which aligns well with the theoretical GMLS error bound of $O((\log N/ N)^{{l}/{d}})$ in (\ref{eq:gmlrt}).

\comment{This agreement is expected
since the normal $\boldsymbol{\hat{%
n}}$\ and the two tangent $\boldsymbol{\hat{t}}_{1},\boldsymbol{\hat{t}}_{2}$%
\ involve first-order derivatives of the parametrized surface $\bar{q%
}(\tilde{\theta}_{1},\tilde{\theta}_{2})$.}


\section{GMLS using intrinsic formulation\label{sec:intf}}
In this section, we approximate differential operators for vector fields using  intrinsic formulation. For 2D manifolds in $\mathbb{R}^3$, we first review GMLS for the local Monge
patch parametrization of each stencil $S_{{\mathbf{x}_{0}}}$ \cite%
{jones2023generalized,gross2020meshfree} in Section~\ref{sec:mgchart}, and subsequently derive the numerical scheme to approximate the Bochner Laplacian acting on a vector field in Section~\ref{sec:ABL}. We then generalize our approach to approximating the Bochner Laplacian  on $d$-dimensional manifolds embedded in $\mathbb{R}^n$  in Section~\ref{3p3highint}, where we also provide a complexity analysis. In Section~\ref{sec:CD}, we approximate the covariant derivative for 2D manifolds, noting that the formulation naturally generalizes to higher dimensions. For readers not familiar with intrinsic differential geometry, please refer to Section~\ref{app:intr_DG} for a short review.

\comment{
Notice again that each $\boldsymbol{t}_{j}$%
is not continuous and thereafter nor
is each corresponding component $%
u^{j}$ due to Hairy ball theorem stated
above. While  each $\boldsymbol{t}_{j}$
and correspondingly each component $u^{j}$%
are continuous, their linear
combination, the vector field $u$, is continuous
on a manifold.
since may not be cont. on some manifolds such
as sphere, we choose an arbitrary randomly-oriented tangent vector
basis $t_j$.}

\subsection{Monge parametrization}\label{sec:mgchart}


We first review the Monge parametrization \cite%
{monge1809application,pressley2010elementary} for a two-dimensional manifold
in $\mathbb{R}^{3}$, which has been used in \cite%
{liang2013solving,jones2023generalized,gross2020meshfree} to approximate the
local metric tensor and other local geometric quantities. Assume that we are
given a true tangent vector basis $\{\boldsymbol{t}_{1}({\mathbf{x}_{0}}),%
\boldsymbol{t}_{2}({\mathbf{x}_{0}})\}$ for known manifolds (see Fig. \ref{fig1_sketch}(a)) or given the approximate tangent vector basis $\{%
\boldsymbol{\hat{t}}_{1}({\mathbf{x}_{0}}),\boldsymbol{\hat{t}}_{2}({\mathbf{%
x}_{0}})\}$\ as in (\ref{eq:that}) for unknown manifolds (see Fig. \ref{fig1_sketch}(c)) at the base point $%
{\mathbf{x}_{0}}$. For the ease of discussion, we use $\{\boldsymbol{t}_{1}({%
\mathbf{x}_{0}}),\boldsymbol{t}_{2}({\mathbf{x}_{0}})\}$ to denote the tangent
vector basis for both cases of known and unknown manifolds in this section. Thereafter one can obtain the normal vector $\boldsymbol{n}(%
{\mathbf{x}_{0}})=\boldsymbol{t}_{1}({\mathbf{x}_{0}})\times \boldsymbol{t}%
_{2}({\mathbf{x}_{0}})$.

A local coordinate chart for the manifold near the base point
could be defined using the embedding map $\boldsymbol{q}$,%
\begin{equation*}
\boldsymbol{q}(\theta _{1},\theta _{2};{\mathbf{x}_{0}})={\mathbf{x}_{0}}%
+\theta _{1}\boldsymbol{t}_{1}({\mathbf{x}_{0}})+\theta _{2}\boldsymbol{t}%
_{2}({\mathbf{x}_{0}})+q\left( \theta _{1},\theta _{2}\right) \boldsymbol{n}(%
{\mathbf{x}_{0}}),
\end{equation*}%
where $(\theta _{1},\theta _{2})$ are local coordinates of $\mathbf{x}$ as
defined in (\ref{eq:thet12}) and $q\left( \theta _{1},\theta _{2}\right)=%
\boldsymbol{n}({\mathbf{x}_{0}})\cdot (\mathbf{x}-{\mathbf{x}_{0}})$. Here $%
q $ is a smooth function over the $(\theta _{1},\theta _{2})$-plane\ which
denotes the \textquotedblleft height\textquotedblright\ of the manifold
surface (see Fig. \ref{fig1_sketch}(d)). This is known as the Monge
parametrization of manifold \cite%
{monge1809application,pressley2010elementary}.

Then the \textquotedblleft height\textquotedblright\ $q$ is approximated as a
polynomial function of the local coordinates $(\theta _{1},\theta _{2})$
using the GMLS approach following the similar procedure introduced in Section %
\ref{sec:unkM}. Here one key difference is the degrees of all
\textquotedblleft intrinsic\textquotedblright\ polynomials in the approximation
space ranging from $2$ to $l$, that is, the approximation space is $\mathbb{P}_{%
\mathbf{x}_{0}}^{l}\backslash \mathbb{P}_{\mathbf{x}_{0}}^{1}=\mathrm{span}%
(\{\theta _{1}^{\alpha _{1}}\theta _{2}^{\alpha _{2}}\}_{2\leq \alpha
_{1}+\alpha _{2}\leq l})$. Thus, the approximation of $q$ is the optimal
solution to the least-squares problem (\ref{eqn:int_LS}):
\begin{equation}
\hat{q}=\sum_{2\leq \alpha _{1}+\alpha _{2}\leq l}a_{\alpha _{1},\alpha
_{2}}\theta _{1}^{\alpha _{1}}\theta _{2}^{\alpha _{2}}.  \label{eqn:qhat}
\end{equation}%
Note here that the constant term and two linear terms w.r.t. $(\theta
_{1},\theta _{2})$ are absent in $\hat{q}$ (Fig. \ref%
{fig1_sketch}(d)) since the parameterized surface $%
\hat{q}$  should pass through the base point ${\mathbf{x}_{0}}$ and possess
the tangent space $T_{\mathbf{x}_{0}}M=\mathrm{span}\{\boldsymbol{t}_{1}({%
\mathbf{x}_{0}}),\boldsymbol{t}_{2}({\mathbf{x}_{0}})\}$.


Observe that the coordinate representation for tangent vector bases in the
vicinity of ${\mathbf{x}_{0}}$\ are $(1,0,\partial _{\theta _{1}}\hat{q}%
)^{\top }$\ and $(0,1,\partial _{\theta _{2}}\hat{q})^{\top }$\ in the local
basis $\{\boldsymbol{t}_{1},\boldsymbol{t}_{2},\boldsymbol{n}\}$ (Fig. \ref%
{fig1_sketch}(d)). For $\mathbf{x}$ close to ${\mathbf{x}_{0}}$, the tangent
vector basis $\{\boldsymbol{t}_{1},\boldsymbol{t}_{2}\}\in T_{\mathbf{x}}M$
in the Cartesian coordinate system of ambient space $\mathbb{R}^3$ are approximated as
\begin{equation}
\partial _{k}:=\frac{\partial }{\partial \theta _{k}}=\boldsymbol{t}_{k}({%
\mathbf{x}_{0}})+\frac{\partial \hat{q}}{\partial \theta _{k}}\boldsymbol{n}(%
{\mathbf{x}_{0}}),\text{ \ \ }k=1,2.  \label{eqn:pak}
\end{equation}%
The local Riemannian metric $\boldsymbol{g} = (g_{ij})_{i,j=1}^2$ for ${\mathbf{x}}\in S_{{\mathbf{x}_{0}}}$\ is
approximated as%
\begin{equation}
\boldsymbol{g}=\left[
\begin{matrix}
1+\left( \partial _{\theta _{1}}\hat{q}\right) ^{2} & \partial _{\theta _{1}}%
\hat{q}\partial _{\theta _{2}}\hat{q} \\
\partial _{\theta _{1}}\hat{q}\partial _{\theta _{2}}\hat{q} & 1+\left(
\partial _{\theta _{2}}\hat{q}\right) ^{2}%
\end{matrix}%
\right] .  \label{eq:gm}
\end{equation}%
As usual, we use $g^{ij}$ to denote the inverse metric tensor. 

\subsection{Approximation of the Bochner Laplacian}\label{sec:ABL}

\comment{In Section~\ref{sec:mgc}, we apply
GMLS to locally approximate the Bochner Laplacian $\Delta _{B}$ acting on a
vector field $u=\tilde{u}^{i}\partial _{i}$\ at the base point ${\mathbf{x}%
_{0}}$\ using the $K$-nearest neighbors in each $S_{{\mathbf{x}_{0}}}$,
where $\{\partial _{1},\partial _{2}\}$ is the associated local
tangent-vector basis of Monge parametrization (Fig. \ref{fig1_sketch}(d)).\
We repeat this procedure for each node in $\{\mathbf{x}_{i}\}_{i=1}^{N}$ so
that we have $N$ Monge patches associated with$\ N$ distinct local bases $%
\{\partial _{1},\partial _{2}\}$. Then we need a global tangent vector basis
in order to represent the vector field $u$ and the Bochner $\Delta _{B}u$\
at each point in $\{\mathbf{x}_{i}\}_{i=1}^{N}$ with one unified basis. In
particular, we choose $\left\{ \boldsymbol{t}_{1},\boldsymbol{t}_{2}\right\}
$ as the global basis (Fig. \ref{fig1_sketch}(a)) and then we can write the
vector field $u=u^{j}\boldsymbol{t}_{j}$ pointwisely on $\{\mathbf{x}%
_{i}\}_{i=1}^{N}$. In Section~\ref{sec:tranf}, we will discuss how to
transform the components of the vector field and the Laplacian matrix from a
local basis to a global one in order to approximate $\Delta _{B}u$\ on $\{%
\mathbf{x}_{i}\}_{i=1}^{N}$ in the form of matrix-vector product.}

In this section, we assume that $M$ is a two-dimensional manifold
in $\mathbb{R}^{3}$. For any vector field given on the point cloud $\mathbf{X}_M\subset M$,  we write it as $\boldsymbol{u}=u^{j}%
\boldsymbol{t}_{j}$ where the components $\{u^{1},u^2\}$\ are defined on the
nodes $\{\mathbf{x}_{i}\}_{i=1}^{N}$ and  $\left\{ \boldsymbol{t}_{1}(\mathbf{x}_{i}),\boldsymbol{t}_{2}(%
\mathbf{x}_{i})\right\} _{i=1}^{N}$ are global bases of tangent vector space at each node. Here and throughout this paper, we use the Einstein summation convention. Note that there is no need to do alignments of these bases (as shown in Fig. \ref{fig1_sketch}(a)). To obtain a finite-difference type approximation of the Bochner Laplacian at a base point $\mathbf{x}_0$, we intend to find the weights $\{%
\boldsymbol{w}^{k}\in \mathbb{R}^{2\times 2}\}_{k=1}^{K}$ such that,
\begin{equation}
\Delta_{B}\boldsymbol{u}({\mathbf{x}_{0}})\approx\sum_{k=1}^{K}\boldsymbol{w}^{k}\boldsymbol{u}|_{{\mathbf{x}_{0,k}}}=(\boldsymbol{w}%
^{1},\ldots ,\boldsymbol{w}^{K})\mathbf{u}_{{\mathbf{x}_{0}}}, \label{eqn:LB}
\end{equation}%
where $\{{\mathbf{x}_{0,k}}\}_{k=1}^{K}$ is the $K$-nearest neighbors of $\mathbf{x}_0$, $\boldsymbol{u}|_{{\mathbf{x}_{0,k}}}=(u^{1}({\mathbf{x}_{0,k}}),u^{2}({\mathbf{x}%
_{0,k}}))^{\top }$ and $\ \mathbf{u}_{{%
\mathbf{x}_{0}}}=(u^{1}({\mathbf{x}_{0,1}}),u^{2}({\mathbf{x}_{0,1}}),\ldots
,u^{1}({\mathbf{x}_{0,K}}),u^{2}({\mathbf{x}_{0,K}}))^{\top }\in \mathbb{R}%
^{2K\times 1}.$

In Section~\ref{sec:mgc}, we derive the intrinsic formulation for the Bochner Laplacian in the local Monge coordinate system  (\ref{eqn:pak}). In Section~\ref{sec:comp}, we apply the GMLS to discretize the operator
in the local Monge coordinate.
In Section~\ref{sec:MaBL}, we discuss how to obtain the weights in (\ref{eqn:LB}) and assemble these weights into a sparse $2N\times 2N$ matrix where
$2$ corresponds to the intrinsic dimension of the surface.


\subsubsection{Approximation of geometric quantities in the local Monge coordinate}\label{sec:mgc}

The computations for geometric quantities and function derivatives  involved in the Bochner Laplacian $\Delta_B\boldsymbol{u}$ could be complicated. Fortunately, here we only focus on the values of $\Delta_B\boldsymbol{u}$ at the base point ${\mathbf{x}_{0}}$ and their calculations can be much simplified in the Monge coordinate system as introduced in Section %
\ref{sec:mgchart}. In particular, we use the following properties of a Monge coordinate system:
\begin{proposition}
\label{prop:norm} Let $\{U_0;\xi^{1},\ldots ,\xi^{d}\}$ be a Monge
coordinate system at $\mathbf{x}_{0}$\ where $U_0\subset M$ is an open
neighborhood of $\mathbf{x}_{0}$ and $M$ is a $d$-dimensional manifold with Riemannian metric $\boldsymbol{g}$. Denote  the Christoffel symbol by $\Gamma _{ij}^{k}$. Then%
\newline
\noindent (1) For all $1\leq i,j\leq d$, $g_{ij}(\mathbf{x}_{0})=\delta
_{ij} $.\newline
\noindent (2) For all $1\leq i,j,k\leq d$, $\Gamma _{ij}^{k}(\mathbf{x}
_{0})=0$.\newline
\noindent (3) For all $1\leq i,j,k\leq d$, $\partial _{k}g_{ij}(\mathbf{x}
_{0})=0$.
\end{proposition}
\noindent These properties can be easily verified by noticing that in the expression
of the metric (\ref{eq:gm}), the polynomial $\hat{q}$ as defined in (\ref%
{eqn:qhat}) has no constant term or linear terms w.r.t. $\{\theta
_{1},\theta _{2}\}$.

The Bochner Laplacian acting on a  vector field $\boldsymbol{u}=\tilde{u}^i\partial_i$ is defined as
$$
\Delta_B\boldsymbol{u}=g^{ij}\tilde{u}^k_{;ij}\partial_k,
$$
where the components are computed by
$$
\tilde{u}^i_{;jk}=\frac{\partial \tilde{u}^i_{;j}}{\partial \theta_k}+\tilde{u}^p_{;j}\Gamma_{kp}^i-\tilde{u}^i_{;p}\Gamma^p_{kj},\quad  \tilde{u}^i_{;j}=\frac{\partial \tilde{u}^i}{\partial\theta_j}+\tilde{u}^k\Gamma^i_{jk}.
$$
Applying Proposition \ref{prop:norm} and noticing that the tangent plane
bases in (\ref{eqn:pak}) are $\partial _{k}|_{\mathbf{x}_{0}}=\boldsymbol{t}%
_{k}({\mathbf{x}_{0}})$, one can simplify the Bochner Laplacian $\Delta
_{B}\boldsymbol{u}(\mathbf{x}_{0})$ as,
\begin{eqnarray}
&&\Delta _{B}\boldsymbol{u}(\mathbf{x}_{0})=g^{ij}\tilde{u}_{;ij}^{k}\partial _{k}|_{%
\mathbf{x}_{0}}=\left[ \left( \tilde{u}_{;11}^{1}+\tilde{u}_{;22}^{1}\right)
\partial _{1}+\left( \tilde{u}_{;11}^{2}+\tilde{u}_{;22}^{2}\right) \partial
_{2}\right] \big{|}_{\mathbf{x}_{0}},  \label{eqn:ux0b} \\
&=&\big\{\left[ \partial _{11}\tilde{u}^{1}+\partial _{22}\tilde{u}^{1}+%
\tilde{u}^{1}(\partial _{1}\Gamma _{11}^{1}+\partial _{2}\Gamma _{12}^{1})+%
\tilde{u}^{2}(\partial _{1}\Gamma _{21}^{1}+\partial _{2}\Gamma _{22}^{1})%
\right] \boldsymbol{t}_{1}  \notag \\
&&+\left[ \tilde{u}^{1}(\partial _{1}\Gamma _{11}^{2}+\partial _{2}\Gamma
_{12}^{2})+\partial _{11}\tilde{u}^{2}+\partial _{22}\tilde{u}^{2}+\tilde{u}%
^{2}(\partial _{1}\Gamma _{21}^{2}+\partial _{2}\Gamma _{22}^{2})\right]
\boldsymbol{t}_{2}\big\}\big|_{\mathbf{x}_{0}},  \notag
\end{eqnarray}%
\comment{
where the components $\tilde{u}_{;st}^{i}=\frac{\partial ^{2}\tilde{u}^{i}}{\partial \theta
_{s}\partial \theta _{t}}+\frac{\partial \tilde{u}^{m}}{\partial \theta _{t}}%
\Gamma _{ms}^{i}+\tilde{u}^{m}\frac{\partial \Gamma _{ms}^{i}}{\partial
\theta _{t}}+\frac{\partial \tilde{u}^{k}}{\partial \theta _{s}}\Gamma
_{kt}^{i}+\tilde{u}^{m}\Gamma _{ms}^{k}\Gamma _{kt}^{i}-\frac{\partial
\tilde{u}^{i}}{\partial \theta _{k}}\Gamma _{st}^{k}-\tilde{u}^{m}\Gamma
_{st}^{k}\Gamma _{km}^{i}$ and can be simplified at $\mathbf{x}_{0}$ as $%
\tilde{u}_{;st}^{i}|_{\mathbf{x}_{0}}=(\frac{\partial ^{2}\tilde{u}^{i}}{%
\partial \theta _{s}\partial \theta _{t}}+\tilde{u}^{m}\frac{\partial \Gamma
_{ms}^{i}}{\partial \theta _{t}})|_{\mathbf{x}_{0}}$ using Proposition \ref%
{prop:norm}. }
where we simplify the components $\tilde{u}_{;jk}^{i}$  at $\mathbf{x}_{0}$  using Proposition \ref{prop:norm} as
$$
\tilde{u}^i_{;jk}|_{\mathbf{x}_{0}}=\Big(\frac{\partial \tilde{u}^i_{;j}}{\partial\theta_k}+\tilde{u}^p_{;j}\Gamma_{kp}^i-\tilde{u}^i_{;p}\Gamma^p_{kj}\Big)\Big|_{\mathbf{x}_{0}}=\frac{\partial \tilde{u}_{;j}^i}{\partial\theta_k}\Big|_{\mathbf{x}_{0}}=\Big(\frac{\partial ^{2}\tilde{u}^{i}}{%
\partial \theta _{j}\partial \theta _{k}}+\tilde{u}^{m}\frac{\partial \Gamma
_{mj}^{i}}{\partial \theta _{k}}\Big)\Big|_{\mathbf{x}_{0}}.
$$
For the Christoffel symbol $\Gamma _{kl}^{i}=\frac{1}{2}%
g^{im}\left( \frac{\partial g_{mk}}{\partial \theta _{l}}+\frac{\partial
g_{ml}}{\partial \theta _{k}}-\frac{\partial g_{kl}}{\partial \theta _{m}}%
\right) $, one can calculate its derivative at the base $\mathbf{x}_{0}$
using Proposition \ref{prop:norm} as%
\begin{equation*}
\frac{\partial \Gamma _{kl}^{i}}{\partial \theta _{s}}\bigg|_{\mathbf{x}%
_{0}}=\frac{1}{2}g^{im}\left( \frac{\partial g_{mk}}{\partial \theta
_{s}\partial \theta _{l}}+\frac{\partial g_{ml}}{\partial \theta
_{s}\partial \theta _{k}}-\frac{\partial g_{kl}}{\partial \theta
_{s}\partial \theta _{m}}\right) \bigg|_{\mathbf{x}_{0}}.
\end{equation*}%
Writing $g_{ij}$ in terms of the function $\hat{q}$ (see equation (\ref{eq:gm})), one can further simplify the derivative of the Christoffel symbol. For example, one can compute
\begin{eqnarray}
\partial _{1}\Gamma _{12}^{1}|_{\mathbf{x}_{0}} &=&\frac{1}{2}g^{11}\left(
\frac{\partial g_{11}}{\partial \theta _{1}\partial \theta _{2}}+\frac{%
\partial g_{12}}{\partial \theta _{1}\partial \theta _{1}}-\frac{\partial
g_{12}}{\partial \theta _{1}\partial \theta _{1}}\right) |_{\mathbf{x}_{0}}=%
\frac{1}{2}\frac{\partial }{\partial \theta _{2}}\left( 2\frac{\partial \hat{%
q}}{\partial \theta _{1}}\frac{\partial ^{2}\hat{q}}{\partial \theta
_{1}\partial \theta _{1}}\right) |_{\mathbf{x}_{0}}  \notag \\
&=&\left( \frac{\partial ^{2}\hat{q}}{\partial \theta _{1}\partial \theta
_{2}}\frac{\partial ^{2}\hat{q}}{\partial \theta _{1}\partial \theta _{1}}%
\right) |_{\mathbf{x}_{0}}=2a_{2,0}a_{1,1},  \label{eq:3der}
\end{eqnarray}%
where $a_{\alpha _{1},\alpha
_{2}}$ are the coefficients of $\hat{q}$ as in (\ref{eqn:qhat}). Similarly, one can obtain all the derivatives of Christoffel symbols at $%
\mathbf{x}_{0}$\ for 2D manifolds in $\mathbb{R}^3$:%
\begin{equation*}
\begin{array}{llll}
\partial _{1}\Gamma _{11}^{1}=4(a_{2,0})^{2}, & \partial _{2}\Gamma
_{11}^{1}=2a_{1,1}a_{2,0}, & \partial _{1}\Gamma _{11}^{2}=2a_{2,0}a_{1,1},
& \partial _{2}\Gamma _{11}^{2}=4a_{0,2}a_{2,0}, \\
\partial _{1}\Gamma _{12}^{1}=2a_{2,0}a_{1,1}, & \partial _{2}\Gamma
_{12}^{1}=(a_{1,1})^{2}, & \partial _{1}\Gamma _{12}^{2}=(a_{1,1})^{2}, &
\partial _{2}\Gamma _{12}^{2}=2a_{0,2}a_{1,1}, \\
\partial _{1}\Gamma _{22}^{1}=4a_{0,2}a_{2,0}, & \partial _{2}\Gamma
_{22}^{1}=2a_{0,2}a_{1,1}, & \partial _{1}\Gamma _{22}^{2}=2a_{0,2}a_{1,1},
& \partial _{2}\Gamma _{22}^{2}=4(a_{0,2})^{2},%
\end{array}%
\end{equation*}%
where the restriction notations $|_{\mathbf{x}_{0}}$ are all omitted. Then
one can calculate $\Delta _{B}\boldsymbol{u}(\mathbf{x}_{0})$ in (\ref{eqn:ux0b}) as
\begingroup\makeatletter\def\f@size{9}\check@mathfonts
\begin{eqnarray}
\Delta _{B}\boldsymbol{u}(\mathbf{x}_{0}) &=&\left( \partial _{11}\tilde{u}^{1}+\partial
_{22}\tilde{u}^{1}+\tilde{u}^{1}\left[ 4(a_{2,0})^{2}+(a_{1,1})^{2}\right] +%
\tilde{u}^{2}\left[ 2a_{2,0}a_{1,1}+2a_{0,2}a_{1,1}\right] \right)
\boldsymbol{t}_{1}({\mathbf{x}_{0}})  \notag \\
&+&\left( \tilde{u}^{1}[2a_{2,0}a_{1,1}+2a_{0,2}a_{1,1}]+\partial _{11}%
\tilde{u}^{2}+\partial _{22}\tilde{u}^{2}+\tilde{u}^{2}\left[
4(a_{0,2})^{2}+(a_{1,1})^{2}\right] \right) \boldsymbol{t}_{2}({\mathbf{x}%
_{0}}).  \label{eqn:lbu0}
\end{eqnarray}
\endgroup
We notice that the Bochner Laplacian in (\ref{eqn:lbu0}) involves only the coefficients $\{a_{\alpha_1,\alpha_2}\}$ that satisfy ${\alpha_1+\alpha_2=2}$, which correspond to the monomials in (\ref{eqn:qhat}) of degree 2.

\subsubsection{Approximation of function and its derivatives in the local Monge coordinate\label{sec:comp}}
Let $\boldsymbol{u}=\tilde{u}^{i}\partial _{i}$ be a vector field defined locally on the stencil $S_{\mathbf{x}_0}$, where $\left\{
\partial _{1},\partial_2\right\}$ is the local basis defined in (\ref{eqn:pak}). Our goal here
is to approximate the Bochner Laplacian $\Delta _{B}$ acting on the vector
field $\boldsymbol{u}$ at ${\mathbf{x}_{0}}$ by a linear combination of
the vector values on $S_{\mathbf{x}_0}$, i.e., $\{\boldsymbol{\tilde{u}}|_{{\mathbf{x}_{0,k}%
}}:=(\tilde{u}^{1}({\mathbf{x}_{0,k}}),\tilde{u}^{2}({\mathbf{x}_{0,k}}%
))^{\top }\in \mathbb{R}^{2\times 1}\}_{k=1}^{K}$. Thus, we are seeking the
weights $\{\boldsymbol{\tilde{w}}^{k}\in \mathbb{R}^{2\times 2}\}_{k=1}^{K}$
such that,
\begin{equation}
\Delta _{B}\boldsymbol{u}({\mathbf{x}_{0}})\approx \sum_{k=1}^{K}\boldsymbol{\tilde{w}}%
^{k}\boldsymbol{\tilde{u}}|_{{\mathbf{x}_{0,k}}}=(\boldsymbol{\tilde{w}}^{1},\ldots ,%
\boldsymbol{\tilde{w}}^{K})\mathbf{\tilde{u}}_{{\mathbf{x}_{0}}}\in \mathbb{R%
}^{2\times 1},  \label{eqn:Bochin}
\end{equation}%
where $\mathbf{\tilde{u}}_{{\mathbf{x}_{0}}}:=(\tilde{u}^{1}({\mathbf{x}%
_{0,1}}),\tilde{u}^{2}({\mathbf{x}_{0,1}}),\ldots ,\tilde{u}^{1}({\mathbf{x}%
_{0,K}}),\tilde{u}^{2}({\mathbf{x}_{0,K}}))^{\top }\in \mathbb{R}^{2K\times
1}$.

The GMLS approximation of the vector field  $\boldsymbol{u}=%
\tilde{u}^{i}\partial _{i}$ gives the GMLS solution to each component, i.e.,
\begin{equation}
\tilde{u}^{1}=\sum_{0\leq \alpha _{1}+\alpha _{2}\leq l}\beta _{\alpha
_{1},\alpha _{2}}\theta _{1}^{\alpha _{1}}\theta _{2}^{\alpha _{2}},\text{ \
\ }\tilde{u}^{2}=\sum_{0\leq \alpha _{1}+\alpha _{2}\leq l}\gamma _{\alpha
_{1},\alpha _{2}}\theta _{1}^{\alpha _{1}}\theta _{2}^{\alpha _{2}}.
\label{eqn:ui}
\end{equation}
%
Then one can calculate the functions and
their derivatives at $\mathbf{x}_{0}$ as%
\begin{eqnarray*}
\tilde{u}^{1}|_{\mathbf{x}_{0}} &=&\beta _{0,0},\text{ \ }\partial _{11}%
\tilde{u}^{1}|_{\mathbf{x}_{0}}=2\beta _{2,0},\text{ \ }\partial _{22}\tilde{%
u}^{1}|_{\mathbf{x}_{0}}=2\beta _{0,2}, \\
\tilde{u}^{2}|_{\mathbf{x}_{0}} &=&\gamma _{0,0},\text{ \ }\partial _{11}%
\tilde{u}^{2}|_{\mathbf{x}_{0}}=2\gamma _{2,0},\text{ \ }\partial _{22}%
\tilde{u}^{2}|_{\mathbf{x}_{0}}=2\gamma _{0,2}.
\end{eqnarray*}%
Substituting all above into $\Delta _{B}\boldsymbol{u}(\mathbf{x}_{0})$ in (\ref{eqn:lbu0}%
), one arrives at%
\begin{eqnarray}
&&\Delta _{B}\boldsymbol{u}(\mathbf{x}_{0})=(\Delta _{B}\boldsymbol{u}(\mathbf{x}_{0}))_{1}%
\boldsymbol{t}_{1}({\mathbf{x}_{0}})+(\Delta _{B}\boldsymbol{u}(\mathbf{x}_{0}))_{2}%
\boldsymbol{t}_{2}({\mathbf{x}_{0}})  \label{eqn:uab} \\
&=&\left[ 2\beta _{2,0}+2\beta _{0,2}+\beta _{0,0}\left(
4(a_{2,0})^{2}+(a_{1,1})^{2}\right) +\gamma _{0,0}\left(
2a_{2,0}a_{1,1}+2a_{0,2}a_{1,1}\right) \right] \boldsymbol{t}_{1}({\mathbf{x}%
_{0}})  \notag \\
&+&\left[ \beta _{0,0}\left( 2a_{2,0}a_{1,1}+2a_{0,2}a_{1,1}\right) +2\gamma
_{2,0}+2\gamma _{0,2}+\gamma _{0,0}\left(
4(a_{0,2})^{2}+(a_{1,1})^{2}\right) \right] \boldsymbol{t}_{2}({\mathbf{x}%
_{0}})  \notag \\
&:= &\left( \boldsymbol{A}_{11}\boldsymbol{\beta }+\boldsymbol{A}_{12}%
\boldsymbol{\gamma }\right) \boldsymbol{t}_{1}({\mathbf{x}_{0}})+\left(
\boldsymbol{A}_{21}\boldsymbol{\beta }+\boldsymbol{A}_{22}\boldsymbol{\gamma
}\right) \boldsymbol{t}_{2}({\mathbf{x}_{0}}),  \notag
\end{eqnarray}%
where the first component of $\Delta _{B}\boldsymbol{u}(\mathbf{x}_{0})$ is $(\Delta
_{B}\boldsymbol{u}(\mathbf{x}_{0}))_{1}:=\boldsymbol{A}_{11}\boldsymbol{\beta }+%
\boldsymbol{A}_{12}\boldsymbol{\gamma }$ and the second component is $%
(\Delta _{B}\boldsymbol{u}(\mathbf{x}_{0}))_{2}:=\boldsymbol{A}_{21}\boldsymbol{\beta }+%
\boldsymbol{A}_{22}\boldsymbol{\gamma }$. Here $\boldsymbol{\beta }=\left(
\beta _{\alpha _{1}(j),\alpha _{2}(j)}\right) _{j=1}^{m}\in \mathbb{R}%
^{m\times 1}$\ and $\boldsymbol{\gamma }=\left( \gamma _{\alpha
_{1}(j),\alpha _{2}(j)}\right) _{j=1}^{m}\in \mathbb{R}^{m\times 1}$ are
both $m$ by $1$ column vectors and $\{\boldsymbol{A}_{ij}\in \mathbb{R}%
^{1\times m}\}_{i,j=1}^{2}$ are all $1$ by $m$ row vectors. Specifically, $%
\boldsymbol{A}_{11}$ has three nonzero entries $%
4(a_{2,0})^{2}+(a_{1,1})^{2},2,2$\ correponding to $\beta _{0,0},\beta
_{2,0},\beta _{0,2}$, respectively, and $\boldsymbol{A}_{12}$ has only one
nonzero entry $2a_{2,0}a_{1,1}+2a_{0,2}a_{1,1}$ correponding to $\gamma
_{0,0}$. $\boldsymbol{A}_{21}$ and $\boldsymbol{A}_{22}$\ are defined
similarly.

\comment{
\begin{remark}
The highest degree of the surface in (\ref{%
eqn:qhat}) should be taken at
least one order
higher
than
the
degree
of
the
component functions in (\ref{eqn:ui})
in
order
to
achieve
the
desired
accuracy. This is so
because
the
Bochner
Laplacian
in (\ref{eqn:ux0b})
involves
with
the
derivatives
of
Christoffel
symbols
corresponding
to
the
third
derivatives of the
surface $\hat{q}$
whereas
the
Bochner
involves
with
only the second
derivatives
of
component
functions $\tilde{u}^{1}$
and $\tilde{u}^{2}$.
\end{remark}
}




We now look for the weights $\{\boldsymbol{\tilde{w}}^{k}\}_{k=1}^{K}$
satisfying (\ref{eqn:Bochin}) using the above intrinsic formulation (\ref%
{eqn:uab}). Notice from (\ref{eqn:ui})\ that the coefficients $%
\boldsymbol{\beta }$ and $\boldsymbol{\gamma }$\ can be solved from normal
equations by regressions $\boldsymbol{\beta }=\boldsymbol{\Phi }^{\dag }%
\mathbf{\tilde{u}}_{{\mathbf{x}_{0}}}^{1}$ and $\boldsymbol{\gamma }=%
\boldsymbol{\Phi }^{\dag }\mathbf{\tilde{u}}_{{\mathbf{x}_{0}}}^{2}$,
respectively, where $\boldsymbol{\Phi }^{\dag }\in \mathbb{R}^{m\times K}$
is given in (\ref{eq:phinv}) and $\mathbf{\tilde{u}}_{{\mathbf{x}_{0}}}^{j}=(%
\tilde{u}^{j}({\mathbf{x}_{0,1}}),...,\tilde{u}^{j}(\mathbf{x}_{0,K}))^{\top
}\in \mathbb{R}^{K\times 1}$ $(j=1,2)$. Then one can write $(\Delta _{B}\boldsymbol{u}(%
\mathbf{x}_{0}))_{1}$ as:
\begin{equation}
(\Delta _{B}\boldsymbol{u}(\mathbf{x}_{0}))_{1}=\boldsymbol{A}_{11}\boldsymbol{\Phi }%
^{\dag }\mathbf{\tilde{u}}_{{\mathbf{x}_{0}}}^{1}+\boldsymbol{A}_{12}%
\boldsymbol{\Phi }^{\dag }\mathbf{\tilde{u}}_{{\mathbf{x}_{0}}}^{2}:=%
\boldsymbol{\tilde{w}}_{11}^{\bullet }\mathbf{\tilde{u}}_{{\mathbf{x}_{0}}%
}^{1}+\boldsymbol{\tilde{w}}_{12}^{\bullet }\mathbf{\tilde{u}}_{{\mathbf{x}_{0}%
}}^{2},  \label{eqn:lu1}
\end{equation}%
where $\boldsymbol{\tilde{w}}_{11}^{\bullet }=(\tilde{w}_{11}^{1},\tilde{w}%
_{11}^{2},\ldots ,\tilde{w}_{11}^{K})\in \mathbb{R}^{1\times K}$ and $%
\boldsymbol{\tilde{w}}_{12}^{\bullet }=(\tilde{w}_{12}^{1},\tilde{w}%
_{12}^{2},\ldots ,\tilde{w}_{12}^{K})\in \mathbb{R}^{1\times K}$ correspond
to the first row of $(\boldsymbol{\tilde{w}}^{1},\ldots ,\boldsymbol{\tilde{w%
}}^{K})\in \mathbb{R}^{2\times 2K}$.  One can follow the same
procedure to obtain the second row of $(\boldsymbol{\tilde{w}}^{1},\ldots ,%
\boldsymbol{\tilde{w}}^{K})$. Last one can stack all the components $\{%
\tilde{w}_{ij}^{k}\}_{i,j=1,2}^{k=1,\ldots ,K}$ to form the weights $(%
\boldsymbol{\tilde{w}}^{1},\ldots ,\boldsymbol{\tilde{w}}^{K})\in \mathbb{R}%
^{2\times 2K}$ as desired in (\ref{eqn:Bochin}).

%
%

%
%

\subsubsection{Matrix assembly of the Bochner Laplacian}\label{sec:MaBL}
In this section, we discuss how to calculate the weights $\{\boldsymbol{w}^{k}\}$ in (\ref{eqn:LB}) from the weights $\{\boldsymbol{\tilde{w}}^{k}\}$ in (\ref{eqn:Bochin}) and assemble all $\{\boldsymbol{w}^{k}\}$  into the appropriate rows and columns of a sparse $%
2N\times 2N$ Laplacian matrix $\mathbf{L}_{B}$.
So far, we have obtained the components$\ \{\boldsymbol{\tilde{u}}|_{{\mathbf{x}_{0,k}}}\in
\mathbb{R}^{2\times 1}\}_{k=1}^{K}$ and the weights $\{\boldsymbol{\tilde{w}}%
^{k}\}_{k=1}^{K}$ to approximate $\Delta _{B}\boldsymbol{u}({\mathbf{x}_{0}})$ (see equation (\ref%
{eqn:Bochin}))\ in the local bases $\left\{ (\partial _{1}|_{\mathbf{x}%
_{0,k}},\partial _{2}|_{\mathbf{x}_{0,k}})\right\}
_{k=1}^{K}$ of each stencil $S_{{\mathbf{x}_{0}}}$. In order to construct a
sparse $2N\times 2N$\ Bochner Laplacian matrix, we need to transform the
weights from the local bases $\left\{ \partial _{1}|_{\mathbf{x}%
_{0,k}},\partial _{2}|_{\mathbf{x}_{0,k}}\right\} _{k=1}^{K}$ in each $S_{{%
\mathbf{x}_{0}}}$ to the unified global basis $\{\boldsymbol{t}_{1}(\mathbf{x%
}_{0,k}),\boldsymbol{t}_{2}(\mathbf{x}_{0,k})\}_{k=1}^{K}$.

\comment{Let $u=u^{j}%
\boldsymbol{t}_{j}$ be a vector field \comment{with the global basis $\{\boldsymbol{t}
_{1},\boldsymbol{t}_2\}$\ and components $\{u^{1},u^2\}$} defined on the
nodes $\{\mathbf{x}_{i}\}_{i=1}^{N}$. We intend to find the weights $\{%
\boldsymbol{w}^{k}\in \mathbb{R}^{2\times 2}\}_{k=1}^{K}$ such that,
\begin{equation}
\sum_{k=1}^{K}\boldsymbol{w}^{k}u({\mathbf{x}_{0,k}})=(\boldsymbol{w}%
^{1},\ldots ,\boldsymbol{w}^{K})\mathbf{u}_{{\mathbf{x}_{0}}}\approx \Delta
_{B}u({\mathbf{x}_{0}})\approx \sum_{k=1}^{K}\boldsymbol{\tilde{w}}^{k}%
\tilde{u}({\mathbf{x}_{0,k}}),  \label{eqn:LBw}
\end{equation}%
where $\{u({\mathbf{x}_{0,k}})=(u^{1}({\mathbf{x}_{0,k}}),u^{2}({\mathbf{x}%
_{0,k}}))^{\top }\in \mathbb{R}^{2\times 1}\}_{k=1}^{K}$\ and$\ \mathbf{u}_{{%
\mathbf{x}_{0}}}=(u^{1}({\mathbf{x}_{0,1}}),u^{2}({\mathbf{x}_{0,1}}),\ldots
,u^{1}({\mathbf{x}_{0,K}}),u^{2}({\mathbf{x}_{0,K}}))^{\top }\in \mathbb{R}%
^{2K\times 1}$. }
Let $\{\boldsymbol{w}%
^{k}\}_{k=1}^{K}$ be the weights as we seek in (\ref{eqn:LB}), then
\begin{equation}
\sum_{k=1}^{K}\boldsymbol{w}^{k}\boldsymbol{u}|_{{\mathbf{x}_{0,k}}}=(\boldsymbol{w}%
^{1},\ldots ,\boldsymbol{w}^{K})\mathbf{u}_{{\mathbf{x}_{0}}}\approx \Delta
_{B}\boldsymbol{u}({\mathbf{x}_{0}})\approx \sum_{k=1}^{K}\boldsymbol{\tilde{w}}^{k}%
\boldsymbol{\tilde{u}}|_{{\mathbf{x}_{0,k}}}.  \label{eqn:LBw}
\end{equation}
In particular, we let
\begin{equation*}
\boldsymbol{w}^{k}\boldsymbol{u}|_{{\mathbf{x}_{0,k}}}\approx \boldsymbol{\tilde{w}}^{k}%
\boldsymbol{\tilde{u}}|_{{\mathbf{x}_{0,k}}},\text{ \ \ }\forall k=1,\ldots ,K.
\end{equation*}%
We note that $\boldsymbol{u}=u^{j}\boldsymbol{t}_{j}\approx \tilde{u}^{i}\partial
_{i} $ in each $S_{{\mathbf{x}_{0}}}$, that is,%
\begin{equation*}
\left[
\begin{matrix}
\boldsymbol{t}_{1}({\mathbf{x}_{0,k}}) & \boldsymbol{t}_{2}({\mathbf{x}_{0,k}%
})%
\end{matrix}%
\right] \boldsymbol{u}|_{{\mathbf{x}_{0,k}}}\approx \left[
\begin{matrix}
\partial _{1}|_{{\mathbf{x}_{0,k}}} & \partial _{2}|_{{\mathbf{x}_{0,k}}}%
\end{matrix}%
\right] \boldsymbol{\tilde{u}}_{{\mathbf{x}_{0,k}}}\in \mathbb{R}^{3\times 1},
\end{equation*}%
where $\mathrm{span}\left\{ \partial _{1},\partial _{2}\right\} \approx
\mathrm{span}\{\boldsymbol{t}_{1},\boldsymbol{t}_{2}\}$ in the neighborhood
of ${\mathbf{x}_{0}}$ in which the small deviation is induced by the local
approximation of Monge parametrization. In general, the bases $\left\{
\partial _{1},\partial _{2}\right\} $ in (\ref{eqn:pak}) are not orthogonal
at the neighboring $\{\mathbf{x}_{0,k}\}_{k=2}^{K}$ due to nonflat geometry
in the Taylor expansion, $g_{ij}({\mathbf{x}_{0}})=\delta _{ij}+\sum_{k,l}%
\frac{1}{3}R_{iklj}({\mathbf{x}_{0}})\theta _{k}\theta _{l}+O(\Vert \theta
\Vert ^{3})$. Then we can solve the least-squares problem to obtain the
weight $\boldsymbol{w}^{k}$ for $\forall k=1,\ldots ,K$:%
\begin{equation*}
\boldsymbol{w}^{k}\approx \boldsymbol{\tilde{w}}^{k}\left( \left[
\begin{matrix}
\partial _{1}^{\top }|_{{\mathbf{x}_{0,k}}} \\
\partial _{2}^{\top }|_{{\mathbf{x}_{0,k}}}%
\end{matrix}%
\right] \left[
\begin{matrix}
\partial _{1}|_{{\mathbf{x}_{0,k}}} & \partial _{2}|_{{\mathbf{x}_{0,k}}}%
\end{matrix}
\right] \right) ^{-1}\left[
\begin{matrix}
\partial _{1}^{\top }|_{{\mathbf{x}_{0,k}}} \\
\partial _{2}^{\top }|_{{\mathbf{x}_{0,k}}}%
\end{matrix}%
\right] \left[
\begin{matrix}
\boldsymbol{t}_{1}({\mathbf{x}_{0,k}}) & \boldsymbol{t}_{2}({\mathbf{\ x}%
_{0,k}})%
\end{matrix}%
\right] .
\end{equation*}%
Then, we repeat the above procedure to obtain the weights for each point in $%
\{\mathbf{x}_{i}\}_{i=1}^{N}$. Last, we assemble all these weights $\{%
\boldsymbol{w}^{k}\}$ into the appropriate rows and columns of a sparse $%
2N\times 2N$\ Laplacian matrix $\mathbf{L}_{B}$ as a pointwise
approximation to the Bochner Laplacian $\Delta _{B}$ on the point cloud $\{%
\mathbf{x}_{i}\}_{i=1}^{N}$.

\subsection{Intrinsic formulation in higher dimensional spaces} \label{3p3highint}
We now extend our approach to the general
$d$-dimensional manifolds embedded in $n$-dimensional ambient spaces ($d\ll n$). In numerical implementation, we focus on small $d$ ranging around
$1\sim 3$  and relatively large $n$ on the order of $O(10)$. In Section \ref{sec:intr_high_d}, we first derive
the intrinsic formulation for the Bochner Laplacian in the local Monge
parametrization. However, a direct implementation of the formulation in Section \ref{sec:intr_high_d} results in a computationally expensive algorithm with a complexity of $O(Nn^{3})$ (due to the full QR factorization to  generate normal directions). To reduce the
complexity for manifold approximation, we introduce a practical algorithm in Section \ref{sec:comp_high_d} that reduces the complexity
to $O(Nn)$.
In Section \ref{sec:intr_alg}, we present the complete algorithm along with its complexity analysis.

\comment{
Let us consider the intrinsic formulation for the Bochner Laplacian in
general cases of $d$ dimensional manifolds embedded in $n$ dimensional
ambient spaces. In numerical implementation, we again only consider small $d$ ranging around
$1\sim 4$ {\color{red}(due to the curse of dimensionality)} and relatively large $n$ ranging around $O(10)$. In Section \ref{sec:intr_high_d}, we derive
the intrinsic formulation for the Bochner Laplacian in the local Monge
parametrization. However, directly following the formulation, our algorithm
is computational expensive with the complexity of $O(Nn^{3})$. To reduce the
complexity, we then present an practical algorithm which has the complexity
of $O(Nn)$ in Section \ref{sec:comp_high_d},
}

\subsubsection{Intrinsic formulation in the Monge parametrization}

\label{sec:intr_high_d}
For the ease of discussion, we again use $\left\{ \boldsymbol{t}_{i}({%
	\mathbf{x}_{j}})\in \mathbb{R}^{n}\right\} _{i=1,\ldots ,d}^{j=1,\ldots ,N}$
to denote the global tangent vector basis on the point cloud $\left\{ {%
	\mathbf{x}_{j}}\right\} _{j=1}^{N}$\ for both known and unknown manifolds.
Then
one can define a local coordinate chart for the manifold near the base point,

\begin{equation}
	{\mathbf{x}}(\theta _{1},\ldots ,\theta _{d})={\mathbf{x}_{0}}%
	+\sum_{i=1}^{d}\theta _{i}\boldsymbol{t}_{i}({\mathbf{x}_{0}}%
	)+\sum_{s=1}^{n-d}q^{s}\left( \theta _{1},\ldots ,\theta _{d}\right)
	\boldsymbol{n}_{s}({\mathbf{x}_{0}}),  \label{eqn:xtq}
\end{equation}%
where $\theta _{i}\left( {\mathbf{x}}\right) =\boldsymbol{t}_{i}({\mathbf{x}%
	_{0}})\cdot (\mathbf{x}-{\mathbf{x}_{0}})$ and $q^{s}\left( {\mathbf{x}}%
\right) =\boldsymbol{n}_{s}({\mathbf{x}_{0}})\cdot (\mathbf{x}-{\mathbf{x}_{0}}).$ The tangent vector bases in the Monge coordinate system are
$	\partial _{k}=\frac{\partial }{\partial \theta _{k}}=\boldsymbol{t}_{k}({%
		\mathbf{x}_{0}})+\sum_{j=1}^{n-d}\frac{\partial q^{j}}{\partial \theta _{k}}%
	\boldsymbol{n}_{j}({\mathbf{x}_{0}}),\text{ \ for }k=1,\ldots ,d\text{.}$
Each component $g_{ij}$\ of the local Riemannian metric $\boldsymbol{g}%
=(g_{ij})_{i,j=1}^{d}$ becomes
\begin{equation}
	g_{ij}=\delta _{ij}+\sum_{s=1}^{n-d}\frac{\partial q^{s}}{\partial \theta
		_{i}}\frac{\partial q^{s}}{\partial \theta _{j}}.  \label{eqn:gij3}
\end{equation}

Let the vector field $\boldsymbol{u}=\tilde{u}^{i}\partial _{i}$ in the
local Monge coordinate system. Then the Bochner Laplacian takes the form
\begin{equation*}
	\Delta _{B}\boldsymbol{u}(\mathbf{x}_{0})=g^{kj}\tilde{u}_{;kj}^{i}\partial
	_{i}|_{\mathbf{x}_{0}}=\sum_{k}\tilde{u}_{;kk}^{i}\partial _{i}|_{\mathbf{x}%
		_{0}}=\sum_{k}\left( \partial _{kk}\tilde{u}^{i}+\tilde{u}^{l}\frac{\partial
		\Gamma _{kl}^{i}}{\partial \theta _{k}}\right) \partial _{i}|_{\mathbf{x}%
		_{0}}.
\end{equation*}%
In the second equality, we use the fact that the metric $\boldsymbol{g}$ is the identity matrix at $\mathbf{x}_0$.  By using the metric (\ref{eqn:gij3}), the derivative of the Christoffel
symbol at the base point can be calculated as,%
\begin{small}
\begin{eqnarray}
	&&\quad  \frac{\partial \Gamma _{kl}^{i}}{\partial \theta _{r}}\Big\vert_{\mathbf{x}_{0}} =%
	\frac{1}{2}g^{im}\left( \frac{\partial g_{mk}}{\partial \theta _{r}\partial
		\theta _{l}}+\frac{\partial g_{ml}}{\partial \theta _{r}\partial \theta _{k}}%
	-\frac{\partial g_{kl}}{\partial \theta _{r}\partial \theta _{m}}\right) \Big\vert_{%
		\mathbf{x}_{0}} = \sum_{s=1}^{n-d}\partial _{ir}q^{s}\partial _{kl}q^{s}|_{\mathbf{x}_{0}}, \label{eqn:derT}  
\end{eqnarray}%
\end{small}
where $\partial _{kl}q^{s}=\frac{\partial ^{2}q^{s}}{\partial \theta
_{k}\partial \theta _{l}}$. Then the Bochner Laplacian can be written as%
\begin{equation}
	\Delta _{B}\boldsymbol{u}(\mathbf{x}_{0})=\sum_{k}\left( \partial _{kk}%
	\tilde{u}^{i}+\tilde{u}^{l}\sum_{s=1}^{n-d}\partial _{ik}q^{s}\partial
	_{kl}q^{s}\right) \partial _{i}|_{\mathbf{x}_{0}}.  \label{eqn:Bochdn}
\end{equation}%
The results of the Hodge Laplacian and another vector Laplacian can be found in  Section~\ref{app:B}.

\subsubsection{Numerical implementation}
\label{sec:comp_high_d}
Note that the formulation (\ref{eqn:Bochdn}) requires computing the $n-d$\ normal vectors $\left\{
\boldsymbol{n}_{i}({\mathbf{x}_{0}})\right\} _{i=1}^{n-d}$ at the base
point which are obtained using the  QR decomposition algorithm with a computational cost of $O(n^{3})$. To avoid explicitly computing these normal vectors, we instead apply the
GMLS to each component of ${\mathbf{x}}=(x^{1},\ldots ,x^{n})$ in the
stencil $S_{{\mathbf{x}_{0}}}=\left\{ {\mathbf{x}}_{0,k}\right\}
_{k=1}^{K}= \{ x_{0,k}^{s}\} _{k=1,\ldots ,K}^{s=1,\ldots ,n}$ rather than using $q^s$ as defined in (\ref{eqn:xtq}). For $s=1,\ldots ,n$, the $s$th component of the Taylor expansion of $\mathbf{x}(\boldsymbol{\theta})$ with respect to $\boldsymbol{\theta}=(\theta_1,...,\theta_d)$ at $\boldsymbol{\theta}=\boldsymbol{0}$ is given by
\begin{equation}
	x^{s}=x_{0}^{s}+\sum_{i=1}^{d}t_{i}^{s}\theta _{i}+\sum_{2\leq
		\left\vert \boldsymbol{\alpha }\right\vert \leq l}b_{\alpha _{1},\ldots
		,\alpha _{d}}^{s}\theta _{1}^{\alpha _{1}}\cdots \theta _{d}^{\alpha
		_{d}}+O(\left\Vert \boldsymbol{\theta }\right\Vert ^{l+1}),
		\label{eqn:pqs1}
\end{equation}%
Here $t_i^s$ is the $s$th component of the tangent vectors $\boldsymbol{t}_{i}({\mathbf{x%
	}_{0}})$ and $\mathbf{x}(\mathbf{0})=(x_0^1,...,x_0^n)$. We then define $p^s$ as
\begin{equation}
p^{s}=x^s-x_{0}^{s}-\sum_{i=1}^{d}t_{i}^{s}\theta _{i}=\sum_{2\leq \left\vert \boldsymbol{\alpha }\right\vert
		\leq l}b_{\boldsymbol{\alpha }}^{s}\boldsymbol{\theta }^{\boldsymbol{\alpha }}.
		\label{eqn:xps}
\end{equation}
where $b_{\boldsymbol{\alpha }}^{s}$ will be determined by the GMLS approximation. For the terms where $|\boldsymbol{\alpha }|=2$ in (\ref{eqn:xps}), we can write them in the following quadratic form,%
\begin{equation}
	\sum_{\left\vert \boldsymbol{\alpha }\right\vert =2}b_{\boldsymbol{\alpha }%
	}^{s}\boldsymbol{\theta }^{\boldsymbol{\alpha }}=\left[
	\begin{array}{c}
		\theta _{1} \\
		\vdots  \\
		\theta _{d}%
	\end{array}%
	\right] ^{\top }\left[
	\begin{array}{ccc}
		c_{11}^{s} & \cdots  & c_{1d}^{s} \\
		\vdots  & \ddots  & \vdots  \\
		c_{d1}^{s} & \cdots  & c_{dd}^{s}%
	\end{array}%
	\right] \left[
	\begin{array}{c}
		\theta _{1} \\
		\vdots  \\
		\theta _{d}%
	\end{array}%
	\right] ,\text{ \ for }s=1,\ldots ,n,  \label{eqn:bthe}
\end{equation}%
where the coefficient matrix is symmetric with $c_{ii}^{s}=\{b_{%
	\boldsymbol{\alpha }}^{s}|\alpha _{i}=2$ and the others $\alpha _{j}=0$ for $%
j\neq i\}$ and $c_{ik}^{s}=$ $\{\frac{1}{2}b_{\boldsymbol{\alpha }%
}^{s}|\alpha _{i}=\alpha _{k}=1$ for $i\neq k$ and the others $\alpha _{j}=0$
for $j\neq i$ and $j\neq k\}$.

On the other hand, let the $s$th normal component $q^s$ in the
Monge coordinate (\ref{eqn:xtq}) be
\begin{equation}
	q^{s}=\sum_{2\leq \left\vert \boldsymbol{\alpha }\right\vert \leq l}a_{%
		\boldsymbol{\alpha }}^{s}\boldsymbol{\theta }^{\boldsymbol{\alpha }%
	},\text{ \ for }%
	s=1,\ldots ,n-d,\label{eqn:qs_high}
\end{equation}%
which generalizes equation (\ref{eqn:qhat}) from  dimension $2$  to
dimension $d$. Then we have the following proposition with its proof provided in Section~\ref{sec:app1_3}.

\begin{proposition}\label{prop3p2}
	Let $\left\{ p^{i}\right\} _{i=1}^{n}$ and $\left\{ q^{j}\right\}
	_{j=1}^{n-d}$ be defined as in (\ref{eqn:xps}) and (\ref{eqn:qs_high}), respectively. Then for any integers $1\leq i,r,k,l\leq
	d,$
	\begin{equation}
		\sum_{s=1}^{n-d}\partial _{ir}q^{s}\partial _{kl}q^{s} = \sum_{s=1}^{n}\partial _{ir}p^{s}\partial
		_{kl}p^{s},
		\label{eqn:pqe}
	\end{equation}%
	where the derivative $\partial _{ir}=\frac{\partial ^{2}}{\partial \theta
		_{i}\partial \theta _{r}}$ is with respect to $\left\{ \theta _{i}\right\}
	_{i=1}^{d}$. Moreover, at the base point $\mathbf{x}_{0}$, it holds that%
	\begin{equation}
			\sum_{s=1}^{n-d}\partial _{ir}q^{s}\partial _{kl}q^{s}|_{\mathbf{x}%
			_{0}}=\sum_{s=1}^{n}4 c_{ir}^{s}c_{kl}^{s} ,
			\label{eqn:pqe2}
	\end{equation}
	where $c_{ir}^s$ is defined  in (\ref{eqn:bthe}).
\end{proposition}
By using the equality (\ref{eqn:pqe2}), the Bochner Laplacian in (\ref{eqn:Bochdn}) can be written as
\begin{equation}
	\Delta _{B}\boldsymbol{u}(\mathbf{x}_{0})=\sum_{k}\left( \partial _{kk}%
	\tilde{u}^{i}+\tilde{u}^{l}\sum_{s=1}^{n}4 c_{ir}^{s}c_{kl}^{s} \right) \partial _{i}|_{\mathbf{x}_{0}}.  \label{eqn:Bochds}
\end{equation}
We note that  the Bochner Laplacian involves only the coefficients of the GMLS approximation of $p^s$, corresponding to ${|\boldsymbol{\alpha}|=2}$ in (\ref{eqn:bthe}).

\subsubsection{Algorithm and complexity analysis}
\label{sec:intr_alg}
The complete procedure of the intrinsic GMLS approximation of the Bochner Laplacian (in the case when $d\ll n$) and its complexity analysis are outlined in Algorithm~\ref{algo:intrin-Boch}. Notably, one can follow the similar procedure to approximate other vector Laplacians acting on vector fields  using the intrinsic method (See Section~\ref{app:B}). 

\begin{algorithm}[ht]
	\caption{Intrinsic GMLS of the Bochner Laplacian when $d\ll n$}

	\begin{algorithmic}[1]
		\STATE {\bf Input:} A point cloud $\{\mathbf{x}_{i}\}_{i=1}^{N}\subset M$, bases of (estimated) tangent vectors at each node $\left\{ \boldsymbol{t}_{1}(\mathbf{x}_{i}),\ldots,\boldsymbol{t}_{d}(%
		\mathbf{x}_{i})\right\} _{i=1}^{N}$, the degree $l$ of  polynomials, and a parameter $K(>m)$ nearest neighbors where $%
		m=\left(
		\begin{array}{c}
			l+d \\
			d%
		\end{array}%
		\right) $ is the number of monomial basis functions $\left\{ \boldsymbol{%
			\theta }^{\boldsymbol{\alpha }}|0\leq \left\vert \boldsymbol{\alpha }%
		\right\vert \leq l\right\} $.
		\STATE Set $\mathbf{L}_B$ to be a sparse $dN\times dN$ matrix with $d^2NK$ nonzeros.
		\FOR{$i\in \{1,2,...,N\}$}
		\STATE Find the $K$ nearest neighbors of the point $\mathbf{x}_i$ in the stencil $S_{\mathbf{x}_i}=\{\mathbf{x}_{i,k}\}_{k=1}^K$.
		\STATE Find the GMLS
		approximation of $p^{s}=x^{s}-x_{0}^{s}-\sum_{i=1}^{d}t_{i}^{s}\theta _{i}$ as in (\ref{eqn:xps})
		using the optimal solution $\hat{p}^{s}=\sum_{2\leq \left\vert \boldsymbol{%
				\alpha }\right\vert \leq l}b_{\boldsymbol{\alpha }}^{s}\boldsymbol{\theta }^{%
			\boldsymbol{\alpha }}$ to the least-squares problem with the complexity of $O(NnKm^{2}+Nnm^{3})$.
		\STATE Construct the local tangent-vector bases $\{\partial _{j}|_{\mathbf{x}%
			_{i,k}}\}_{j=1,\ldots ,d}$ over $S_{\mathbf{x}_i}$ with the complexity $O(NnKm^{2}d)$ based on taking the derivative of (\ref{eqn:pqs1}).
		\STATE  Compute the coefficients of the Bochner Laplacian $\{\boldsymbol{\tilde{w}}^{k}\}_{k=1}^{K}$ in the local Monge coordinate  based on (\ref{eqn:Bochds}) with the
		complexity $O(Nnd^{3})$.
		\STATE Calculate the weights $\{\boldsymbol{w}^{k}\}_{k=1}^{K}$ with respect to the global tangent-vector bases with
		the complexity $O(NKd^{2}n+NKd^{4})$ by the formula,%
		\begin{small}
		\begin{equation*}
			\boldsymbol{w}^{k}=\boldsymbol{\tilde{w}}^{k}\left( \left[ %
			\setlength{\arraycolsep}{0.3pt}%
			\begin{array}{c}
				\partial _{1}^{\top }|_{{\mathbf{x}_{0,k}}} \\
				\vdots \\
				\partial _{d}^{\top }|_{{\mathbf{x}_{0,k}}}%
			\end{array}%
			\right] \left[ \setlength{\arraycolsep}{0.3pt}
			\begin{array}{ccc}
				\partial _{1}|_{{\mathbf{x}_{0,k}}} & \cdots & \partial _{d}|_{{\mathbf{x}%
						_{0,k}}}%
			\end{array}%
			\right] \right) ^{-1}\left( \left[ \setlength{\arraycolsep}{0.3pt}
			\begin{array}{c}
				\partial _{1}^{\top }|_{{\mathbf{x}_{0,k}}} \\
				\vdots \\
				\partial _{d}^{\top }|_{{\mathbf{x}_{0,k}}}%
			\end{array}%
			\right] \left[ \setlength{\arraycolsep}{0.3pt}
			\begin{array}{ccc}
				\boldsymbol{t}_{1}({\mathbf{x}_{0,k}}) & \cdots & \boldsymbol{t}_{d}({%
					\mathbf{\ x}_{0,k}})%
			\end{array}%
			\right] \right) .
		\end{equation*}
		\end{small}
		\STATE Arrange the weights $\{\boldsymbol{w}^{k}\}_{k=1}^{K}$ into corresponding rows and columns of $\mathbf{L}_B$.
		\ENDFOR
		\STATE {\bf Output:} The approximate operator matrix $\mathbf{L}_B$.
	\end{algorithmic}

	\label{algo:intrin-Boch}
\end{algorithm}


\subsection{Approximation of covariant derivative}\label{sec:CD}
Note that the covariant derivative of
the vector field $\boldsymbol{u}=\tilde{u}^{i}\partial _{i}$ along itself could be written as
$$
\nabla_{\boldsymbol{u}}\boldsymbol{u}=\tilde{u}_{;k}^i\tilde{u}^k\partial_i=\tilde{u}^{k}\left( \partial _{k}\tilde{u}%
^{i}+\tilde{u}^{j}\Gamma _{kj}^{i}\right) \partial _{i}.
$$
By using the local Monge
coordinate system introduced in Section~\ref{sec:mgchart}, the approximation of $\nabla_{\boldsymbol{u}}\boldsymbol{u}$ on 2D manifolds could be derived directly as
\begin{eqnarray*}
\nabla _{\boldsymbol{u}}\boldsymbol{u}(\mathbf{x}_{0})&=&\tilde{u}^{k}\left( \partial _{k}\tilde{u}%
^{i}+\tilde{u}^{j}\Gamma _{kj}^{i}\right) \partial _{i}\big{|}_{\mathbf{x}%
_{0}} \\
&\approx&\left( \beta _{0,0}\beta _{1,0}+\gamma _{0,0}\beta _{0,1}\right)
\boldsymbol{t}_{1}({\mathbf{x}_{0}})+\left( \beta _{0,0}\gamma _{1,0}+\gamma
_{0,0}\gamma _{0,1}\right) \boldsymbol{t}_{2}({\mathbf{x}_{0}}),
\end{eqnarray*}%
where the coefficients $\beta _{\alpha _{1},\alpha _{2}}$ and $\gamma
_{\alpha _{1},\alpha _{2}}$ for $0\leq \alpha _{1}+\alpha _{2}\leq l$\ are
defined in (\ref{eqn:ui}). Here, we again use Proposition \ref{prop:norm} to simplify the calculations. We remark that this formulation naturally extends to higher dimensions.

\comment{
\subsection{Computational complexity}\label{sec:CC}
In the end of this section, we discuss the computational cost for approximating the Bochner Laplacian using the intrinsic formulation.

{\color{red} computational cost of bochner laplacian at a base point include the coefficients and regression matrix. }

Suppose that the coefficients $\boldsymbol{a}=\left\{ a_{%
\boldsymbol{\alpha }}\right\} _{2\leq |\boldsymbol{\alpha }|\leq l}$\ and the regression matrix $\left( \boldsymbol{\Phi }%
^{\top }\boldsymbol{W}\boldsymbol{\Phi }\right) ^{-1}\boldsymbol{\Phi }%
^{\top }\boldsymbol{W}$ have been pre-computed. Then we can see from (\ref{eqn:lu1}) that the computational cost for calculating all the weights $(\boldsymbol{\tilde{w}}^{1},\ldots ,\boldsymbol{\
\tilde{w}}^{K})$ of the Bochner Laplacian is $O(d^{2}mK)=O(d^{2}\left(\begin{matrix}l+d\\ d\end{matrix}\right)K)$  where $d=2$\ is the intrinsic dimension of the surface.
{\color{red} computational cost of bochner laplacian over point cloud (matrix assembly)} }

\section{GMLS using extrinsic formulation}
\label{sec:extf}

In this section, we approximate operators on a point cloud using the extrinsic formulation. Section~\ref{sec:revext} reviews the representations of the
differential operators on 2D surfaces as tangential derivatives \cite%
{fuselier2013high,Natasha2015Solving,shankar2015radial,lehto2017radial,harlim2023radial}%
, which could be formulated as the projection of the appropriate derivatives
in the ambient space. In Section~\ref{sec:revdisH}, we approximate the
extrinsic formulation of the gradient and the Bochner Laplacian of vector
fields using the GMLS approach \cite%
{liang2013solving,flyer2016role,gross2020meshfree,jones2023generalized,jiang2024generalized}%
.  The vector fields here are represented using the bases from ambient space  so that the resulting discrete Bochner Laplacian is of size $%
3N\times 3N$ where $N$ is the number of data points and $3$ corresponds to the ambient dimension. In Section~\ref{sec:dimred}, we introduce a coordinate transformation technique to convert the $3N\times 3N$ matrix to a $2N\times 2N$ matrix where $2$ corresponds to the intrinsic dimension of the
surface. In Section~\ref{hihext}, we directly present the generalized results of the Bochner Laplacian for $d$-dimensional manifolds in $\mathbb{R}^n$. In Section~\ref{sec:acd}, we derive the approximation formula for covariant derivatives, specifically for 2D manifolds in $\mathbb{R}^3$.

%
%
%
%
%
%
%


\subsection{Review of extrinsic formulation of differential operators on
manifolds}

\label{sec:revext}

We first review some basic concepts and define the conventional notations in extrinsic
differential geometry for 2D manifolds in $\mathbb{R}^3$. For any point $\mathbf{x}\in M$, the local
parametrization $\iota :O\subseteq \mathbb{R}^{2}\longrightarrow M\subseteq
\mathbb{R}^{3}$, is defined through the following map, $\boldsymbol{\Xi |}%
_{\iota ^{-1}(\mathbf{x})}\boldsymbol{\ \longmapsto }\mathbf{X|}_{\mathbf{x}%
} $. Here, $O$ denotes a open domain that contains the point $\iota ^{-1}(%
\mathbf{x})$, which we denoted as $\boldsymbol{\Xi }_{\iota ^{-1}(\mathbf{x}%
)}$ in the canonical coordinates $\left( \frac{\partial }{\partial \xi ^{1}},%
\frac{\partial }{\partial \xi ^{2}}\right) \Big\vert_{\iota ^{-1}(\mathbf{x}%
)}$ and $\mathbf{X|}_{\mathbf{x}}$ is the embedded point represented in the
ambient coordinates $\left( \frac{\partial }{\partial X^{1}},\frac{\partial
}{\partial X^{2}},\frac{\partial }{\partial X^{3}}\right) \Big\vert_{\mathbf{%
x}}$. The local coordinate $\left( \xi ^{1},\xi ^{2}\right) $ is chosen
arbitrarily which can also be the same as the local Monge coordinate $(\theta
_{1},\theta _{2})$ in Section~\ref{sec:mgchart}. The basic idea for the extrinsic formulation is to
rewrite the surface derivatives in local $\left( \xi ^{1},\xi ^{2}\right) $
as the projection of the derivatives in the ambient space $\left(
X^{1},X^{2},X^{3}\right) $.

Let $T_{\mathbf{x}}M$ be the tangent space at point $\mathbf{x}\in M$ and
denote a set of orthonormal tangent-vector basis by $\left\{ \boldsymbol{t}%
_{1}(\mathbf{x}),\boldsymbol{t}_{2}(\mathbf{x})\right\} $. Then the
projection matrix $\mathbf{P}$ which projects vectors in $\mathbb{R}^{3}$ to
$T_{\mathbf{x}}M$ could be written as $\mathbf{P}=\sum_{i=1}^{2}\boldsymbol{t%
}_{i}\boldsymbol{t}_{i}^{\top }$ at any $\mathbf{x}\in M$ (see e.g. \cite%
{harlim2023radial} for more detailed properties about $\mathbf{P}$).
Subsequently, the surface gradient acting on a smooth function $%
f:M\rightarrow \mathbb{R}$ evaluated at $\mathbf{x}\in M$ in the Cartesian
coordinates can be given by \cite%
{fuselier2013high,Natasha2015Solving,shankar2015radial,lehto2017radial},
\begin{equation}
\mathrm{grad}_{g}f(\mathbf{x})=\mathbf{P}\overline{\mathrm{grad}}_{\mathbb{R}%
^{3}}F(\mathbf{x})=\Big(\sum_{j=1}^{2}\boldsymbol{t}_{j}\boldsymbol{t}%
_{j}^{\top }\Big)\overline{\mathrm{grad}}_{\mathbb{R}^{3}}F(\mathbf{x}),
\label{eq:gdF}
\end{equation}%
where $\overline{\mathrm{grad}}_{\mathbb{R}^{3}}=[\partial _{X^{1}},\partial
_{X^{2}},\partial _{X^{3}}]^{\top }$ is the usual gradient operator defined
in $\mathbb{R}^{3}$ and $\mathrm{grad}_{g}f:=g^{ij}\frac{\partial f}{%
\partial \xi ^{i}}\frac{\partial }{\partial \xi ^{j}}$\ is defined w.r.t.
the Riemannian metric $\boldsymbol{g}$ induced by $M$ from $\mathbb{R}^{3}$. Here, $F$ is
a smooth extension of $f$ to an open subset in $\mathbb{R}^{3}$\ such that $%
F|_{M}=f$. Note that the identity (\ref{eq:gdF})\ holds true since for any
smooth function $f\in C^{\infty }(M)$, such an extension $F$\ exists and the
identity does not depend on the chosen extension \cite{lee2018introduction}.
In the following, we will abuse the notation to use $f$ as both a function
defined on a surface and its extension to an open set in $\mathbb{R}^{3}$. Let
$\left\{ \mathbf{e}_{1},\mathbf{e}_{2},\mathbf{e}_{3}\right\} $ be the
standard orthonormal basis corresponding to the Cartesian coordinate $%
\left\{ X^{1},X^{2},X^{3}\right\} $ in $\mathbb{R}^{3}$. Then we can rewrite
(\ref{eq:gdF}) in component form as
\begin{equation}
\mathrm{grad}_{g}f(\mathbf{x})=\left[
\begin{matrix}
\left( \mathbf{e}_{1}\cdot \mathbf{P}\right) \cdot \overline{\mathrm{grad}}_{%
\mathbb{R}^{3}}f(\mathbf{x}) \\
\left( \mathbf{e}_{2}\cdot \mathbf{P}\right) \cdot \overline{\mathrm{grad}}_{%
\mathbb{R}^{3}}f(\mathbf{x}) \\
\left( \mathbf{e}_{3}\cdot \mathbf{P}\right) \cdot \overline{\mathrm{grad}}_{%
\mathbb{R}^{3}}f(\mathbf{x})%
\end{matrix}%
\right] :=\left[
\begin{matrix}
\mathcal{G}_{1}f(\mathbf{x}) \\
\mathcal{G}_{2}f(\mathbf{x}) \\
\mathcal{G}_{3}f(\mathbf{x})%
\end{matrix}%
\right] ,  \label{eqn:Gs}
\end{equation}%
where $\mathcal{G}_{j}=\left( \mathbf{e}_{j}\cdot \mathbf{P}\right) \cdot
\overline{\mathrm{grad}}_{\mathbb{R}^{3}}$ $(j=1,2,3)$ is a differential
operator.

In a nutshell, the basic idea of the extrinsic formulation for the covariant
derivative follows from the definition of the tangential connection on a
submanifold $M$ of $\mathbb{R}^{n}$ (see e.g. Example 4.9 of \cite%
{lee2018introduction}). That is, for smooth vector fields $\boldsymbol{u},\boldsymbol{y}\in \mathfrak{X}%
\left( M\right) $, the tangential connection can be defined as%
\begin{equation}
\nabla _{\boldsymbol{u}}\boldsymbol{y}=\mathbf{P}\left( \bar{\nabla}_{\boldsymbol{U}}\boldsymbol{Y}|_{M}\right) ,
\label{eqn:nabuy}
\end{equation}%
where $\nabla $ is the Levi-Civita connection on $M$, $\bar{\nabla}:%
\mathfrak{X}\left( \mathbb{R}^{n}\right) \times \mathfrak{X}\left( \mathbb{R}%
^{n}\right) \rightarrow \mathfrak{X}\left( \mathbb{R}^{n}\right) $ is the
Euclidean connection (directional derivative) on $\mathbb{R}^{n}$\ mapping $(\boldsymbol{U},\boldsymbol{Y})$\ to $\bar{\nabla}%
_{\boldsymbol{U}}\boldsymbol{Y}$, $\mathbf{P}:T\mathbb{R}%
^{n}\rightarrow TM$ is the orthogonal projection onto $TM$, and $\boldsymbol{U}$ and $\boldsymbol{Y}$
are extensions of $u$ and $y$ onto an open subset in $\mathbb{R}^{n}$ such
that $\boldsymbol{U}|_{M}=\boldsymbol{u}$ and $\boldsymbol{Y}|_{M}=\boldsymbol{y}$, respectively. One can show that such extensions exist and the value is  independent of the chosen
extension and thus $\nabla _{\boldsymbol{u}}\boldsymbol{y}$ is well defined. More details about geometry can also be found in
\cite{do1992riemannian,morita2001geometry,crane2010trivial,2015Discrete} and
references therein.



We now review the extrinsic formulation in our previous work \cite{harlim2023radial} for the
gradient and the Bochner Laplacian of a 2D vector field when it has an
extension in $\mathbb{R}^3$. For a vector field $\boldsymbol{u}\in \mathfrak{X}(M),$ its
gradient in the local coordinate is defined to be a $(2,0)$ tensor field in $\mathfrak{X}%
(M)\times \mathfrak{X}(M)$, ${\mathrm{grad}}_{g}\boldsymbol{u}{=\sharp \nabla \boldsymbol{u}}%
=g^{kj}\tilde{u}_{;k}^{i}\frac{\partial }{\partial \xi ^{i}}\otimes \frac{\partial }{%
\partial \xi ^{j}}$, where $\sharp$ is the standard musical isomorphism notation to raise the index.
Let $\boldsymbol{U}=U^{k}\frac{\partial }{\partial X^{k}}$ be
a smooth extension of $\boldsymbol{u}=u^{i}\frac{\partial }{\partial \xi ^{i}}$ to an
open subset in $\mathbb{R}^{3}$\ satisfying $\boldsymbol{U}|_{M}=\boldsymbol{u}$. In fact, $\{\frac{%
\partial }{\partial X^{1}},\frac{\partial }{\partial X^{2}},\frac{\partial }{%
\partial X^{3}}\}=\left\{ \mathbf{e}_{1},\mathbf{e}_{2},\mathbf{e}%
_{3}\right\} $ is the standard orthonormal basis for the Cartesian
coordinate system. Let $\boldsymbol{U}=(U^{1},U^{2},U^{3})^{\top }\in \mathbb{R}^{3\times
1}$ be the coordinate representation and let $\frac{\partial }{\partial
\mathbf{X}}=(\frac{\partial }{\partial X^{1}},\frac{\partial }{\partial X^{2}%
},\frac{\partial }{\partial X^{3}})^{\top }\in \mathbb{R}^{3\times 1}$ be
the gradient operator $\overline{{\mathrm{grad}}}_{\mathbb{R}^{3}}$. Then
one can obtain the extrinsic formulation for the gradient of a vector field
(see e.g. \cite{harlim2023radial,gross2018trace}):
\begin{eqnarray}
{\mathrm{grad}}_{g}\boldsymbol{u} &=&\mathbf{P}\left( \overline{{\mathrm{grad}}}_{\mathbb{%
R}^{3}}\boldsymbol{U}\right) \mathbf{P}=\mathbf{P}\frac{\partial \boldsymbol{U}}{%
\partial \mathbf{X}}\mathbf{P}=\mathbf{P}\left[
\begin{matrix}
\mathcal{G}_{1}U^{1} & \mathcal{G}_{2}U^{1} & \mathcal{G}_{3}U^{1} \\
\mathcal{G}_{1}U^{2} & \mathcal{G}_{2}U^{2} & \mathcal{G}_{3}U^{2} \\
\mathcal{G}_{1}U^{3} & \mathcal{G}_{2}U^{3} & \mathcal{G}_{3}U^{3}%
\end{matrix}%
\right]  \label{eqn:gradu} \\
&:=&\left( \mathcal{H}_{1}\boldsymbol{U},\mathcal{H}_{2}\boldsymbol{U},\mathcal{H%
}_{3}\boldsymbol{U}\right),  \notag
\end{eqnarray}%
where for any $\mathbf{x}\in M$,  the projection matrix $\mathbf{P}=[P_{ij}]_{i,j=1}^{3}$ is
symmetric, $\overline{{\mathrm{grad}}}_{\mathbb{R}^{n}}\boldsymbol{U}=\frac{%
\partial \boldsymbol{U}}{\partial \mathbf{X}}=[\frac{\partial U^{i}}{\partial
X^{j}}]_{i,j=1}^{3}$ is a $3\times 3$ matrix, and $\mathcal{H}_{s}:=\mathbf{P%
}\mathrm{diag}\left( \mathcal{G}_{s},\mathcal{G}_{s},\mathcal{G}_{s}\right) $
has the form,
\begin{equation}
\mathcal{H}_{s}\boldsymbol{U}:=\mathbf{P}\left[
\setlength{\arraycolsep}{1.3pt}
\begin{matrix}
\mathcal{G}_{s} &  &  \\
& \mathcal{G}_{s} &  \\
&  & \mathcal{G}_{s}%
\end{matrix}%
\right] \left[
\begin{matrix}
U^{1} \\
U^{2} \\
U^{3}%
\end{matrix}%
\right] =
\left[
\begin{matrix}
	P_{11} & P_{12} & P_{13} \\
	P_{21} & P_{22} & P_{23} \\
	P_{31} & P_{32} & P_{33}%
\end{matrix}%
\right]
\left[
\setlength{\arraycolsep}{1.3pt}
\begin{matrix}
	\mathcal{G}_{s} &  &  \\
	& \mathcal{G}_{s} &  \\
	&  & \mathcal{G}_{s}%
\end{matrix}%
\right] \left[
\begin{matrix}
	U^{1} \\
	U^{2} \\
	U^{3}%
\end{matrix}%
\right]
\in \mathbb{R}^{3\times 1}.  \label{eqn:pV}
\end{equation}%

Using the above notations, one can eventually obtain the extrinsic formulation for
the Bochner Laplacian acting on a vector field (see Section 2.6 of \cite%
{harlim2023radial}),
\begin{equation}
\Delta _{B}\boldsymbol{u}:=g^{ij}u_{;ij}^{k}\frac{\partial }{\partial \xi ^{k}}=\mathcal{H}%
_{1}\mathcal{H}_{1}\boldsymbol{U}+ \mathcal{H}_{2}\mathcal{H}_{2}\boldsymbol{U} +\mathcal{H}_{3}%
\mathcal{H}_{3}\boldsymbol{U}:=\bar{\Delta}_{B}\boldsymbol{U}.  \label{eqn:LBH}
\end{equation}
Detailed derivations of the above formulation  can be found in \cite{harlim2023radial} or Section~\ref{app:newA}. Here,
the operator $\bar{\Delta}_B$ with overline defined on the right hand side of (\ref{eqn:LBH})
indicates that it involves with the extended vector field $\boldsymbol{U}$ on an open neighborhood of $M$ in $\mathbb{R}^n$ while one has $\bar{\Delta}_{B}\boldsymbol{U}  = \Delta _{B}\boldsymbol{u}$ restricted on $M$ as shown in \cite{harlim2023radial}. One can refer to Chap. 4 of \cite{lee2018introduction} for more about extensions of variables and operators.

\subsection{Review of local approximation of differential operators on
manifolds\label{sec:revdisH}}


We now apply the GMLS approach \cite{liang2013solving,flyer2016role,suchde2019meshfree,gross2020meshfree,jiang2024generalized,jones2023generalized} to approximate the Bochner Laplacian in (\ref%
{eqn:LBH}). In order to approximate the differential operators $\mathcal{H}_{s}$ in %
\eqref{eqn:pV} and subsequently $\bar{\Delta}_{B}$\ in \eqref{eqn:LBH}, we
first briefly review the local GMLS approximation to the differential
operator $\mathcal{G}_{s}$ in \eqref{eqn:Gs} (see e.g. \cite%
{jiang2024generalized}).

%
%
%


For a smooth function $f\in C^{\infty }(M)$, let $\mathbf{f}_{{\mathbf{x}_{0}}}=(f({\mathbf{x}_{0,1}}),...,f(\mathbf{x}%
_{0,K}))^{\top }$  be the restriction of $f$ on the
stencil $S_{{\mathbf{x}_{0}}}=\{{\mathbf{x}_{0,k}}\}_{k=1}^{K}$. Using the
GMLS approach as introduced in Section \ref{sec:gmlsm}, we can approximate
the differential operator $\mathcal{G}_{s}$ as,
\begin{equation*}
\mathcal{G}_{s}f({\mathbf{x}_{0,k}})\approx (\mathcal{G}_{s}\mathcal{I}_{%
\mathbb{P}}\mathbf{f}_{{\mathbf{x}_{0}}})({\mathbf{x}_{0,k}})=\sum_{0\leq
\alpha _{1}+\alpha _{2}\leq l}b_{\alpha _{1},\alpha _{2}}\mathcal{G}%
_{s}\theta _{1}^{\alpha _{1}}({\mathbf{x}_{0,k}})\theta _{2}^{\alpha _{2}}({%
\mathbf{x}_{0,k}}),\quad
\end{equation*}%
for any $k=1,...,K\ $and $s=1,2,3,$ where $\theta _{i}(\mathbf{x})$ $(i=1,2)$
are the local coordinates in (\ref{eq:thet12})$\ $and $b_{\alpha _{1},\alpha
_{2}}$\ are the coefficients for $(\mathcal{I}_{\mathbb{P}}\mathbf{f}_{{%
\mathbf{x}_{0}}})(\mathbf{x})=\sum b_{\alpha _{1},\alpha _{2}}\theta
_{1}^{\alpha _{1}}(\mathbf{x})\theta _{2}^{\alpha _{2}}(\mathbf{x})$ as in (%
\ref{eq:phinv}). For each $s=1,2,3$, the above relation can be written in
matrix form as,
\begin{eqnarray}
\left[
\begin{matrix}
\mathcal{G}_{s}f({\mathbf{x}_{0,1}}) \\
\vdots  \\
\mathcal{G}_{s}f({\mathbf{x}_{0,K}})%
\end{matrix}%
\right]  &\approx &\underbrace{\left[
\begin{matrix}
(\mathcal{G}_{s}\theta _{1}^{\alpha _{1}(1)}\theta _{2}^{\alpha _{2}(1)})(%
\mathbf{x}_{0,1}) & \cdots  & (\mathcal{G}_{s}\theta _{1}^{\alpha
_{1}(m)}\theta _{2}^{\alpha _{2}(m)})(\mathbf{x}_{0,1}) \\
\vdots  & \ddots  & \vdots  \\
(\mathcal{G}_{s}\theta _{1}^{\alpha _{1}(1)}\theta _{2}^{\alpha _{2}(1)})(%
\mathbf{x}_{0,K}) & \cdots  & (\mathcal{G}_{s}\theta _{1}^{\alpha
_{1}(m)}\theta _{2}^{\alpha _{2}(m)})(\mathbf{x}_{0,K})%
\end{matrix}%
\right] }_{\mathbf{B}_{s}}\underbrace{\left[
\begin{matrix}
b_{\boldsymbol{\alpha }(1)} \\
\vdots  \\
b_{\boldsymbol{\alpha }(m)}%
\end{matrix}%
\right] }_{\boldsymbol{b}}  \notag \\
&=&\mathbf{B}_{s}\boldsymbol{\Phi }^{\dag }\mathbf{f}_{{\mathbf{x}_{0}}}:=%
\mathbf{G}_{s}\mathbf{f}_{{\mathbf{x}_{0}}},  \label{eq:BPhi}
\end{eqnarray}%
where $\mathbf{B}_{s}$ is a $K$ by $m$ matrix with $[\mathbf{B}_{s}]_{ij}=(%
\mathcal{G}_{s}\theta _{1}^{\alpha _{1}(j)}\theta _{2}^{\alpha _{2}(j)})(%
\mathbf{x}_{0,i})$, $\boldsymbol{\Phi }^{\dag }$ is a $m$ by $K$ matrix
defined in (\ref{eq:phinv}) and $\mathbf{G}_{s}=\mathbf{B}_{s}\boldsymbol{\Phi }^{\dag }$ is a $K$ by $K$ matrix for
approximating $\mathcal{G}_{s}$ locally in the stencil $S_{{\mathbf{x}_{0}}}$%
.

Let $\mathbf{P}(\mathbf{x}_{0,k})=[P_{ij}(\mathbf{x}_{0,k})]_{i,j=1}^{3}$ be the projection matrix at the point $\mathbf{x}_{0,k}$.  Then we can define a tensor projection matrix $\mathbf{P}^{\otimes }\in \mathbb{R}%
^{3K\times 3K}$ as
\begin{eqnarray}
\mathbf{P}^{\otimes } &=&\sum_{k=1}^K\mathbf{P}(\mathbf{x}_{0,k}) \otimes \left[ \delta _{kk}\right] _{K\times K}=\sum_{k=1}^{K}\left[
\setlength{\arraycolsep}{2.3pt}
\begin{matrix}
P_{11}(\mathbf{x}_{0,k}) & P_{12}(\mathbf{x}_{0,k}) & P_{13}\left( \mathbf{x}%
_{0,k}\right) \\
P_{21}(\mathbf{x}_{0,k}) & P_{22}(\mathbf{x}_{0,k}) & P_{23}(\mathbf{x}%
_{0,k}) \\
P_{31}\left( \mathbf{x}_{0,k}\right) & P_{32}(\mathbf{x}_{0,k}) &
P_{33}\left( \mathbf{x}_{0,k}\right)%
\end{matrix}%
\right] \otimes \left[ \delta _{kk}\right] _{K\times K}  \notag \\
&=&\left[
\begin{matrix}
\mathrm{diag}(\mathbf{p}_{11}) & \mathrm{diag}(\mathbf{p}_{12}) & \mathrm{%
diag}(\mathbf{p}_{13}) \\
\mathrm{diag}(\mathbf{p}_{21}) & \mathrm{diag}(\mathbf{p}_{22}) & \mathrm{%
diag}(\mathbf{p}_{23}) \\
\mathrm{diag}(\mathbf{p}_{31}) & \mathrm{diag}(\mathbf{p}_{32}) & \mathrm{%
diag}(\mathbf{p}_{33})%
\end{matrix}%
\right] _{3K\times 3K},  \label{eq:dkk}
\end{eqnarray}%
where $\otimes $ is the Kronecker product, $\mathbf{p}_{ij}=\left( P_{ij}(%
\mathbf{x}_{0,1}),\ldots ,P_{ij}(\mathbf{x}_{0,K})\right) ^{\top }\in
\mathbb{R}^{K\times 1}$, and $\delta _{kk}\in \mathbb{R}^{K\times K}$ has
only one nonzero value $1$ in the $k$th row and $k$th column and has values $%
0$ elsewhere. Note that $\mathbf{P}^{\otimes }$ is a symmetric positive
definite projection matrix that maps a vector field in $\mathbb{R}^{3}$ onto
the tangent space of $M$ for the $K$ neighboring points in $S_{{\mathbf{x}%
_{0}}}$.

Then the differential operator $%
\mathcal{H}_{s}$\ in (\ref{eqn:pV}) can be approximated as,
\begin{equation}
\mathbf{H}_{s}=\mathbf{P}^{\otimes }\left[
\begin{matrix}
\mathbf{G}_{s} &  &  \\
& \mathbf{G}_{s} &  \\
&  & \mathbf{G}_{s}%
\end{matrix}%
\right] _{3K\times 3K},  \label{eqn:Hi}
\end{equation}
and the Bochner Laplacian in \eqref{eqn:LBH} can be subsequently approximated
in the stencil $S_{{\mathbf{x}_{0}}}$ as,
\begin{eqnarray}
(\Delta _{B}\boldsymbol{u})|_{S_{{\mathbf{x}_{0}}}} &=&(\bar{\Delta}_{B}\boldsymbol{U})|_{S_{{\mathbf{x}%
_{0}}}}=\sum_{\ell =1}^{3}(\mathcal{H}_{\ell }\mathcal{H}_{\ell }\boldsymbol{U})|_{S_{{%
\mathbf{x}_{0}}}}\approx \sum_{\ell =1}^{3}\mathbf{H}_{\ell }\mathbf{H}%
_{\ell }\mathbf{U}_{{\mathbf{x}_{0}}}  \notag \\
&=&\sum_{\ell =1}^{3}\mathbf{P}^{\otimes }\left[
\begin{matrix}
\mathbf{G}_{\ell } &  &  \\
& \mathbf{G}_{\ell } &  \\
&  & \mathbf{G}_{\ell }%
\end{matrix}%
\right] \mathbf{P}^{\otimes }\left[
\begin{matrix}
\mathbf{G}_{\ell } &  &  \\
& \mathbf{G}_{\ell } &  \\
&  & \mathbf{G}_{\ell }%
\end{matrix}%
\right] \mathbf{P}^{\otimes }\mathbf{U}_{{\mathbf{x}_{0}}},  \label{eqn:Sx0}
\end{eqnarray}%
where $\mathbf{U}_{{\mathbf{x}_{0}}}=((\mathbf{U}_{{\mathbf{x}_{0}}%
}^{1})^{\top },(\mathbf{U}_{{\mathbf{x}_{0}}}^{2})^{\top },(\mathbf{U}_{{%
\mathbf{x}_{0}}}^{3})^{\top })^{\top }\in \mathbb{R}^{3K\times 1}$ with each component
$\mathbf{U}_{{\mathbf{x}_{0}}}^{\ell }=(U^{\ell }({\mathbf{x}_{0,1}}%
),$ $...,$  $U^{\ell }(\mathbf{x}_{0,K}))^{\top }$ $\in \mathbb{R}^{K\times 1}$ for $\ell =1,2,3$,
and $\mathbf{U}_{{\mathbf{x}_{0}}}=\mathbf{P}^{\otimes }\mathbf{%
U}_{{\mathbf{x}_{0}}}$ since $\mathbf{U}_{{\mathbf{x}_{0}}}$ always lives in
the tangent bundle of $M$.

So far, we are following most procedures in our previous work \cite%
{harlim2023radial}, where the only difference is that now we use the GMLS
regression approach, which yields a sparse Laplacian matrix, instead of
a global RBF interpolation in \cite{harlim2023radial}, which yields a dense
matrix.\ However, the discrete Bochner Laplacian matrix is still of size $%
3N\times 3N$ which depends on the ambient dimension $n=3$ containing $9NK$ total nonzero entries. For other problems
with small intrinsic dimension $d$ and large ambient dimension $n$, the
computation would be expensive for solving PDEs involving the Bochner
Laplacian. In the next section, we introduce a dimension reduction technique to
reduce the size of the Bochner Laplacian from $3N\times 3N$ to $2N\times 2N$, thus reducing the number of nonzero entries from $3^2NK$ to $2^2NK$, where $2$ corresponds to
the intrinsic dimension.

%


\subsection{Reduction from ambient dimension to intrinsic dimension}

\label{sec:dimred}

For the dimension reduction, the key observation here is that the projection
matrix has the decomposition as $\mathbf{P(x)=T(x)T(x)}^{\top }$ (see e.g.
Proposition 2.1 of \cite{harlim2023radial}) for any $\mathbf{x}\in M$, where $%
\mathbf{T}(\mathbf{x})=[\boldsymbol{t}_{1}(\mathbf{x}),\boldsymbol{t}_{2}(%
\mathbf{x})]\in \mathbb{R}^{3\times 2}$ and $\boldsymbol{t}_{i}(\mathbf{x}%
)=(t^1_i(\mathbf{x}),t^2_{i}(\mathbf{x}),t^3_{i}(\mathbf{x}))^\top \in \mathbb{R}%
^{3\times 1}$ ($i=1,2$) are the same global tangent vector basis used in Section~\ref{sec:ABL}. Then, the
projection matrix $\mathbf{P}^{\otimes }\in \mathbb{R}^{3K\times 3K}$ for $%
K $ neighboring points in $S_{{\mathbf{x}_{0}}}$\ has the decomposition as $%
\mathbf{P}^{\otimes }=\mathbf{T}^{\otimes }\mathbf{T}^{\otimes \top },$
where $\mathbf{T}^{\otimes }\in \mathbb{R}^{3K\times 2K}$ is given by
\begin{equation}
\mathbf{T}^{\otimes }=\sum_{k=1}^{K}\left[ \mathbf{T}(\mathbf{x}_{0,k})%
\right] _{3\times 2}\otimes \left[ \delta _{kk}\right] _{K\times K}=\left[
\begin{matrix}
\mathrm{diag}(\boldsymbol{t}_{1,1}) & \mathrm{diag}(\boldsymbol{t}_{2,1}) \\
\mathrm{diag}(\boldsymbol{t}_{1,2}) & \mathrm{diag}(\boldsymbol{t}_{2,2}) \\
\mathrm{diag}(\boldsymbol{t}_{1,3}) & \mathrm{diag}(\boldsymbol{t}_{2,3})%
\end{matrix}%
\right] .\label{eqn:tensor-tangent}
\end{equation}%
Here, $\delta _{kk}$\ is defined as in (\ref{eq:dkk}) and $\boldsymbol{t}_{i,j}=\left( t_i^j(\mathbf{x}_{0,1}),\ldots ,t_i^j(%
\mathbf{x}_{0,K})\right) ^{\top }\in \mathbb{R}^{K\times 1}$ is the $j$th entry of the $i$th tangent vector $\boldsymbol{t}_i$ for $i=1,2$ and
$j=1,2,3$. Then the Bochner Laplacian (\ref{eqn:Sx0}) in the stencil $S_{{%
\mathbf{x}_{0}}}$ can be written as%
\begin{eqnarray}
\Delta _{B}\boldsymbol{u}|_{S_{{\mathbf{x}_{0}}}} &\approx& \sum_{\ell =1}^{3}\mathbf{T}%
^{\otimes }\left( \mathbf{T}^{\otimes \top }\left[
\setlength{\arraycolsep}{0.3pt}
\begin{matrix}
\mathbf{G}_{\ell } &  &  \\
& \mathbf{G}_{\ell } &  \\
&  & \mathbf{G}_{\ell }%
\end{matrix}%
\right] \mathbf{T}^{\otimes }\right) \left( \mathbf{T}^{\otimes \top }\left[
\setlength{\arraycolsep}{0.3pt}
\begin{matrix}
\mathbf{G}_{\ell } &  &  \\
& \mathbf{G}_{\ell } &  \\
&  & \mathbf{G}_{\ell }%
\end{matrix}%
\right] \mathbf{T}^{\otimes }\right)  \left( \mathbf{T}^{\otimes \top }%
\mathbf{U}_{{\mathbf{x}_{0}}}\right) \notag\\
  &:=& \mathbf{T}^{\otimes }(\sum_{\ell =1}^{3}\mathbf{R}_{\ell }\mathbf{R}%
_{\ell })\left( \mathbf{T}^{\otimes \top }\mathbf{U}_{{\mathbf{x}_{0}}%
}\right) \label{eqn:TTB} .
\end{eqnarray}%
Here, $\mathbf{R}_{\ell }:=\mathbf{T}^{\otimes \top }\left( \mathbf{I}%
_{3}\otimes \mathbf{G}_{\ell }\right) \mathbf{T}^{\otimes }\in \mathbb{R}%
^{2K\times 2K}$ with $\mathbf{I}_{3}$ being a $3\times 3$ identity matrix, $%
\mathbf{T}^{\otimes \top }\mathbf{U}_{{\mathbf{x}_{0}}}\in \mathbb{R}%
^{2K\times 1}$ is the coordinate representation of the vector field $\boldsymbol{u}|_{S_{{%
\mathbf{x}_{0}}}}$\ w.r.t. the global basis $\{\boldsymbol{t}_{j}(\mathbf{x}%
_{0,k})\}_{j=1,2}^{k=1,\ldots ,K}$\ and $\sum_{\ell =1}^{3}\mathbf{R}_{\ell }%
\mathbf{R}_{\ell }\in \mathbb{R}^{2K\times 2K}$\ is the Bochner Laplacian
matrix that maps $\boldsymbol{u}|_{S_{{\mathbf{x}_{0}}}}$\ to $(\Delta _{B}\boldsymbol{u})|_{S_{{%
\mathbf{x}_{0}}}}$\ w.r.t. the global basis $\{\boldsymbol{t}_{j}(\mathbf{x}%
_{0,k})\}$.


In order to find the weights $\{\boldsymbol{w}^{k}\}_{k=1}^{K}$ as in (\ref%
{eqn:LB}) to approximate $\Delta _{B}\boldsymbol{u}\left( {\mathbf{x}_{0}}\right) $, we
need to further analyze the structure of the differential matrix $\mathbf{R}%
_{\ell }$. From a geometric view point, the left multiplication of $\mathbf{T}%
^{\otimes \top }\ $and the right multiplication of $\mathbf{T}^{\otimes }$
in the\ matrix $\mathbf{R}_{\ell }=\mathbf{T}^{\otimes \top }\left( \mathbf{I%
}_{3}\otimes \mathbf{G}_{\ell }\right) \mathbf{T}^{\otimes }$ correspond to
the pointwise left projection and the pointwise right projection of the $%
3K\times 3K$\ matrix $\mathbf{I}_{3}\otimes \mathbf{G}_{\ell }$\ onto the
tangent spaces of $\{{\mathbf{x}_{0,k}}\}_{k=1}^{K}$. That is, the entries
in the $s^{\text{th}}$ and ($K+s$)$^{\text{th}}$\ rows of $\mathbf{R}_{\ell
} $\ (w.r.t. the point $\mathbf{x}_{0,s}$) and in the $r^{\text{th}}$\ and ($%
K+r$)$^{\text{th}}$\ columns\ (w.r.t. the point $\mathbf{x}_{0,r}$) can be
calculated by the following projection to form a $2\times 2$\ block
submatrix $\mathbf{\hat{R}}_{sr}^{\ell }$,
\begin{eqnarray}
\left[
\begin{matrix}
\boldsymbol{t}_{1}(\mathbf{x}_{0,s})^{\top } \\
\boldsymbol{t}_{2}(\mathbf{x}_{0,s})^{\top }%
\end{matrix}%
\right] _{2\times 3}&&\left[
\begin{matrix}
\left( \mathbf{G}_{\ell }\right) _{sr} & 0 & 0 \\
0 & \left( \mathbf{G}_{\ell }\right) _{sr} & 0 \\
0 & 0 & \left( \mathbf{G}_{\ell }\right) _{sr}%
\end{matrix}%
\right] _{3\times 3}\left[
\begin{matrix}
\boldsymbol{t}_{1}(\mathbf{x}_{0,r}) & \boldsymbol{t}_{2}(\mathbf{x}_{0,r})%
\end{matrix}%
\right] _{3\times 2}  \notag \\
&=&\left( \mathbf{G}_{\ell }\right) _{sr}\left[
\begin{matrix}
\boldsymbol{t}_{1}(\mathbf{x}_{0,s})^{\top } \\
\boldsymbol{t}_{2}(\mathbf{x}_{0,s})^{\top }%
\end{matrix}%
\right] _{2\times 3}\left[
\begin{matrix}
\boldsymbol{t}_{1}(\mathbf{x}_{0,r}) & \boldsymbol{t}_{2}(\mathbf{x}_{0,r})%
\end{matrix}%
\right] _{3\times 2}:=\mathbf{\hat{R}}_{sr}^{\ell }\in \mathbb{R}^{2\times
2},  \label{eqn:glst}
\end{eqnarray}%
where $\left( \mathbf{G}_{\ell }\right) _{sr}$ denotes the entry in the $s$%
th row and $r$th column of $\mathbf{G}_{\ell }$.

Using all these block submatrices $\hat{\mathbf{R}}_{sr}^{\ell }$ for $1\leq s,r\leq K$ and $\ell =1,2,3$%
, one can find the weights $\{\boldsymbol{w}^{k}\in \mathbb{R}^{2\times
2}\}_{k=1}^{K}$ satisfying
$\Delta_{B}\boldsymbol{u}({\mathbf{x}_{0}})\approx\sum_{k=1}^{K}\boldsymbol{w}^{k}\boldsymbol{u}|_{{\mathbf{x}_{0,k}}}$ in (\ref{eqn:LB}) by
\begin{equation}
(\boldsymbol{w}^{1},\ldots ,\boldsymbol{w}^{K})=\sum_{\ell =1}^{3}\left(
\mathbf{\hat{R}}_{11}^{\ell },\ldots ,\mathbf{\hat{R}}_{1K}^{\ell }\right) %
\left[
\setlength{\arraycolsep}{0.5pt}
\begin{matrix}
\mathbf{\hat{R}}_{11}^{\ell } & \cdots & \mathbf{\hat{R}}_{1K}^{\ell } \\
\vdots & \ddots & \vdots \\
\mathbf{\hat{R}}_{K1}^{\ell } & \cdots & \mathbf{\hat{R}}_{KK}^{\ell }%
\end{matrix}%
\right] :=\sum_{\ell =1}^{3}(\mathbf{\hat{R}}_{\ell }\mathbf{\hat{R}}_{\ell
})_{1}, \label{eqn:wk}
\end{equation}%
where $\mathbf{\hat{R}}_{\ell }$ is row and column equivalent to $\mathbf{R}_{\ell }$, that
is, $\mathbf{\hat{R}}_{\ell }$\ can be obtained from $\mathbf{R}_{\ell }$ in
(\ref{eqn:TTB})\ by taking appropriate row and column interchanges. Here
subscript$-1$ is to denote the elements in the first $2$ rows of the
enclosed $2K\times 2K$\ matrix $\mathbf{\hat{R}}_{\ell }\mathbf{\hat{R}}%
_{\ell }$\ corresponding to the Bochner Laplacian approximation at the base $%
{\mathbf{x}_{0}}$. Note that there is no need for a transformation of the
weights here since only the global basis $\{\boldsymbol{t}_{1}(\mathbf{x}%
_{i}),\boldsymbol{t}_{2}(\mathbf{x}_{i})\}_{i=1}^{N}$\ has been involved in
the computation for the extrinsic formulation of the Bochner Laplacian.
Following the above procedure, we can compute these weights for each base
point ${\mathbf{x}_{i}}$ $(i=1,\ldots ,N)$ and then appropriately arrange
them into a sparse $2N$ by $2N$ Bochner Laplacian matrix.

See Section~\ref{sec:inLlap} and Section~\ref{sec:apphodge} for the extrinsic approximations of two other vector Laplacians acting on vector fields for 2D manifolds in $\mathbb{R}^3$.

\subsection{Extrinsic formulation in higher dimensional spaces}\label{hihext}

Without loss of generality, the above extrinsic formulation can be readily extended to $d$-dimensional manifolds in $\mathbb{R}^n$. For simplicity, we present the generalized formulation of the Bochner Laplacian rather than other vector Laplacians.
The complete procedure for the extrinsic GMLS approximation of the Bochner Laplacian as well as the complexity is outlined in Algorithm~\ref{algo:extr-Boch}.

	
\begin{algorithm}[ht]
	\caption{Extrinsic GMLS of the Bochner Laplacian}
	\begin{algorithmic}[1]
		\STATE {\bf Input:} A point cloud $\{\mathbf{x}_{i}\}_{i=1}^{N}\subset M$, bases of (estimated) tangent vectors at each node $\left\{ \boldsymbol{t}_{1}(\mathbf{x}_{i}),\ldots,\boldsymbol{t}_{d}(%
		\mathbf{x}_{i})\right\} _{i=1}^{N}$, (estimated) projection matrices $\mathbf{P}$,  the degree $l$ of polynomials, and a parameter $K(>m)$ nearest neighbors.
		\STATE Set $\mathbf{L}_B$ to be a sparse $dN\times dN$ matrix with $d^2NK$ nonzeros.
		\FOR{$i\in \{1,2,...,N\}$}
		\STATE Find the $K$ nearest neighbors of the point $\mathbf{x}_i$ in the stencil $S_{\mathbf{x}_i}=\{\mathbf{x}_{i,k}\}_{k=1}^K$.
		\STATE  Construct the differential matrix $\mathbf{G}_\ell$ as in (\ref{eq:BPhi}) with complexity $O(NnK^2m + NnKmd^2)$.
		\STATE Calculate the block submatrices  $\mathbf{\hat{R}}^{\ell }_{sr} \in \mathbb{R}^{d\times d}$ similarly as in (\ref{eqn:glst}) for $1\leq s,r\leq K$ and $\ell=1,\ldots,n$ with complexity $O(Nnd^2K^2)$.
		\STATE Find the weights  $\{\boldsymbol{w}^{k}\}_{k=1}^{K}$ similarly as  in (\ref{eqn:wk}) with complexity $O(Nnd^3K^2)$.	
		\STATE Arrange the weights $\{\boldsymbol{w}^{k}\}_{k=1}^{K}$ into corresponding rows and columns of $\mathbf{L}_B$.
		\ENDFOR
		\STATE {\bf Output:} The approximate operator matrix $\mathbf{L}_B$.
	\end{algorithmic}
	\label{algo:extr-Boch}
\end{algorithm}

\subsection{Approximation of covariant derivative}\label{sec:acd}

We now compute the covariant derivative of a vector field $\boldsymbol{u}$\ along itself at the base
point $\mathbf{x}_{0}$, $\nabla _{\boldsymbol{u}}\boldsymbol{u}(\mathbf{x}_{0})$, in extrinsic
formulation using the $K$ neighboring points of $\mathbf{x}_{0}$ in $S_{%
\mathbf{x}_{0}}$. Let $\boldsymbol{U}=(U^{1},U^{2},U^{3})^{\top }$ be an extension of $%
\boldsymbol{u}=u^{j}\boldsymbol{t}_{j}$ so that $(U^{1},U^{2},U^{3})^{\top }=u^{1}%
\boldsymbol{t}_{1}+u^{2}\boldsymbol{t}_{2}$. The extrinsic formulation for
the covariant derivative can be written as (see e.g. \cite{harlim2023radial}),
\begin{equation}
\nabla _{\boldsymbol{u}}\boldsymbol{u}(\mathbf{x}_{0})=\mathbf{P}\left( \bar{\nabla}_{\boldsymbol{U}}\boldsymbol{U}\right) |_{%
\mathbf{x}_{0}}=\mathbf{TT}^{\top }\left. \left[
\begin{matrix}
\frac{\partial U^{1}}{\partial X^{1}} & \frac{\partial U^{1}}{\partial X^{2}}
& \frac{\partial U^{1}}{\partial X^{3}} \\
\frac{\partial U^{2}}{\partial X^{1}} & \frac{\partial U^{2}}{\partial X^{2}}
& \frac{\partial U^{2}}{\partial X^{3}} \\
\frac{\partial U^{3}}{\partial X^{1}} & \frac{\partial U^{3}}{\partial X^{2}}
& \frac{\partial U^{3}}{\partial X^{3}}%
\end{matrix}%
\right] \left[
\begin{matrix}
U^{1} \\
U^{2} \\
U^{3}%
\end{matrix}%
\right] \right\vert _{\mathbf{x}_{0}},  \label{eq:ddUx}
\end{equation}%
where $\mathbf{T}=[\boldsymbol{t}_{1},\boldsymbol{t}_{2}]\in \mathbb{R}%
^{3\times 2}$ is the global tangent vector basis. Then we approximate $\left[
\frac{\partial U^{r}}{\partial X^{s}}\right] _{r,s=1}^{3}$\ using GMLS
approach. Note that $\frac{\partial }{\partial X^{s}}=\left( \mathbf{e}%
_{s}\cdot \mathbf{I}\right) \cdot \overline{\mathrm{grad}}_{\mathbb{R}^{3}}$
can be obtained from $\mathcal{G}_{s}=\left( \mathbf{e}_{s}\cdot \mathbf{P}%
\right) \cdot \overline{\mathrm{grad}}_{\mathbb{R}^{3}}$\ in (\ref{eqn:Gs})
by replacing the projection $\mathbf{P}$ with the identity $\mathbf{I}$.
Following the similar formula as in (\ref{eq:BPhi}),\ one can construct a $K
$ by $K$ differential matrix $\mathbf{D}_{s}$\ to approximate $\frac{%
\partial }{\partial X^{s}}$ in $S_{\mathbf{x}_{0}}$. Taking the first row of
$\mathbf{D}_{s}$, denoted by $\mathbf{D}_{s,1}\in \mathbb{R}^{1\times K}$,
one can approximate the value of $\frac{\partial U^{r}}{\partial X^{s}}|_{%
\mathbf{x}_{0}}$ at the base $\mathbf{x}_0$ as $\mathbf{D}_{s,1}\mathbf{U}_{{\mathbf{x}_{0}}}^{r}$ for $%
1\leq r,s\leq 3$, where $\mathbf{U}_{{\mathbf{x}_{0}}}^{r}=(U^{r}({\mathbf{x}%
_{0,1}}),...,U^{r}(\mathbf{x}_{0,K}))^{\top }\in \mathbb{R}^{K\times 1}$ has
been defined right after (\ref{eqn:Sx0}). Then the covariant derivative (\ref%
{eq:ddUx}) can be calculated as%
\begin{eqnarray*}
\nabla _{\boldsymbol{u}}\boldsymbol{u}(\mathbf{x}_{0}) &=&\mathbf{T(\mathbf{x}_{0})T}(\mathbf{x}%
_{0})^{\top }\left[
\begin{matrix}
\mathbf{D}_{1,1}\mathbf{U}_{{\mathbf{x}_{0}}}^{1} & \mathbf{D}_{2,1}\mathbf{U}%
_{{\mathbf{x}_{0}}}^{1} & \mathbf{D}_{3,1}\mathbf{U}_{{\mathbf{x}_{0}}}^{1}
\\
\mathbf{D}_{1,1}\mathbf{U}_{{\mathbf{x}_{0}}}^{2} & \mathbf{D}_{2,1}\mathbf{U}%
_{{\mathbf{x}_{0}}}^{2} & \mathbf{D}_{3,1}\mathbf{U}_{{\mathbf{x}_{0}}}^{2}
\\
\mathbf{D}_{1,1}\mathbf{U}_{{\mathbf{x}_{0}}}^{3} & \mathbf{D}_{2,1}\mathbf{U}%
_{{\mathbf{x}_{0}}}^{3} & \mathbf{D}_{3,1}\mathbf{U}_{{\mathbf{x}_{0}}}^{3}%
\end{matrix}%
\right] \left[
\begin{matrix}
U^{1}(\mathbf{x}_{0}) \\
U^{2}(\mathbf{x}_{0}) \\
U^{3}(\mathbf{x}_{0})%
\end{matrix}%
\right]  \\
&=&\mathbf{T(\mathbf{x}_{0})T}(\mathbf{x}_{0})^{\top }\left[
\begin{matrix}
\left( \mathbf{U}_{{\mathbf{x}_{0}}}^{1}\right) ^{\top } \\
\left( \mathbf{U}_{{\mathbf{x}_{0}}}^{2}\right) ^{\top } \\
\left( \mathbf{U}_{{\mathbf{x}_{0}}}^{3}\right) ^{\top }%
\end{matrix}%
\right] _{3\times K}\left[
\begin{matrix}
\mathbf{D}_{1,1}^{\top } & \mathbf{D}_{2,1}^{\top } & \mathbf{D}_{3,1}^{\top }%
\end{matrix}%
\right] _{K\times 3}\left[
\begin{matrix}
U^{1}(\mathbf{x}_{0}) \\
U^{2}(\mathbf{x}_{0}) \\
U^{3}(\mathbf{x}_{0})%
\end{matrix}%
\right] .
\end{eqnarray*}

\section{Numerical results\label{sec:numr}}

To support the performance of both intrinsic and extrinsic GMLS approximations  as well as to validate the error estimates in Remark~\ref{rem:err2}, we present numerical examples, including eigenvalue problems, screened Poisson equations, linear vector diffusion equations, and nonlinear Burgers' equations. 
In Section \ref{sec:eigen}, we numerically
investigate the eigenvalue stability of the vector Laplacian on manifolds. In Section \ref{sec:poisson}, we verify the numerical convergence of the  GMLS approximation for solving the vector screened Poisson equation formulated with the Bochner Laplacian and the Hodge Laplacian. In Section~\ref{sec:lintim} and Section \ref{sec:burgers}, we study
the numerical performance of the linear vector diffusion equation and the nonlinear viscous Burgers' equation, respectively.

In our numerical implementation, for known smooth manifolds, we are given a point cloud associated with its  parametrization and tangent bundle.
For unknown smooth manifolds, we are only given a point cloud and  we apply the algorithm outlined in Section~\ref{sec:unkM} to approximate the tangent vectors of the manifolds without alignment as shown in Fig.~\ref{fig_tang}. Then we apply Algorithm~\ref{algo:intrin-Boch} and Algorithm~\ref{algo:extr-Boch} to obtain intrinsic and extrinsic approximations of vector Laplacians, respectively.  For time-independent problems, we choose the polynomial degree of manifolds in GMLS to be identical to that used for functions. For time-dependent problems, the degree of manifolds is chosen to be 6 which is higher than that used for functions ($l=2,3,4,5$) such that the error from manifold approximation is negligible. To achieve a stable and efficient approximation of the vector Laplacian, we empirically choose the $K$-nearest neighbors around $50 \sim 70$ for 2D manifolds, satisfying $K>(l+2)(l+1)/2$ where degree $l=5$. {The following numerical experiments demonstrate that the GMLS approaches can  provide stable estimators and convergent results with small error bars, given various trials of independent randomly sampled data points.}

\subsection{Eigenvalue stability}\label{sec:eigen}

\begin{figure*}
	\begin{minipage}[hbtp]{0.46\linewidth}
		\centering
		\begin{tabular}{c}
			{\footnotesize Unit Sphere (Bochner)} \\ 
			\includegraphics[scale=0.35]{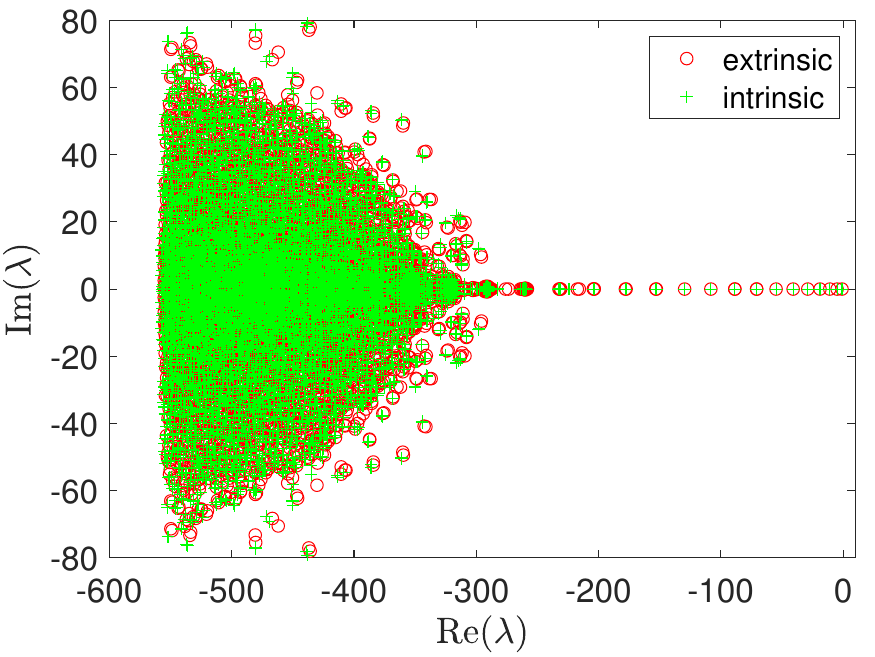} 
		\end{tabular}	
	\end{minipage}
	\begin{minipage}[hbtp]{0.45\linewidth}
		\centering
		\scalebox{0.85}[0.85]{
			\begin{tabular}{|c|c|c|c|}
				\hline
				\multicolumn{4}{|c|}{Unit Sphere (Bochner)} \\ \hline
				k & True & Extrinsic & Intrinsic \\ \hline
				1 & -1 & -0.9999 & -1.0002 \\ \hline
				2 & -1 & -0.9999 & -1.0002 \\ \hline
				3 & -1 & -0.9999 & -1.0002 \\ \hline
				4 & -1 & -1.0000 & -1.0003 \\ \hline
				5 & -1 & -1.0000 & -1.0003 \\ \hline
				6 & -1 & -1.0000 & -1.0003 \\ \hline
				7 & -5 & -5.0000 & -5.0003 \\ \hline
				8 & -5 & -5.0000 & -5.0003 \\ \hline						
		\end{tabular}}
	\end{minipage}	
	\caption{\textbf{Eigenvalues on 2D sphere in $\mathbb{R}^3$}. $K = 50$, $l = 5$, and $N = 6400$ points are randomly sampled. (Left) All eigenvalues lie in the left half  complex plane.  (Right) Comparison of the leading 8 eigenvalues of the Bochner Laplacian between our approximation using intrinsic and extrinsic GMLS methods and the analytic truth.}
	\label{fig:eigenv}
\end{figure*}

In this subsection, we investigate the eigenvalue problem for vector
Laplacians using two examples. For PDE applications, a necessary condition for stability is that all eigenvalues of the discrete vector Laplacian have non-positive real parts.
The first example involves the Bochner Laplacian on a 2D unit sphere embedded in $\mathbb{R}^3$.
The true eigenvalues of the Bochner Laplacian differ from those of the Hodge Laplacian
by a constant  Ricci curvature  on the unit sphere. The non-trivial eigenvalues of the Hodge Laplacian {coincide with} those of the Laplace-Beltrami operator but with double multiplicities (see e.g., \cite{harlim2023radial}). Neither of the two vector Laplacians has zero eigenvalues. In our numerical implementation, data points $\{\mathbf{x}_i\}_{i=1}^N$ are randomly sampled  with a uniform distribution on the unit sphere.
In the left panel of Fig. \ref{fig:eigenv}, we plot all the eigenvalues of the non-symmetric Bochner Laplacian matrices in the complex plane.
We see that all the eigenvalues lie in the left half complex plane with negative real parts, ensuring the eigenvalue stability of the GMLS approximation.
In the right panel of Fig. \ref{fig:eigenv}, we observe that the leading 8 numerical eigenvalues are in good agreement with the analytic true eigenvalues.



\begin{figure*}
	\begin{minipage}[hbtp]{0.46\linewidth}
		\centering
		\begin{tabular}{c}
			{\footnotesize (a) Torus in $\mathbb{R}^3$ (Hodge)} \\ 
			\includegraphics[scale=0.35]{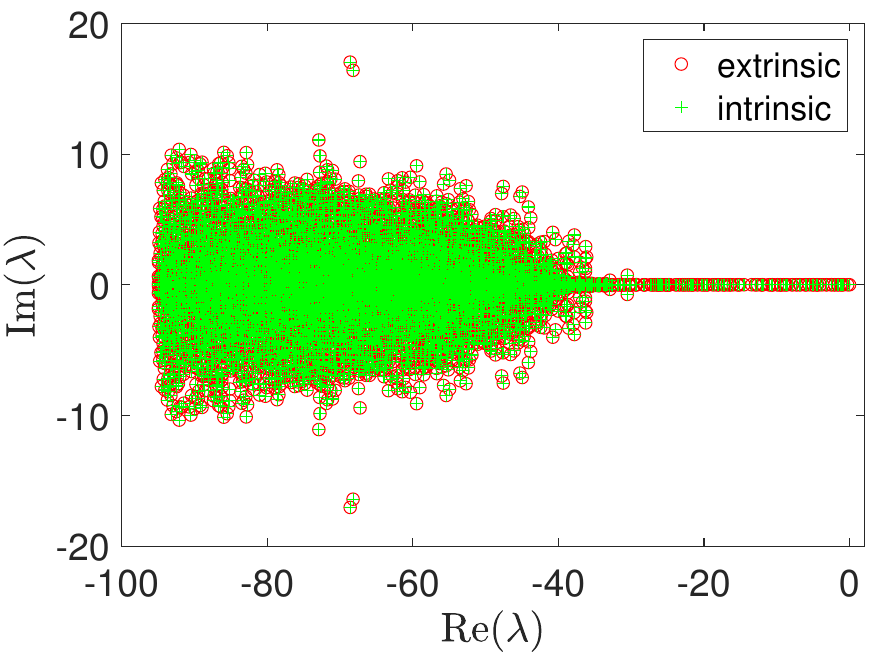} 
		\end{tabular}	
	\end{minipage}
	\begin{minipage}[hbtp]{0.45\linewidth}
		\centering
		\scalebox{0.85}[0.85]{
			\begin{tabular}{|c|c|c|c|}
				\hline
				\multicolumn{4}{|c|}{Torus in $\mathbb{R}^3$ (Hodge)} \\ \hline
				k & True & Extrinsic & Intrinsic \\ \hline
				1 & 0 & -2.7499e-04 & -7.2051e-04 \\ \hline
				2 & 0 & 0.0010 & -0.0050 \\ \hline
				3 & -0.2494 & -0.2480 & -0.2504 \\ \hline
				4 & -0.2494 & -0.2481 & -0.2504 \\ \hline
				5 & -0.2494 & -0.2494 & -0.2549 \\ \hline
				6 & -0.2494 & -0.2495 & -0.2551 \\ \hline
				7 & -0.7946 & -0.7930 & -0.7955 \\ \hline
				8 & -0.7946 & -0.7931 & -0.7956 \\ \hline						
		\end{tabular}}
	\end{minipage}	
	\\
	\begin{minipage}[hbtp]{1\linewidth}
		\vspace{+0.3in}
		\centering
		\scalebox{0.9}[0.9]{
			\begin{tabular}{ccc}	
				{\normalsize (b) Two Leading Eigenvalues} & {\normalsize (c) 1st Eigenvector Field} & {\normalsize (d) 2nd Eigenvector Field}\\
				\includegraphics[scale=0.28]{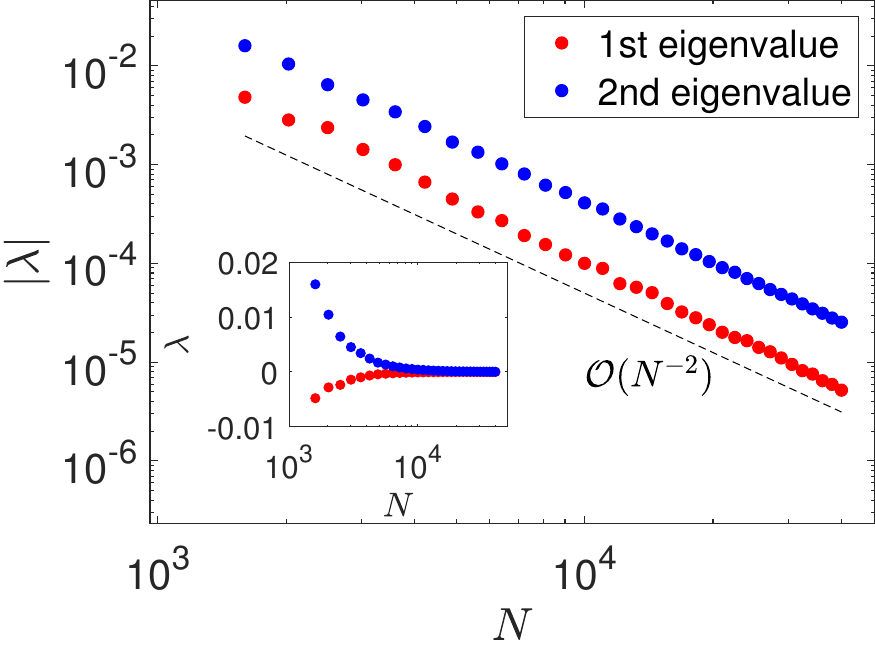} &		
				\includegraphics[scale=0.28]{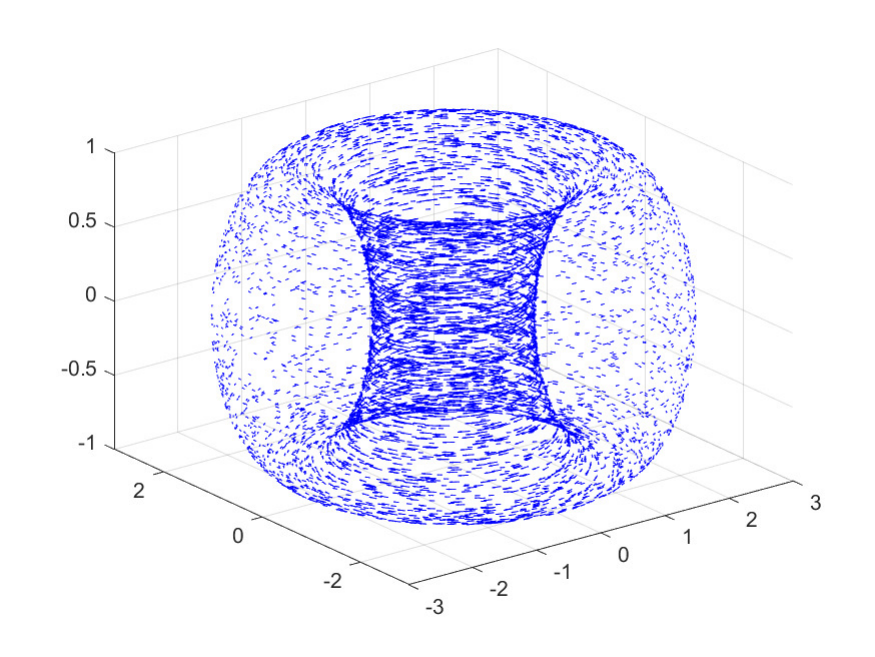} &
				\includegraphics[scale=0.28]{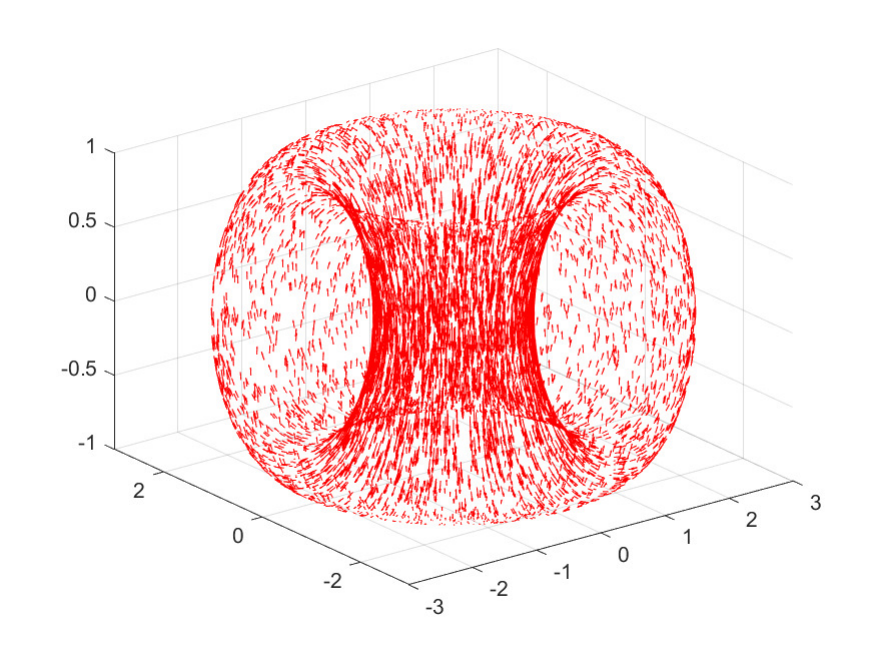}
		\end{tabular}}
	\end{minipage}	
	\vspace{+0.1in}
	\caption{\textbf{Eigenvalues on 2D torus in $\mathbb{R}^3$}. $K = 50$, $l = 5$, and $N = 6400$ points are randomly sampled.
	The first row: (a) Almost all eigenvalues lie in the left half complex plane. (Right) Comparison of the leading 8 eigenvalues of the Hodge Laplacian between the GMLS and the semi-analytic truth. The second row: (b) The first two eigenvalues of our approximation (with one positive and the other negative) converge to zero with $O(N^{-2})$.  (c) and (d) are the two corresponding harmonic fields of our approximation.}
		\label{fig:eigenv_torus}
\end{figure*}

The second example involves the Hodge Laplacian on a torus parametrized by
\begin{equation}
	\mathbf{x} = \left[
	(R+r\cos\theta)\cos\phi ,
	(R+r\cos\theta)\sin\phi ,
	r\sin\theta \right], \quad 0 \leq \theta < 2\pi, 0 \leq \phi < 2\pi,
	\label{eqn:torusr3}
\end{equation}
where $R$ is the distance from the center of the tube to the center of the torus and $r< R$ is the distance from the center of the tube to the surface of the tube. We set $R=2$ and $r=1$.
The eigenvalues of the Hodge Laplacian are identical to those of the Laplace-Beltrami but with double multiplicities. Notably, the Hodge Laplacian on a torus has a zero eigenvalue with multiplicity two, associated with two harmonic fields. The nontrivial true eigenvalues can be approximated semi-analytically by solving the Sturm-Liouville form of the Laplace-Beltrami operator using the finite difference method (see \textit{p}. 68 of \cite{prycel1993SL}).

Numerically, the data $\{\theta_i,\phi_i\}_{i=1}^N$ are generated from a uniform distribution on $\left[0,2\pi\right) \times \left[0,2\pi\right)$ and then mapped onto the point cloud $\{\mathbf{x}_i\}_{i=1}^N$ via the parametrization (\ref{eqn:torusr3}). We observe from the first row of Fig.~\ref{fig:eigenv_torus} that almost all eigenvalues lie in the left half plane except a value of 0.001 which approximates the zero eigenvalue. Next, we study the behavior of the two numerical eigenvalues close to zero. From the second row of Fig.~\ref{fig:eigenv_torus}, we see that these two numerical eigenvalues converge to the true zero value while one of them remains positive. Additionally, the associated numerical eigenvector fields approximate the two harmonic fields of the Hodge Laplacian, one parallel to the $\theta$ direction and the other parallel to the $\phi$ direction.  Before explaining the numerical observation for the positive eigenvalue, we first note that for a closed manifold, the
discrete GMLS Laplace-Beltrami matrix always has a zero eigenvalue, with an associated constant eigenvector in the
null space of the discrete matrix. However, this is not the case for a discrete vector Laplacian. A vector Laplacian may not have a zero eigenvalue (e.g., on the unit sphere), or if it does, the corresponding eigenvector field will not  have constant components. Due to approximation errors, the leading numerical eigenvalue might be slightly greater than zero but it will converge to zero when a true zero  eigenvalue exists, as shown in Fig.~\ref{fig:eigenv_torus}.

In the following PDE applications, we always test for the eigenvalue stability before solving the PDE. If a numerical eigenvalue $\hat{\lambda}_0$ has a small positive real part (as in the torus example), which characterizes the approximation error of the zero eigenvalue in this case, we will regularize the vector Laplacian matrix by subtracting a constant from its diagonal entries, $\mathbf{L}_B-\hat{\lambda}_0 \mathbf{I}$, to guarantee the eigenvalue stability where $\mathbf{I}$ is an identity matrix.


\subsection{Screened Poisson problems}
\label{sec:poisson}
In this subsection, we consider solving the following vector-valued screened Poisson problem for $\boldsymbol{u}\in \mathfrak{X}(M)$,
\begin{equation}
	(a - \Delta_B) \boldsymbol{u} = \boldsymbol{f},\quad  \mathbf{x} \in M,  \label{eqn:poisson}
\end{equation}
where  $a=1>0$ and $\boldsymbol{f}\in \mathfrak{X}(M)$ are defined such that the problems are well-posed. For simplicity, we focus on the Bochner Laplacian in the following examples while Section~\ref{sec:521} presents an example using the Hodge Laplacian. We note that the vector fields $\boldsymbol{f}$ and $\boldsymbol{u}$ must live in the tangent bundle $TM$ of the smooth manifold $M$. In the case of an unknown manifold, $\boldsymbol{f}$  is replaced by $\mathbf{\hat{T}}^\top \boldsymbol{F}$ in our numerical implementation where $\mathbf{\hat{T}}=[\boldsymbol{\hat{t}}_1,\ldots, \boldsymbol{\hat{t}}_d ]$ corresponds to the approximate tangent vectors and $\boldsymbol{F}\in \mathfrak{X}(\mathbb{R}^n)$ is a vector field in ambient space. Thus, the vector field $\mathbf{\hat{T}}^\top \boldsymbol{F}$ and the approximate solution $\hat{\mathbf{u}}=(a-\mathbf{L}_B)^{-1}(\mathbf{\hat{T}}^\top \boldsymbol{F})$ both live in the approximate tangent bundle $\widehat{TM}:=\cup_{\mathbf{x}\in M}\widehat{T_\mathbf{x} M}=\cup_{\mathbf{x}\in M}\mathrm{Span}\{\boldsymbol{\hat{t}}_1(\mathbf{x}),\ldots, \boldsymbol{\hat{t}}_d(\mathbf{x})\}$.
To facilitate comparison with true solutions in the following numerical examples, we prescribe an analytic solution $\boldsymbol{u}$ and then  approximate it, subjected to the manufactured  $\boldsymbol{f}:=\mathbf{\hat{T}}^\top \boldsymbol{F} = \mathbf{\hat{T}}^\top (a - \bar{\Delta}_B) \boldsymbol{U}$, where $(a - \bar{\Delta}_B) \boldsymbol{U}$ is an extension of $(a - \Delta_B) \boldsymbol{u}$, as defined in (\ref{eqn:LBH}).

To verify convergence, we define the forward error (\textbf{FE}) and the inverse error (\textbf{IE}) in the ambient space $\mathbb{R}^n$ as
	\begin{equation}
		\textbf{FE} = \max_{1 \leq i \leq N} \|  \hat{\textbf{T}}(\mathbf{x}_i)(\mathbf{L}_B\mathbf{u}_{X_M})_i - \textbf{T}(\mathbf{x}_i)\Delta_B \boldsymbol{u}(\mathbf{x}_i) \|_2,
		\mathrm{\ } \textbf{IE} = \max_{1 \leq i \leq N} \| (\hat{\textbf{T}}\hat{\mathbf{u}})(\mathbf{x}_i) - (\mathbf{T}\boldsymbol{u})(\mathbf{x}_i) \|_2,
		\nonumber
	\end{equation}
where $\Vert \cdot \Vert_2$ is the Euclidean distance, the true solution $\mathbf{u}_{X_M}=(\boldsymbol{u}(\mathbf{x}_1),\ldots,\boldsymbol{u}(\mathbf{x}_N))^\top\in \mathbb{R}^{dN\times 1}$ is a column vector,   $(\mathbf{L}_B\mathbf{u}_{X_M})_i \in \mathbb{R}^{d\times 1}$ corresponds to the approximate value of $\Delta_B \boldsymbol{u}$ at $\mathbf{x}_i$,
$\Delta_B \boldsymbol{u}(\mathbf{x}_i)\in \mathbb{R}^{d\times 1}$ is the representation in the global basis at $\mathbf{x}_i$, $\hat{\mathbf{u}}$ is the approximate solution, $\mathbf{T}(\mathbf{x})$ $=[\boldsymbol{t}_{1}(\mathbf{x}),$$\ldots,\boldsymbol{t}_{d}(\mathbf{x})]$ $\in \mathbb{R}^{n\times d}$
and $\mathbf{\hat{T}}(\mathbf{x})=[\boldsymbol{\hat{t}}_{1}(\mathbf{x}),\ldots,\boldsymbol{\hat{t}}_{d}(\mathbf{x})]\in \mathbb{R}^{n\times d}$ correspond to the true and the approximate global bases, respectively. Here we push both tangent vectors on (approximate) manifolds forward to tangent vectors in $\mathbb{R}^n$ to make the \textbf{FE} and \textbf{IE} well-defined in the Euclidean norm.


\subsubsection{Hodge Laplacian on 2D sphere in $\mathbb{R}^3$\label{sec:521}}
In this numerical experiment, we solve the screened Poisson problem (\ref{eqn:poisson}) on a 2D unit sphere associated with the Hodge Laplacian $\Delta_H$, which  typically can be addressed using the FEEC framework. Here, we compare the performance of the GMLS approach proposed in this paper with the global RBF method \cite{harlim2023radial} (a meshless approach) and the FEEC method \cite{holst2012geometric} (a mesh-based approach), as discussed in Section \ref{sec:intro}. We set the true solution $\boldsymbol{u} = (x-x^3)\mathbf{e}_1 - (x^2y)\mathbf{e}_2 - (x^2z)\mathbf{e}_3$ $\in \mathfrak{X}(M)$, where $\mathbf{x}=(x,y,z)$ is the Euclidean coordinate. For the GMLS approximation, we fix $K=50$ for $K$-nearest neighbors. In the RBF method, we select $\Psi(\mathbf{x})=\exp(-s^2\Vert \mathbf{x} - \mathbf{x}_0\Vert^2)$ as the kernel function with a shape parameter $s=1$. For the surface FEM solution within the FEEC framework, we choose FEM spaces for differential forms as in \cite{holst2012geometric}, ensuring that the FEM spaces have the same polynomial degree as the surface approximation.  Computations are performed using the NGSolve package \cite{Schoberl2014C11IO}. For a consistent comparison, we generate triangle meshes over the surface with varying maximal mesh sizes. The mesh vertices are then extracted as a point cloud for the GMLS and RBF methods. The \textbf{IE}s are calculated as before at the vertices.
	
\begin{figure*}[htbp]
	{\scriptsize \centering
		\setlength{\tabcolsep}{1pt}
		\begin{tabular}{cccc}
			{\normalsize (a) (extr.) \textbf{%
					GMLS}} & {\normalsize (b) (intr.) \textbf{GMLS}} & {\normalsize (c) \textbf{%
					FEEC $\&$ RBF}}  & {\normalsize (d) \textbf{%
					time}}\\
			\includegraphics[width=1.25
			in, height=1.1 in]{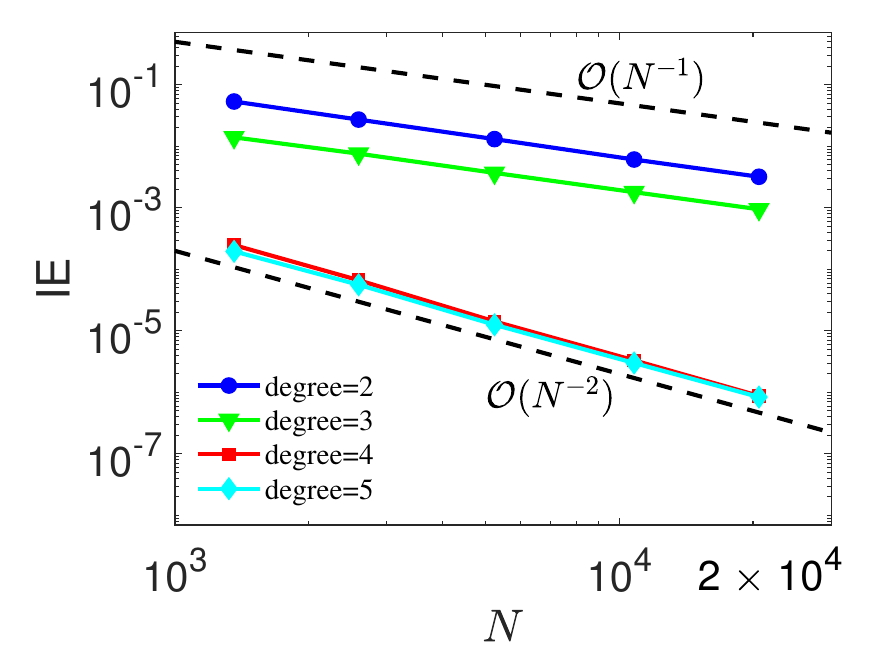} &
			\includegraphics[width=1.25
			in, height=1.1 in]{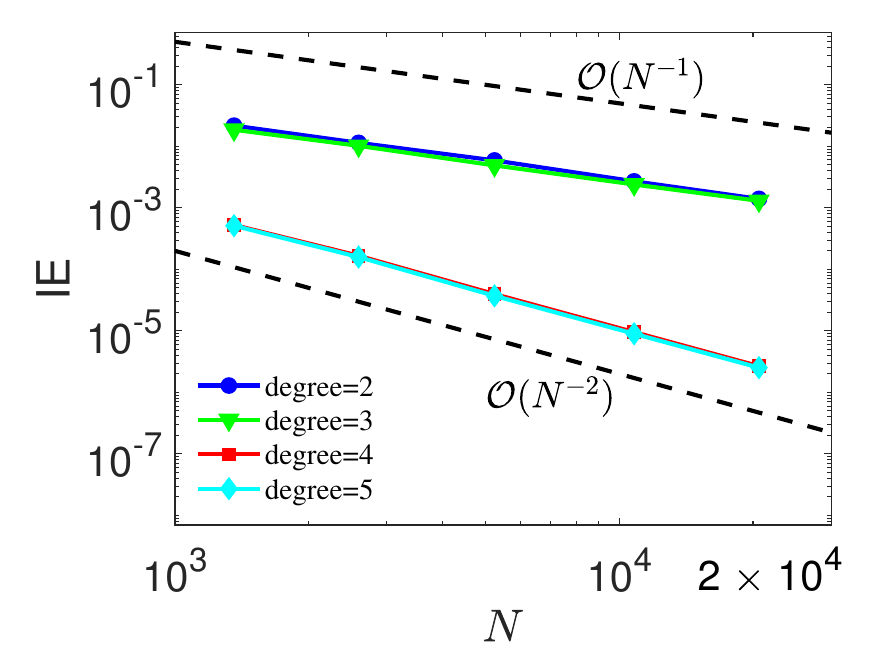} &
			\includegraphics[width=1.25
			in, height=1.1 in]{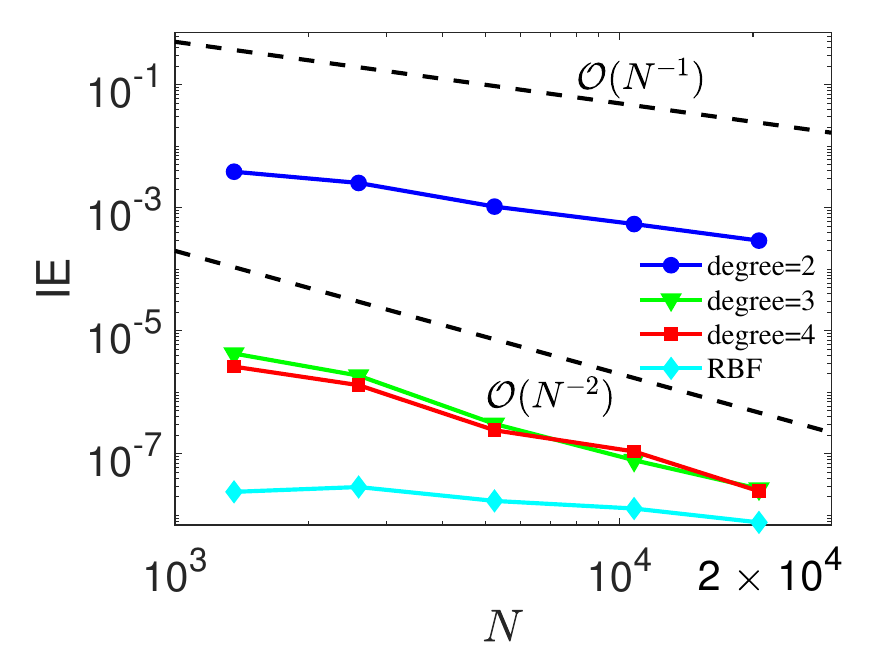} &
			\includegraphics[width=1.25
			in, height=1.1 in]{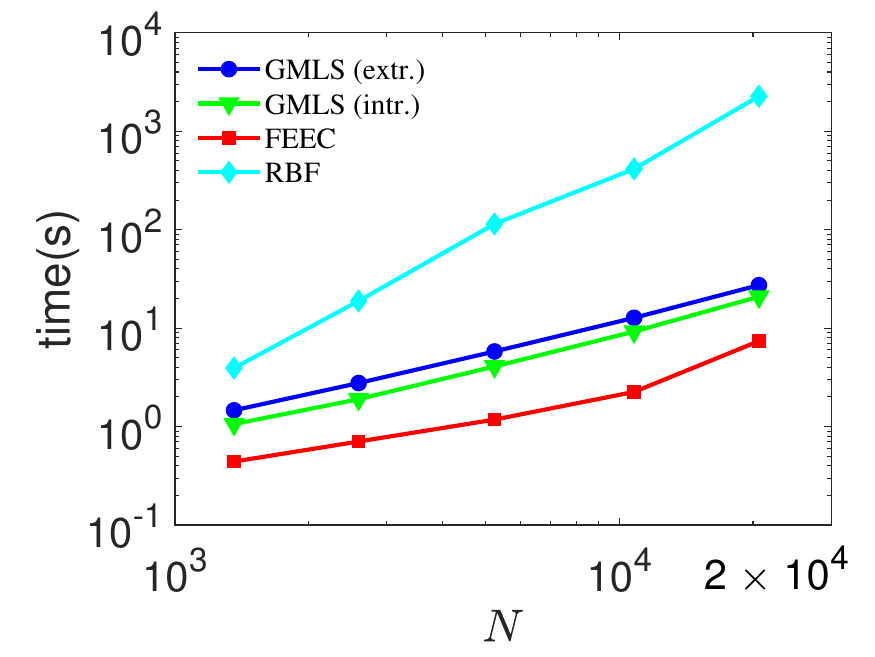}
		\end{tabular}
	}
	\caption{\textbf{Screened Poisson problems on 2D sphere in $\mathbb{R}^3$ associated with the Hodge
	Laplacian}. Random data
	with 6 trials. $K=50$ nearest neighbors are used. In panels (a) and (b), shown are the average \textbf{IE}s  using extrinsic and intrinsic methods, respectively. In panels (c), shown are the \textbf{IE}s of the solutions using FEEC and RBF. In panels (d), shown are the comparison of the time cost among different methods. For time comparison, we use degree 2 for GMLS and FEEC.}
	
	\label{fig:hodge}
\end{figure*}

\comment{	
	(Doug FEEC) Mixed methods
	for the Hodge Laplacian problem:
	\begin{equation}
		\begin{aligned}
			\langle \sigma, \tau\rangle - \langle d\tau,u\rangle &= 0,\ \tau \in V^{k-1}, \\
			\langle d\tau, v\rangle + \langle du,dv \rangle + \langle v,p \rangle &= \langle f,v\rangle, \  v \in V^k, \\
			\langle u,q\rangle &= 0, q \in \mathfrak{H}^k.
		\end{aligned}
\end{equation}}

In Figs. \ref{fig:hodge}(a)-(c),  we plot the $\mathbf{IE}$s as functions of $N$. The $\mathbf{IE}$s of the GMLS solutions decay at rates of $N^{-1}$  for $l=2,3$ and $N^{-2}$ for $l=4,5$, which agrees with the GMLS error estimates  (\ref{eq:gmlrt}). In Fig. \ref{fig:hodge}(c), we observe that for $l=2$ and $l=4$, the FEEC solution converges at the same rate as the GMLS methods but yields smaller errors. This behavior is expected as high-quality surface meshes can be easily constructed in this example. On the other hand, the RBF method provides the most accurate solution.  In Fig. \ref{fig:hodge}(d), we plot computational time as a function of $N$. We can see that the time cost of the RBF method increasing significantly faster than that of the GMLS and the FEEC methods. This is due to the fact that the RBF method produces a full matrix whereas the other two methods result in sparse matrices.

\subsubsection{Bochner Laplacian on 2D torus in $\mathbb{R}^9$}

\begin{figure*}[htbp]
{\scriptsize \centering
	\setlength{\tabcolsep}{1pt}
\begin{tabular}{cccc}
{\normalsize (a) (extr.) \textbf{%
FE}} & {\normalsize (b) (intr.) \textbf{FE}} & {\normalsize (c) (extr.) \textbf{%
IE}} & {\normalsize (d) (intr.) \textbf{%
IE}} \\
\includegraphics[width=1.25 in, height=1.1 in]{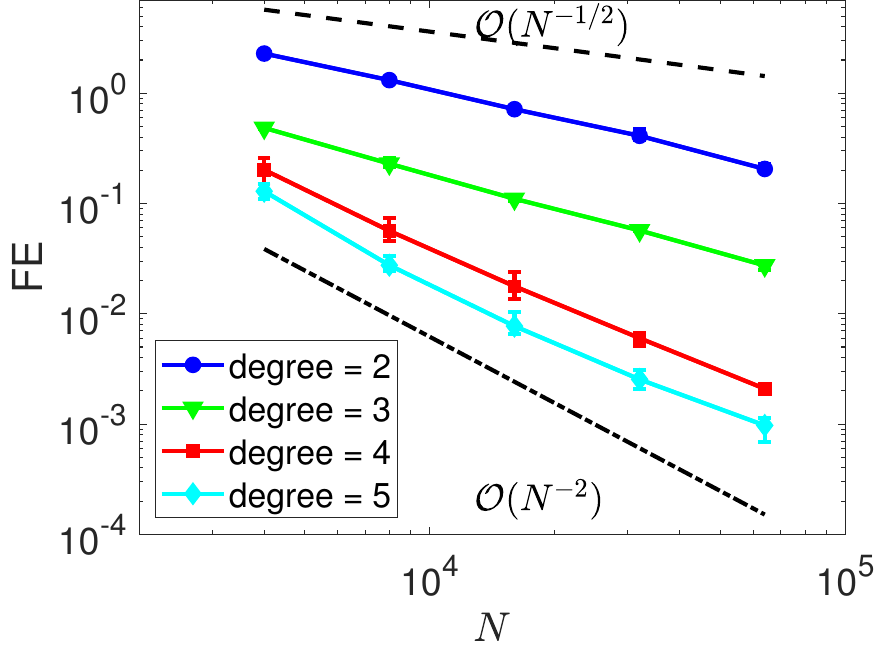} &
\includegraphics[width=1.25
			in, height=1.1 in]{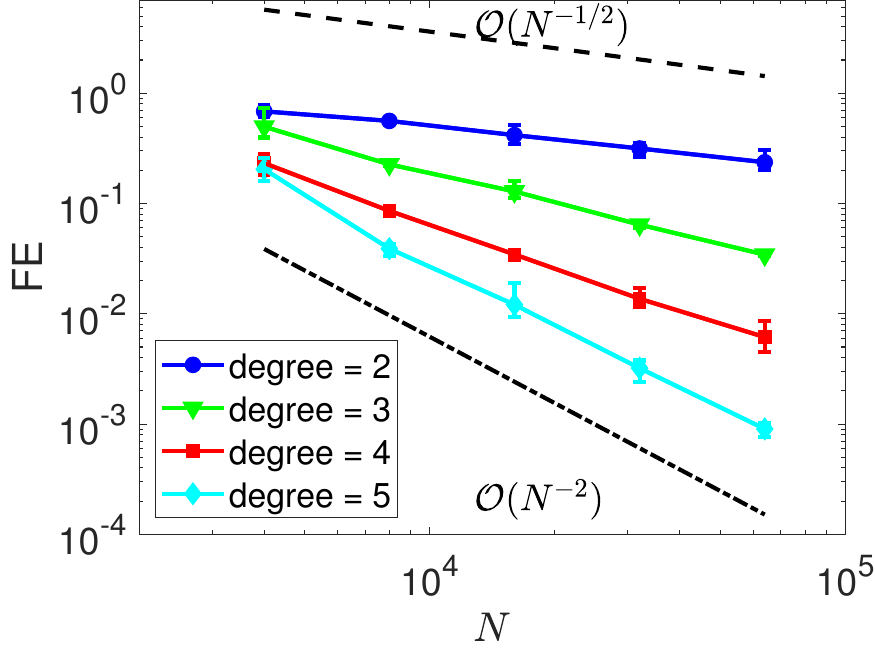} &
\includegraphics[width=1.25
			in, height=1.1 in]{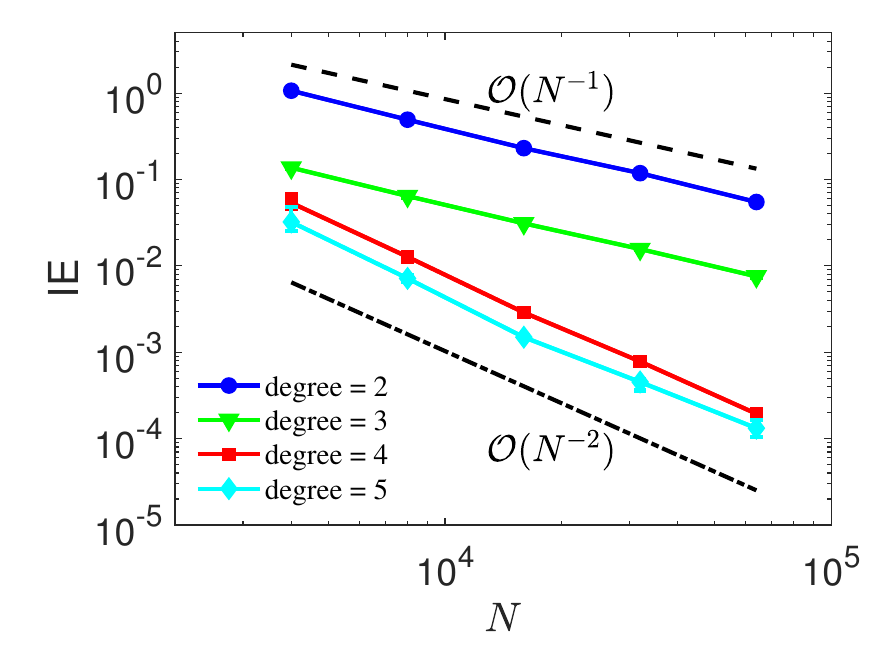} &
\includegraphics[width=1.25
			in, height=1.1 in]{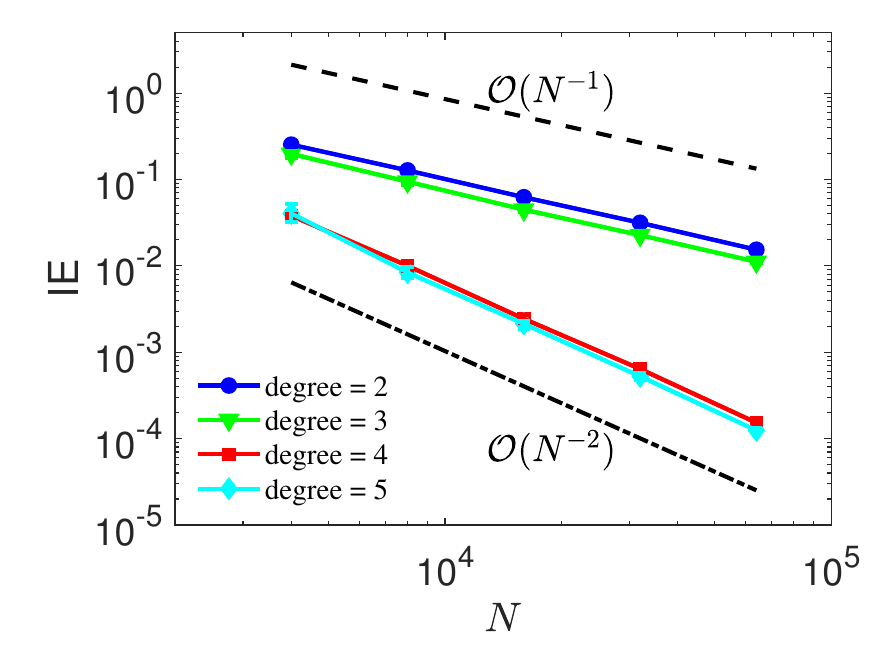}%
\end{tabular}
}
\caption{\textbf{Screened Poisson problems on 2D torus in $\mathbb{R}^9$}. In panels (a) and (b), shown are the average \textbf{FE}s of the approximate Bochner
Laplacian $\mathbf{L}_B$ using extrinsic and intrinsic methods, respectively. In panels (c) and (d), shown are the average \textbf{IE}s of the solutions of the screened Poisson problem (\ref{eqn:poisson}).
}
\label{fig:kkernel}
\end{figure*}


In this numerical experiment, we consider solving the screened Poisson problem (\ref{eqn:poisson}) on a 2D torus in $\mathbb{R}^9$,
\begin{equation}
		\mathbf{x} = \left[\begin{matrix}
		(R+\cos\theta)\cos\phi \\
		(R+\cos\theta)\sin\phi \\
		\vdots \\
		(1/4)(R+\cos\theta)\cos4\phi \\
		(1/4)(R+\cos\theta)\sin4\phi \\
		\sqrt{\sum_{i=1}^4 (1/i)^2}\sin\theta
	\end{matrix}\right]\in \mathbb{R}^9, \quad 0 \leq \theta < 2\pi, \quad 0 \leq \phi < 2\pi,
	\label{eqn:torus}
\end{equation}
where $R=2$ is the distance from the center of the tube to the center of the torus. The true solution $\boldsymbol{u}$  is
set to be $\boldsymbol{u} = \sin \theta \sin \phi \frac{\partial}{\partial \theta} + \sin \theta \cos \phi \frac{\partial}{\partial \phi}$ $\in \mathfrak{X}(M)$.  The point cloud $\{\mathbf{x}_i\}_{i=1}^N$ is generated using the parametrization (\ref{eqn:torus}), where $\{\theta_i,\phi_i\}_{i=1}^N$ are randomly sampled from a uniform distribution on $\left[0,2\pi\right) \times \left[0,2\pi\right)$. We fix $K=50$ in $K$-nearest neighbors.






In Fig. \ref%
{fig:kkernel} , we plot the averages of $\mathbf{FE}$s and $\mathbf{IE}$s
over 4 independent trials as functions of $N$. One can see that the $\mathbf{FE}$s using degree $%
l=2,3,4,5$ decrease on the order of $N^{-(l-1)/2}$ for both the
extrinsic and intrinsic methods. This numerical error rate is consistent with the GMLS error bound in (\ref{eq:gmlrt})
since the Bochner Laplacian is involved with the second-order derivative of component functions. One can further observe that the $\mathbf{IE}$s for $l=3$ and $l=5$ decay with respective rates $N^{-1}$ and $N^{-2}$ as the GMLS theory predicts. Moreover, the $\mathbf{IE}$s possess the super-convergence when degree $l=2$ and $l=4$, that is, the $\mathbf{IE}$s decay on the order of $N^{-1}$ and $N^{-2}$, respectively, which are half-order faster than the corresponding $\mathbf{FE}$s.
This super-convergence phenomenon has also been observed in the scalar case for the Laplace-Beltrami operator of functions \cite{liang2013solving,jiang2024generalized}.


\comment{Plotted in panels (a)(b) are the \textbf{FE}s of the approximate Bochner
Laplacian $\mathbf{L}_B$. Plotted in panels (c)(d) are the \textbf{IE}s of the solution of the Poisson problem. Panels (a)(c) correspond to  extrinsic formulation while panels (b)(d) correspond to
intrinsic formulation.}

\subsubsection{3D flat torus in $\mathbb{R}^{12}$}\label{sec:523}

We consider a $3$-dimensional manifold embedded in $\mathbb{R}^{12}$ with
the following parameterization,
\begin{equation}
	\label{eqn:flattorus}
	\mathbf{x}=\frac{1}{\sqrt{1^2 +2^{2}}}\left(
	\begin{array}{cccc}
		\cos (\phi_{1}), & \sin (\phi_{1}),  & \cos (2\phi_1), & \sin (2\phi_{1}), \\
		\cos (\phi_{2}), & \sin (\phi_{2}),  & \cos (2\phi_2), & \sin (2\phi_{2}), \\
		\cos (\phi_{3}), & \sin (\phi_{3}),  & \cos (2\phi_3), & \sin (2\phi_{3})
	\end{array}%
	\right) ,
\end{equation}%
with $0\leq \phi_{1},\phi_2,\phi_3< 2\pi$. The Riemannian
metric is given by a $3\times 3$ identity matrix $\mathbf{I}_{3}$.  The true solution $\boldsymbol{u}$ is set to be $\boldsymbol{u} =\sin (\phi_1) \sin (\phi_2) \frac{\partial}{\partial \phi_{1}} + \sin (\phi_{2}) \sin (\phi_{3}) \frac{\partial}{\partial \phi_{2}} + \sin (\phi_{3}) \cos (\phi_{1}) \frac{\partial}{\partial \phi_{3}}$. Numerically, the points $\{\mathbf{x}_i\}_{i=1}^N$ are generated from the parametrization \eqref{eqn:flattorus} using randomly sampled $\{\phi_{1}(\mathbf{x}_i),\phi_{2}(\mathbf{x}_i),\phi_{3}(\mathbf{x}_i)\}_{i=1}^N$ with
the uniform distribution on $\left[0,2\pi\right)^3$.
We use the degree $l=3,4$ and $K=75$ nearest neighbors for illustrating the convergence of solutions.

\begin{figure*}
	\centering
	\begin{tabular}{ccc}
		{\normalsize (a) (extr.) $\textbf{IE}$} & {\normalsize (b) (intr.) $\textbf{IE}$} & {\normalsize (c) Pointwise Absolute Error}\\
		\includegraphics[width=1.4
		in, height=1.3 in]{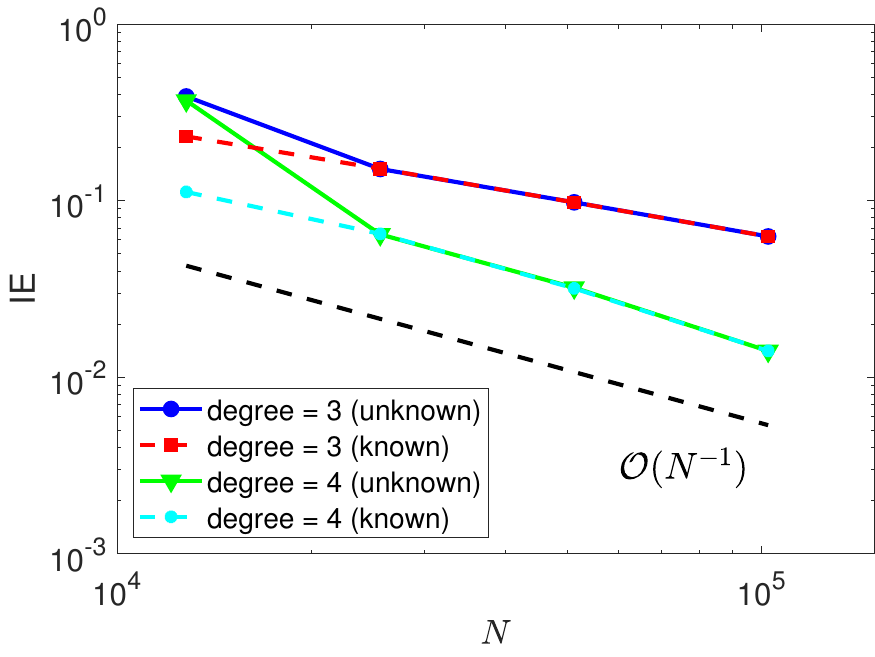} &
		\includegraphics[width=1.4
		in, height=1.3 in]{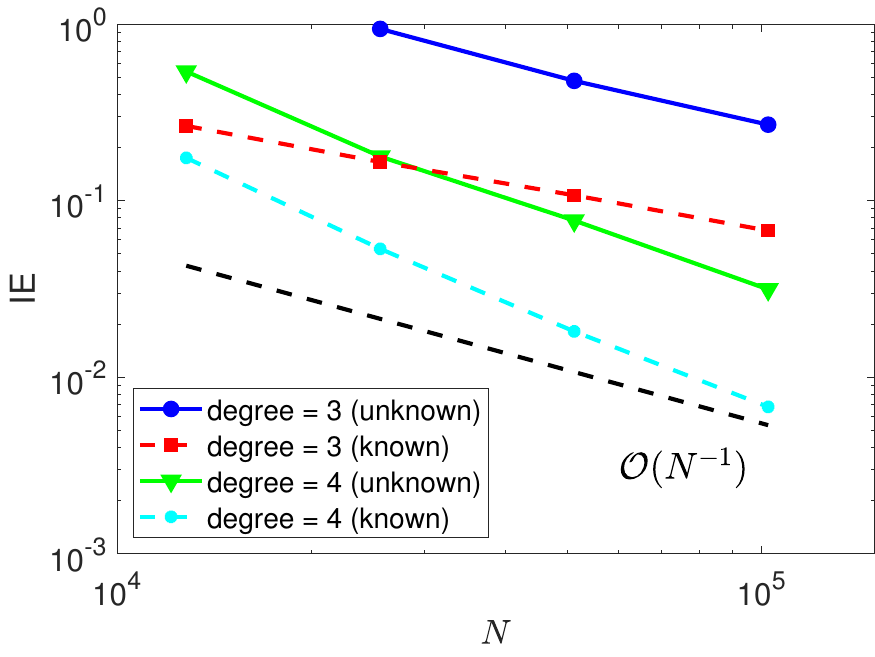} &
		\includegraphics[width=1.4
		in, height=1.3 in]{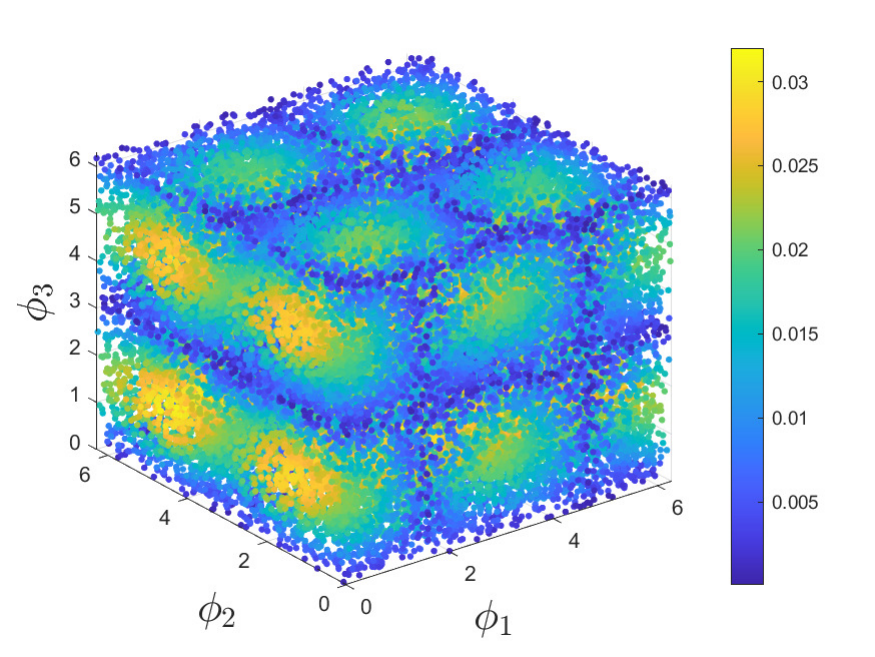}
	\end{tabular}
	\caption{\textbf{Screened Poisson problems on 3D flat torus in $\mathbb{R}^{12}$}. $K=75$ nearest neighbors and degree $l=3,4$ are used. In panels (a) and (b), shown are the \textbf{IE}s of the solutions using our extrinsic and intrinsic methods, respectively. In panel (c), shown is the pointwise absolute error of the solution using extrinsic method for an unknown manifold setup with $N = 51200$.
	}
		\label{fig:flattorus}
\end{figure*}

In Fig. \ref{fig:flattorus}, we illustrate the convergence of $\textbf{IE}s$ for both the extrinsic and intrinsic GMLS methods.
From Figs.~\ref{fig:flattorus}(a) and (b), we observe that both methods provide a stable and convergent approximation to the solution for both known and unknown manifold setups.
Additionally, the extrinsic method exhibits nearly the same errors between known and unknown manifold setups. However, for the intrinsic method, the error for the unknown manifold is much larger than that for the known manifold while the convergence rates hold.
Upon further examination, we identify that the large error gap in the intrinsic method for this 3D example stems from the local regression
approximation of Monge parametrization in Step 8 of Algorithm~\ref{algo:intrin-Boch}.
Specifically, the tangent space of the local GMLS  approximation of Monge parametrization ($\mathrm{span}\{\partial_1,\partial_2,\partial_3\}$) deviates from that of the manifold  ($\mathrm{span }\{\boldsymbol{t}_1,\boldsymbol{t}_2,\boldsymbol{t}_3\}$) at some neighboring points in the local stencil $S_{\mathbf{x}_0}$. To improve the  result for the intrinsic method, one potential approach is to use an interpolation (see e.g. \cite{flyer2016role,jones2023generalized}) instead of a  regression  for approximating the local manifold.
This suggests a direction for future research.



\subsection{Linear time-dependent equations}\label{sec:lintim}
In this subsection, we will show the results of solving the following linear time-dependent vector-valued  PDE (vector diffusion equation):
\begin{equation}
	\boldsymbol{u}_t = \nu \Delta_B \boldsymbol{u} + \boldsymbol{f}, \quad (\mathbf{x},t) \in M \times \left(0,T
	\right] \text{ with } \boldsymbol{u}(\cdot,0) = \boldsymbol{u}_0(\cdot),  \label{eqn:linPDE}
\end{equation}
where the viscosity $\nu=0.1$.
Here, we consider two examples of surfaces that are both isomorphic to a sphere in $\mathbb{R}^3$.
The first example is a red blood cell (RBC) with the parametrization:
\begin{equation}
	\textbf{x} = \left(r\cos\theta\cos\phi,r\cos\theta\sin\phi, \frac{1}{2}\sin\theta \left( c_0+c_2\cos^2\theta+c_4\cos^4\theta\right) \right),
	\label{eqn:rbc}
\end{equation}
where $-\pi/2 \leq \theta \leq \pi/2$, $-\pi \leq \phi < \pi$, $r=3.91/3.39$, $c_0=0.81/3.39$, $c_2=7.83/3.39$ and $c_4=-4.39/3.39$, which is the same as the RBC in \cite{fuselier2013high}. The second example is a bumpy sphere (BSP) with the parametrization:
\begin{equation}
	\textbf{x} = \left(r(\theta,\phi)\sin\theta\cos\phi,r(\theta,\phi)\sin\theta\sin\phi, r(\theta,\phi)\cos\theta \right),
	\label{eqn:psp}
\end{equation}
where $0 \leq \theta \leq \pi$, $0 \leq \phi < 2\pi$ and $r(\theta,\phi) = 1 + 0.2(\sin^7 3\theta)(\sin4\phi)$.

In our numerical experiments, we set the external forcing $\boldsymbol{f}(\textbf{x},t):=\hat{\textbf{T}}^{\top}\textbf{x}$ in both examples to be time-independent. The initial condition is $\boldsymbol{u}_0(\textbf{x}):=0.2\hat{\textbf{T}}^{\top}\textbf{x}$ for RBC and $\boldsymbol{u}_0(\textbf{x}):=0.1\hat{\textbf{T}}^{\top}\textbf{x}$ for BSP.
To solve the time evolution equation (\ref{eqn:linPDE}), we use the GMLS method for the space discretization and the second-order explicit Runge-Kutta method for the time discretization. The data are randomly sampled from the intrinsic coordinate and then mapped onto the surface using  the above parametrization.
To examine the convergence, we define the time-dependent solution error (\textbf{SE}) between true solution $\boldsymbol{u}$ and approximate solution $\hat{\mathbf{u}}$:
\begin{equation*}
\textbf{SE} = \max_{1 \leq i \leq N} \| \hat{\textbf{T}}(\mathbf{x}_i)\hat{\mathbf{u}}(\mathbf{x}_i,t) - \mathbf{T}(\mathbf{x}_i)\boldsymbol{u}(\mathbf{x}_i,t) \|_2.
\end{equation*}

\comment{{\color{red}super-convergence check using bdf2}}

\comment{{\color{red}isomorphic  sampling density}}

\begin{figure*}[htbp]
	\begin{minipage}[!b]{1\linewidth}
		\centering
		\includegraphics[width=5
		in, height=2 in]{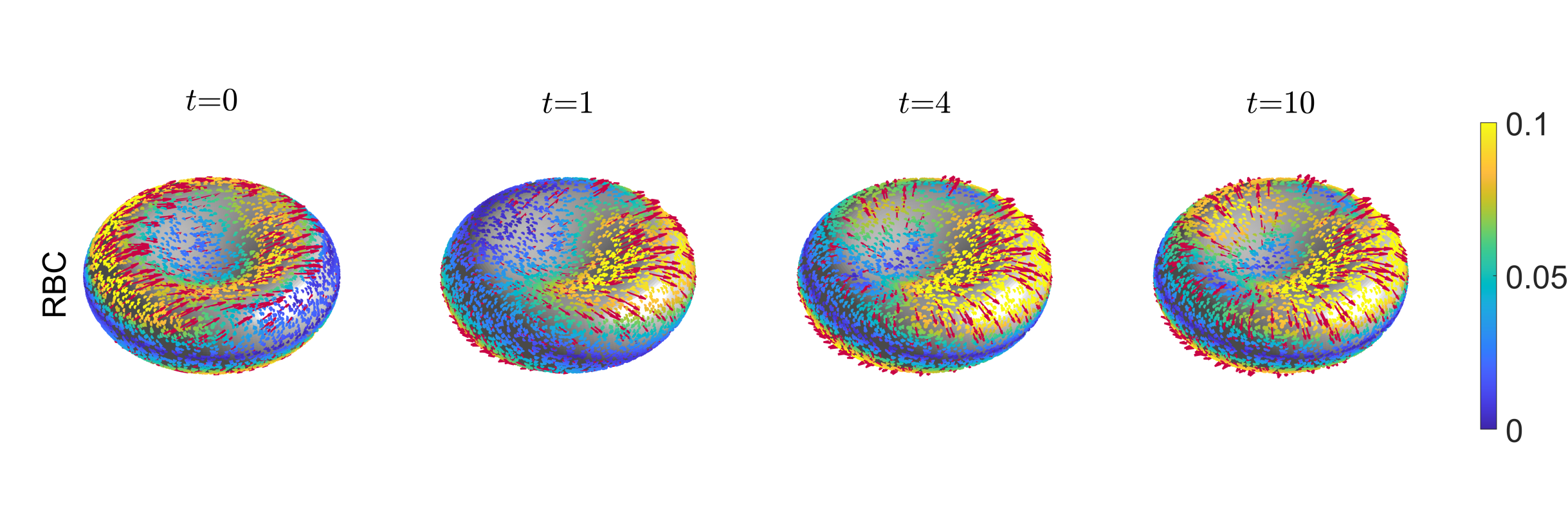}
	\end{minipage}

	\begin{minipage}[!b]{1\linewidth}
		\vspace{-0.1in}
		\centering
		\includegraphics[width=5 in, height=1.7 in]{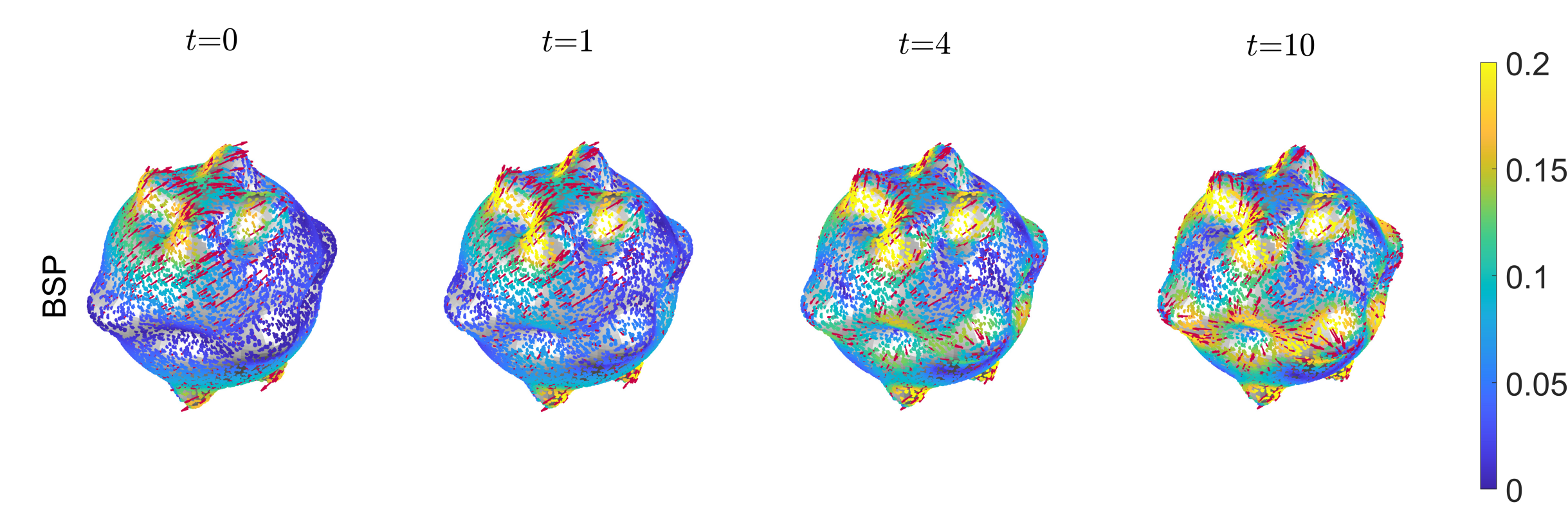}
	\end{minipage}
	\begin{minipage}[tbph]{1\linewidth}
		\centering
		\begin{tabular}{cccc}
			{\normalsize (a)
				RBC (extr.) } & {\normalsize (b) RBC (intr.) } &
			{\normalsize (c)
				BSP (extr.) } & {\normalsize (d) BSP (intr.) }	
			\\
			\includegraphics[width=1.1
			in, height=1. in]{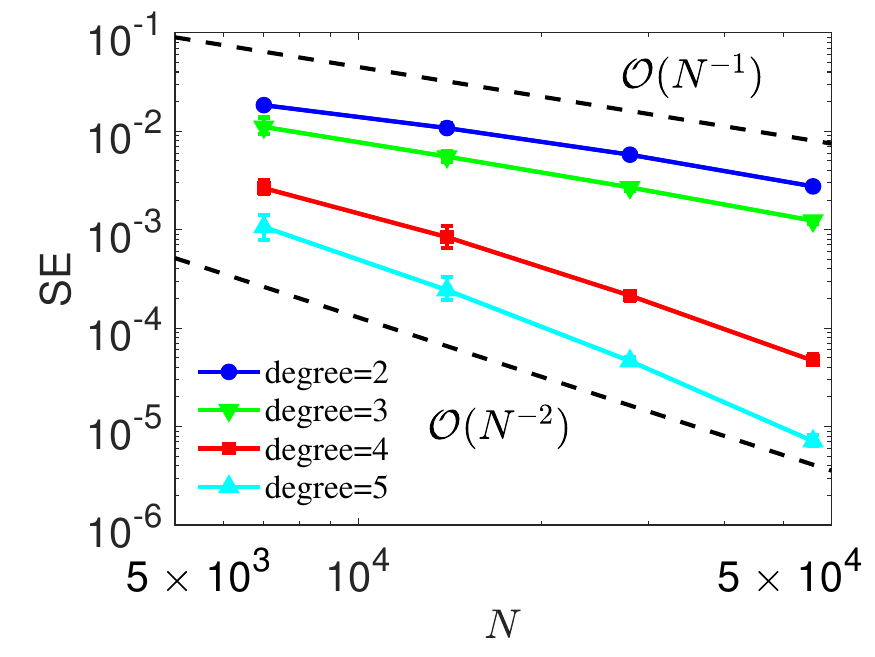} &
			\includegraphics[width=1.1
			in, height=1. in]{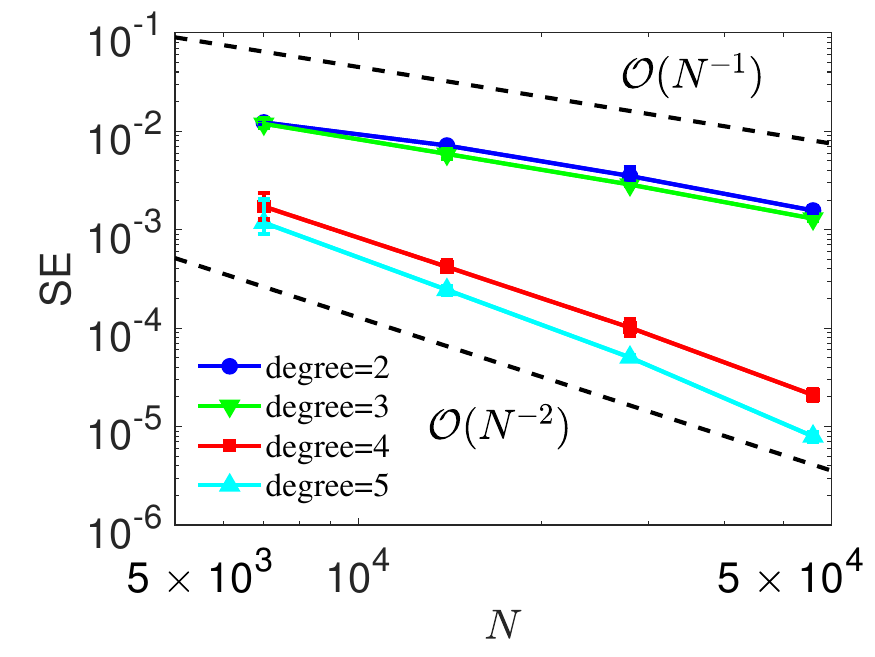} &
			\includegraphics[width=1.1
			in, height=1. in]{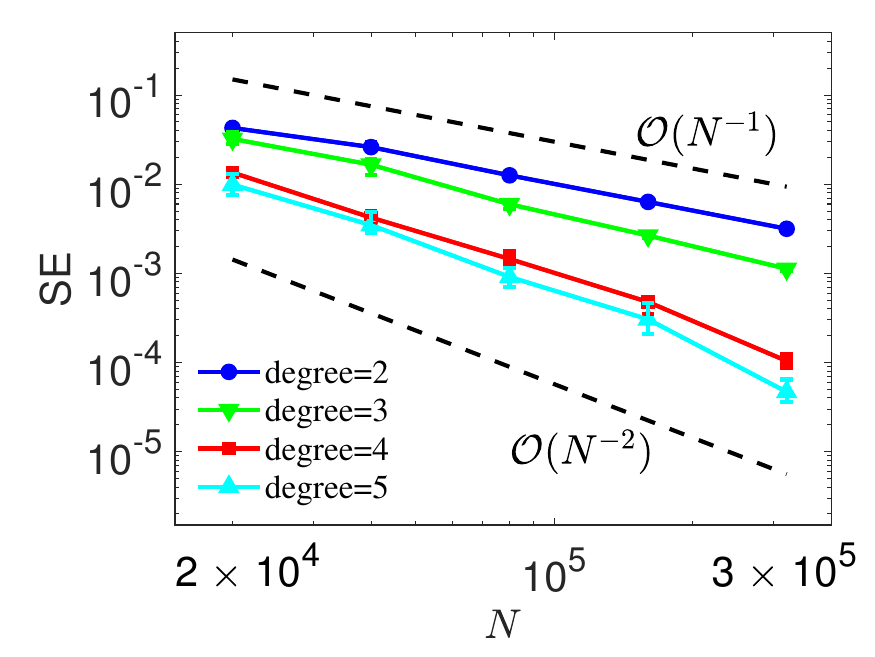} &
			\includegraphics[width=1.1
			in, height=1. in]{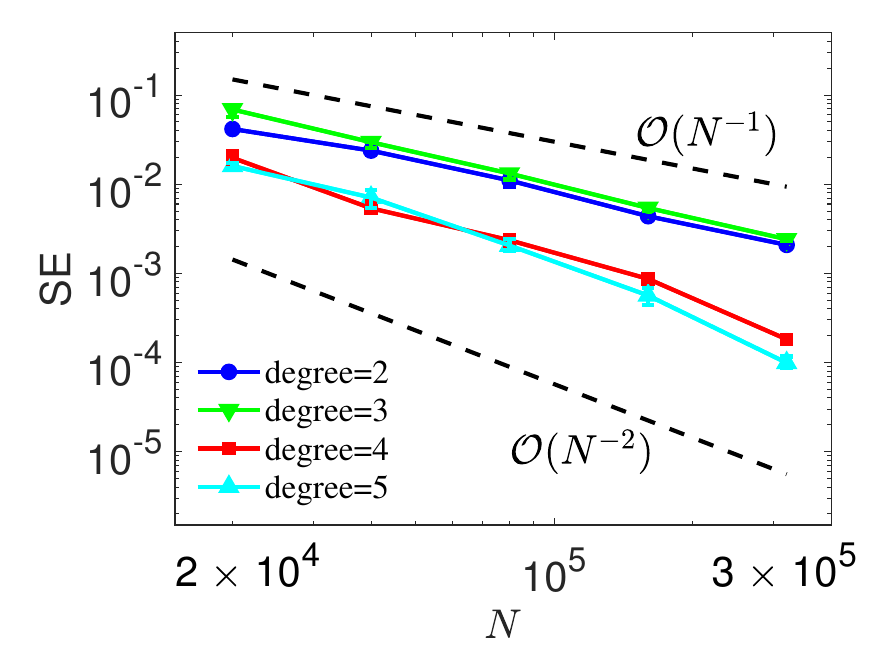}%
		\end{tabular}
	\end{minipage}
	\caption{\textbf{Linear vector diffusion equation on RBC and BSP}. For both examples, we set $K=60$ in $K$ nearest neighbors. The time step $\Delta t=0.001$ for RBC and $\Delta t=0.0001$ for BSP.  Shown in the first two rows are time evolution of the extrinsic GMLS solutions on RBC (the first row) and BSP (the second row) with the color bar  representing the magnitude of the vector fields $\hat{\mathbf{u}}$  and the red arrows indicating the direction of the vector fields $\hat{\mathbf{u}}$. We use $N = 20000$ randomly sampled data points, degree $l=2$ for RBC and $N=160000$ randomly sampled data points, degree $l=4$ for BSP. In the third row, we show the convergence examination of the average \textbf{SE}s at $t=0.05$ with error bars obtained from the standard deviations over 6 independent trials.}
	\label{fig:linear_eqn}
\end{figure*}

In Fig. \ref{fig:linear_eqn}, we show our numerical results on RBC and BSP.  The first two rows illustrate the time evolution of the extrinsic GMLS solutions on surfaces using $N=20000$ data points for RBC and $N = 160000$ data points for BSP at times $t=[0,1,4,10]$. After a sufficiently long time ($t=10$), the solution reaches its steady state, corresponding to the solution of the associated Poisson equation. The last row shows the average \textbf{SE}s for the numerical solutions as functions of $N$  at $t=0.05$.  To compute the reference solution, we apply the same discretization schemes over $N=80000$ randomly sampled points for RBC and $N=640000$ for BSP, with polynomial degree $l=5$, and analytic true tangent space spanned by $\boldsymbol{t}_1,\boldsymbol{t}_2$. To verify the convergence, we compute the numerical solutions based on 6 independent  subsamplings with $N=[7000,14000,28000,56000]$ on RBC and $N=[20000,40000,80000,160000,320000]$ on BSP. From the last row of Fig.~\ref{fig:linear_eqn}, we can still see the super-convergence phenomena for both extrinsic formulation and intrinsic formulation in the two examples. We have also examined such super-convergence phenomena using a second-order backward differentiation formula (BDF2) \cite{curtiss1952integration} for time-discretization.

\comment{
In addition, we also applied our GMLS approach to bumpy sphere examples which have the same parametrization form as in (\ref{eqn:psp}) but with the radial function $r(\theta,\phi)=1+0.05(1-\cos w_\theta \theta)\sin w_\phi \phi$ having perturbations from higher wavenumbers $w_\theta,w_\phi$.
However, we found that the convergence rate becomes slower than that predicted by the GMLS theory when $w_\theta=6,w_\phi=6$ and the Laplacian matrix even suffers from the eigenvalue instability when $w_\theta,w_\phi$ become even larger [not shown here].
Note that larger wavenumbers $w_\theta$ and $w_\phi$ give us a bumpy sphere with larger local curvatures.
From the above numerical results, we believe that the current GMLS approach with the specified weight function in  (\ref{kkernel}) can stabilize the vector Laplacians when manifolds are smooth with small curvatures. Nevertheless, when curvatures become large, the current GMLS approach may not provide a stable approximation to vector Laplacians.
}



\subsection{Viscous Burgers' equation}

\label{sec:burgers} In this subsection, we will show the results of solving  the
Burgers' equation on a torus,
\begin{equation}
\boldsymbol{u}_t  + \nabla_{\boldsymbol{u}} \boldsymbol{u} =  \nu \Delta_B \boldsymbol{u} + \boldsymbol{f}, \quad (\mathbf{x},t) \in M \times \left(0,T
\right] \text{ with } \boldsymbol{u}(\cdot,0) = \boldsymbol{u}_0(\cdot),  \label{eqn:burgers}
\end{equation}
where $\boldsymbol{u}(\cdot,t) \in \mathfrak{X}(M)$ represents a velocity field and $\nu=0.1$ is a
viscosity parameter. In this example, the parametrization of the torus is given in (\ref{eqn:torusr3}) and the points are sampled in the same way as those in subsection \ref{sec:eigen}. For the time discretization, we use the second-order Runge-Kutta method with $\Delta t = 0.0001$. We also examine the following result using Crank-Nicolson/Adams-Bashforth  scheme \cite{kim1985application}.


\begin{figure*}[htbp]
	\begin{minipage}[!b]{1\linewidth}
		\centering
		\includegraphics[width=4.7
		in, height=3 in]{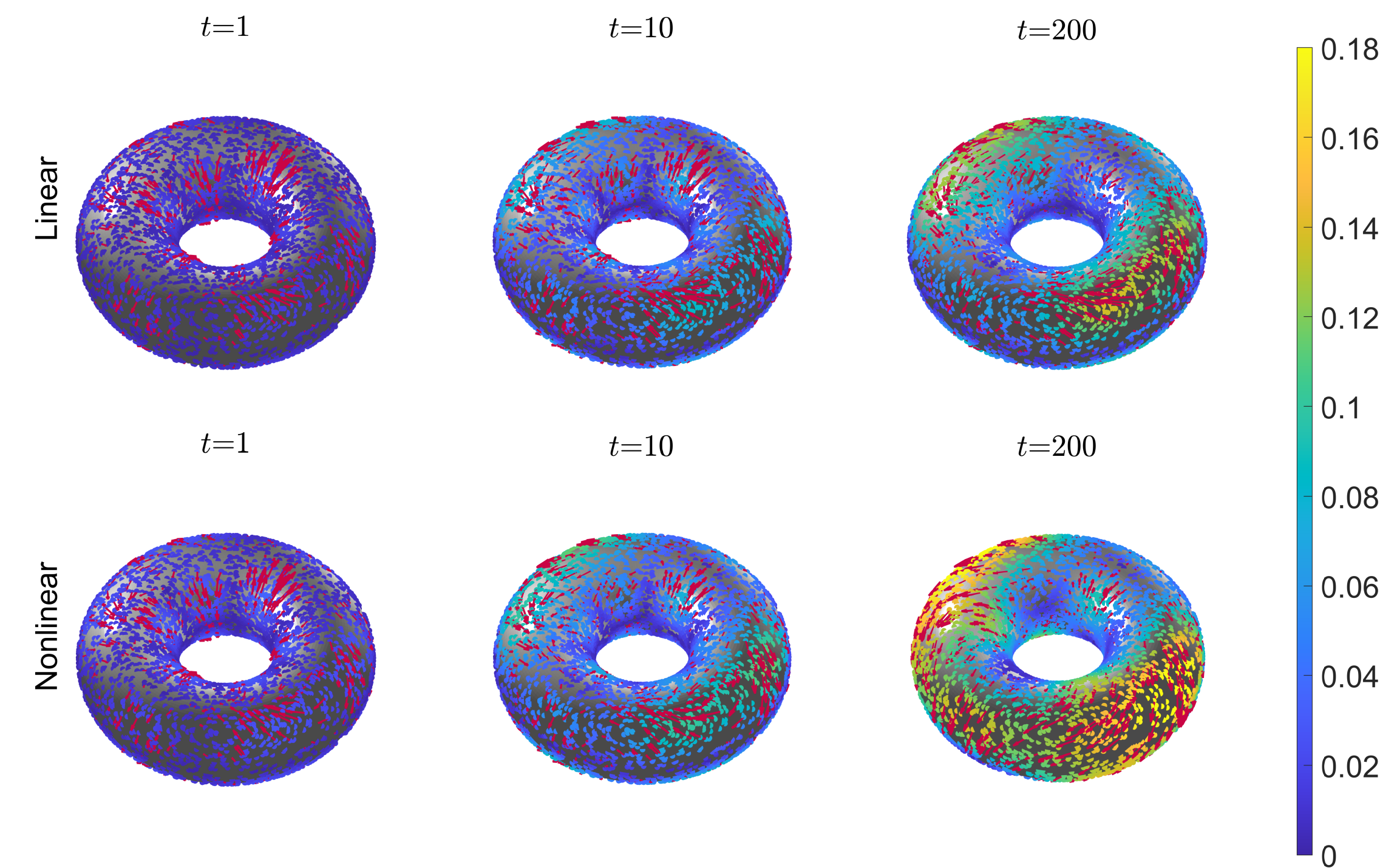}
	\end{minipage}
	\begin{minipage}[!b]{0.99\linewidth}
		\vspace{0.2in}
		\centering
		\begin{tabular}{cccc}
			{\small
				Linear (extr.) } & {\small Linear (intr.) } &
			{\small
				\ \ Nonlinear (extr.) } & {\small \ \ Nonlinear (intr.) }	
			\\
			\includegraphics[width=1.09
			in, height=0.95 in]{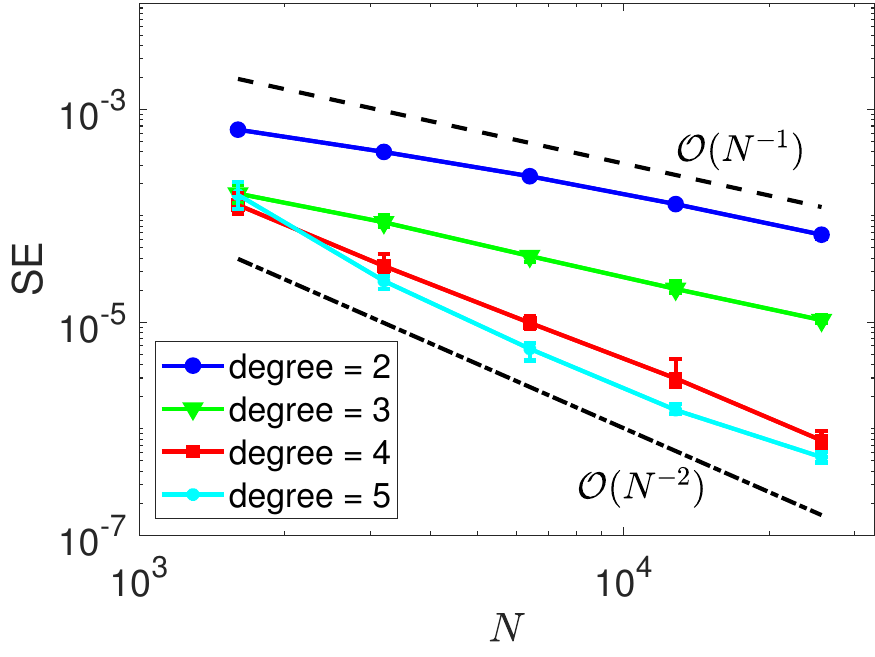} &
			\includegraphics[width=1.09
			in, height=0.95 in]{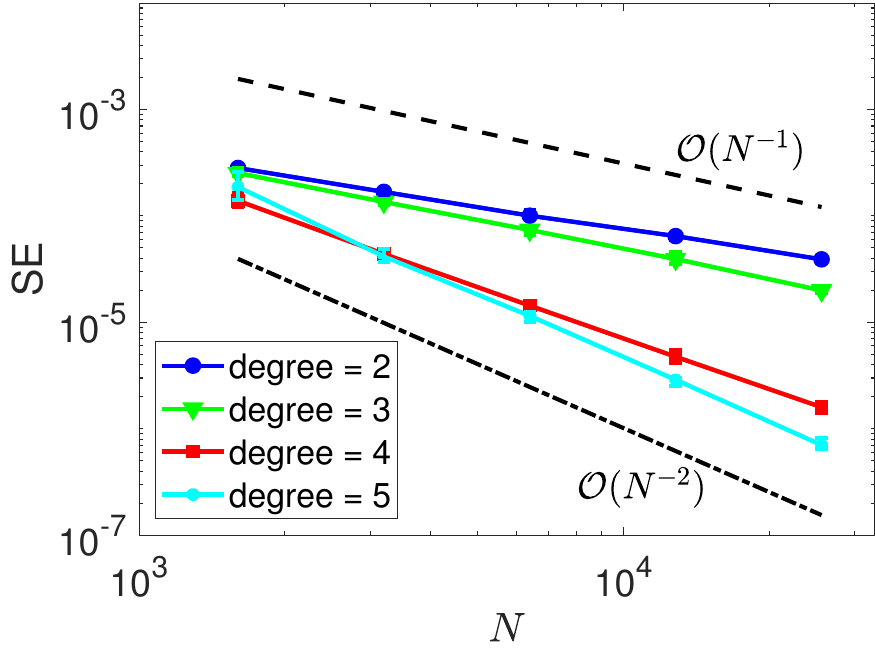} &
			\includegraphics[width=1.09
			in, height=0.95 in]{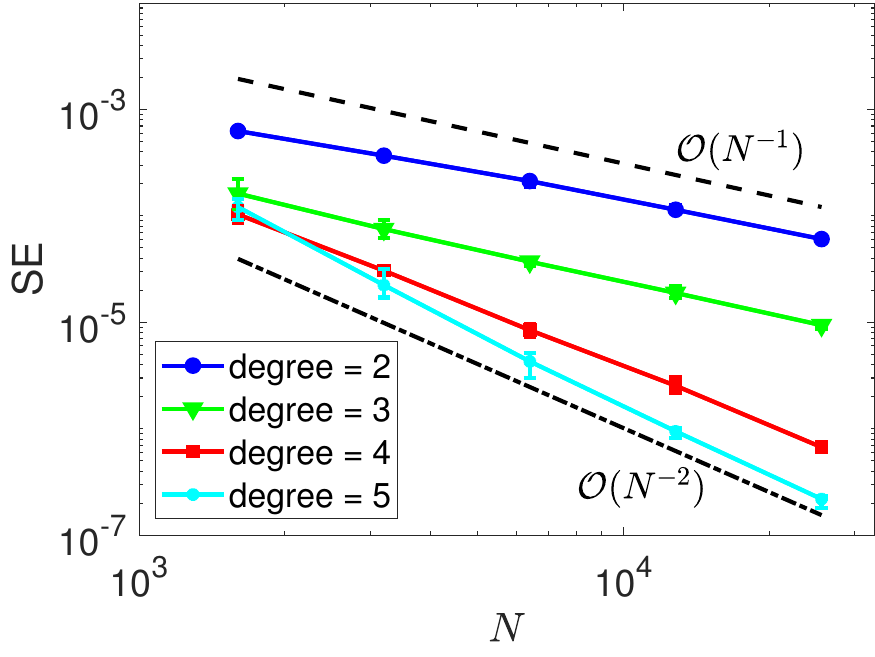} &
			\includegraphics[width=1.09
			in, height=0.95 in]{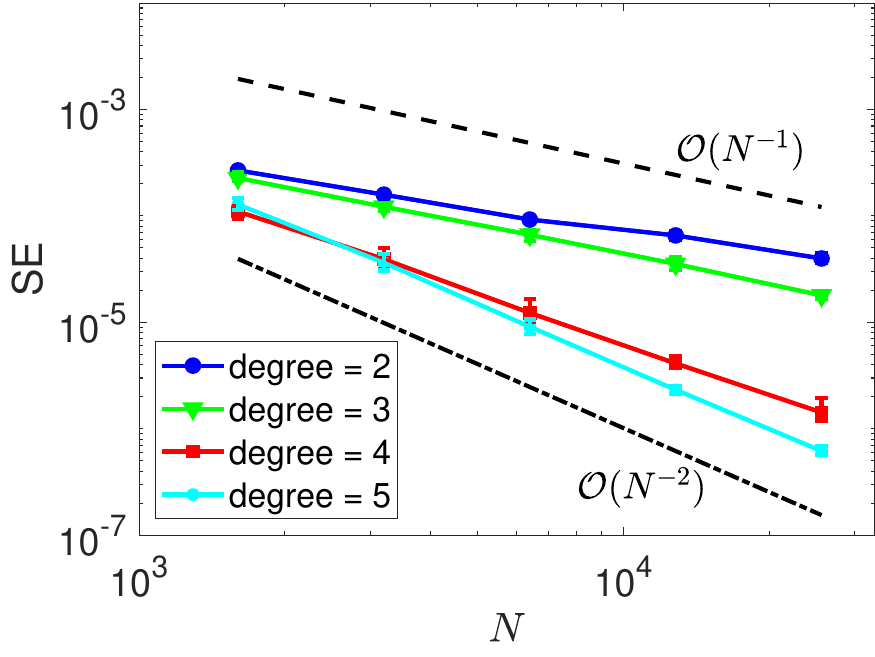}%
		\end{tabular}
	\end{minipage}
\caption{\textbf{Linear vector diffusion equation and nonlinear Burgers' equation on 2D torus in $\mathbb{R}^3$}. Time step $\Delta t=0.0001$
and $K$-nearest neighbor $K=70$. In the first two rows, we show the time evolution of the extrinsic GMLS solutions of linear (the first row) and nonlinear (the second row) equations using $N=6400$, degree $l=2$, with the color bar representing the magnitude of the vector fields $\hat{\mathbf{u}}$  and the red arrows indicating the direction of $\hat{\mathbf{u}}$. Shown in the third row is the convergence examination of the average \textbf{SE}s with error bars obtained from the standard deviations over 6 independent trials at $t=0.05$. }
\label{fig:burgers}
\end{figure*}

In the first two rows of Fig. \ref{fig:burgers}, we compare the evolution of the extrinsic GMLS solutions between the linear equation (without covariant derivative) in (\ref{eqn:linPDE}) and the nonlinear viscous Burgers' equation (with covariant derivative) in (\ref{eqn:burgers}). For both equations, we set the same initial condition $\boldsymbol{u}_0={0}$ to be a zero vector field and use the same external forcing $\boldsymbol{f}:= -\nu \Delta_B \boldsymbol{w}$ with $\boldsymbol{w}=0.05\sin 2 \theta\sin 3\phi\frac{\partial}{\partial\theta}+0.05\sin\theta\cos\phi\frac{\partial}{\partial\phi}$. In
our numerical implementation for the unknown manifold, we project the forcing onto the approximate tangent bundle by multiplying by $\mathbf{\hat{T}}^\top$ as the same in subsection~\ref{sec:poisson}. One can see from the
first two rows of Fig. \ref{fig:burgers} that the solutions exhibit similar dynamical behavior for a relatively short time ($t=10$) since the covariant derivative term has a small contribution for a short period.
However, the solutions behave quite differently at their steady states for a sufficiently long time ($t=200$).

In the third row of Fig. \ref{fig:burgers}, we verify the convergence of the proposed GMLS methods for solving the linear equation and the nonlinear Burgers' equation.  The true solution is set to be $\boldsymbol{u}=0.05\cos t\sin\theta\sin\phi\frac{\partial}{\partial\theta}$ $+0.05\cos t\sin\theta\cos\phi\frac{\partial}{\partial\phi}$ for both equations. One can see that the average \textbf{SE}s at time $t=0.05$ over 6 independent trials decrease at least on the order of $N^{-(l-1)/2}$ when degree $l=2,3,4,5$ as the GMLS theory predicts in (\ref{eq:gmlrt}). Here, the super-convergence can still be observed using the extrinsic method as seen in previous  Section~\ref{sec:poisson} and  Section~\ref{sec:lintim} while it can be barely seen using the intrinsic method.

\section{Conclusions and discussion}\label{sec:con}
In this paper, we extended the GMLS approach using both intrinsic formulation and extrinsic formulation for approximating vector Laplacians and covariant derivative
and further for solving vector-valued PDEs on unknown smooth manifolds, identified by a randomly-sampled data point cloud. The scheme can reach an arbitrary-order algebraic accuracy on smooth manifolds.
For the intrinsic method, we approximated the geometric quantities, functions, and their derivatives in the local Monge coordinate system to simplify the intrinsic formulation of the Bochner Laplacian in  (\ref{eqn:uab}) and (\ref{eqn:Bochdn}).
Then we transformed the representation of GMLS weights
from the  local tangent-vector basis to the  global tangent-vector basis in order to assemble these weights into
the sparse $dN\times dN $ Bochner Laplacian matrix.
For the extrinsic method,
we saved the memory cost and computational cost in approximating the Bochner Laplacian and covariant derivative, which improved our previous results in \cite{harlim2023radial}. In particular, we applied the GMLS approach as a sparse discretization and we also reduced the dimension of Laplacian matrices as well as the nonzero entries in the matrices by applying an appropriate coordinate transformation.

Several open questions remain for future investigation. For our GMLS method, the total complexity is $O(Nn)$, allowing it to accommodate relatively large ambient dimensions.  However, the requirements for storage, stability and convergence impose a strong restriction on the intrinsic dimension $d$.
For $d \leq 3$, as discussed in the 3D example in subsection \ref{sec:523}, an interpolation approach, in particular, Radial Basis Function-Finite Difference (RBF-FD) \cite{flyer2016role,jones2023generalized}, could potentially improve the stability and convergence.
Since the approximation space of the RBF-FD approach  contains both multivariate polynomials and radial basis functions, we expect it to outperform the GMLS approach,  whose approximation space consists solely of multivariate polynomials. For even larger intrinsic dimensions ($d\geq 4$), both GMLS and RBF-FD approaches may not be efficient for solving vector-valued PDE problems on manifolds due to the increased data requirements and memory storage demands. In such cases, neural networks  could be considered for discretization and approximation of PDE solutions, while the framework for approximating vector Laplacians developed in this work is still useful.
Additionally, an important direction for future work is solving vector-valued PDEs with classical boundary conditions on manifolds with boundaries.


















\section*{Acknowledgments}

S. J. was supported by the NSFC Grant No. 12471412 and the HPC Platform of ShanghaiTech University. Q. Y. was supported by the Simons
Foundation grant 601937, DNA. We would like to thank John Harlim for helpful comments.

\appendix

\section{Intrinsic formulation for other vector Laplacians}  \label{app:B}

In this section, we first review some concepts and notations from intrinsic differential geometry \cite{lee2018introduction,taylorPDE1}. Then we discuss the intrinsic GMLS approximation of the L Laplacian (which involves the symmetric stress tensor) and the Hodge Laplacian (which involves the anti-symmetric stress tensor) acting on vector fields. In addition, we also present the proof for Proposition~\ref{prop3p2}.

\subsection{Intrinsic Differential Geometry} \label{app:intr_DG}

Let $M$ be a $d$-dimensional Riemannian manifold with a local coordinate $%
(\xi ^{1},..,\xi ^{d})$ defined on an open subset $O$ in $M$. The Riemannian
metric is $\boldsymbol{g}=g_{ij}d\xi ^{i}\otimes d\xi ^{j}$ with the components $%
g_{ij}=\langle \partial _{i},\partial _{j}\rangle ,$ where $\partial _{i}=%
\frac{\partial }{\partial \xi ^{i}}$ is the $i$th coordinate vector field
and $(d\xi ^{1},...,d\xi ^{d})$ is the associated coordinate coframe. Recall
that $d\xi ^{j}(\partial _{i})=\delta _{i}^{j}$, where $\delta _{i}^{j}$ is
the delta function.

To define the vector Laplacian, we first introduce the concepts for the
gradient operator and the divergence operator. Let $\boldsymbol{u}=\tilde{u}^{i}\partial
_{i}\in \mathfrak{X}(M)$ be a vector field and let $\nabla $ be the
Levi-Civita connection. Then the covariant differentiation is defined as
$$
\nabla \boldsymbol{u}=\tilde{u}_{;j}^{i}\frac{\partial }{\partial \xi ^{i}}\otimes d\xi
^{j}
$$
where $\tilde{u}_{;j}^{i}=\frac{\partial \tilde{u}^{i}}{%
\partial \xi ^{j}}+\tilde{u}^{k}\Gamma _{jk}^{i}$, and the gradient of the vector field is defined to be a $(2,0)$ tensor field as in (\ref%
{eqn:gradu}),
$${\mathrm{grad}}_{g}\boldsymbol{u}={\sharp \nabla \boldsymbol{u}}=g^{kj}\tilde{u}%
_{;k}^{i}\frac{\partial }{\partial \xi ^{i}}\otimes \frac{\partial }{%
\partial \xi ^{j}}.
$$
For a vector field $\boldsymbol{u}$, its divergence can be expressed as 
$$\mathrm{div}_{g}\left( \boldsymbol{u}\right) =\text{tr}\left( \nabla \boldsymbol{u}\right) =\tilde{%
u}_{;i}^{i}.
$$
 Let $\boldsymbol{v}=\tilde{v}^{ij}\partial _{i}\otimes
\partial _{j}\in \mathfrak{X}(M)\times \mathfrak{X}(M)$ be a $(2,0)$ tensor
field. One can define the divergence of a $(2,0)$ tensor field (see e.g.  Ch 2 in \cite{taylorPDE1}):
\begin{equation}
\mathrm{div}_{g}(\boldsymbol{v}):=\mathrm{tr}_{1}^{2}\left( \nabla \boldsymbol{v}\right) =\mathrm{tr}_{1}^{2}\left( \tilde{%
v}_{;k}^{ij}\frac{\partial }{\partial \xi ^{i}}\otimes \frac{\partial }{%
\partial \xi ^{j}}\otimes d\xi ^{k}\right) =\tilde{v}_{;k}^{ik}\partial
_{i}, \label{eqdiv}
\end{equation}%
where the trace operator $\mathrm{tr}^2_1$ is taken on the \textit{last two indices} of
the tensor field $\nabla \boldsymbol{v}$.

We now present the intrinsic formulation for the three Laplacians acting on
vector fields: 1) the Bochner Laplacian, 2) the L Laplacian which involves
the symmetric surface stress tensor, and 3) the Hodge Laplacian. We concern
these three Laplacians here since they are the candidates for the diffusion
operators in Navier-Stokes equations on Riemannian manifolds as discussed in
references \cite{Chan2017,samavaki2020navier,taylorPDE3}.

\noindent 1) The Bochner Laplacian acting on a vector field $\boldsymbol{u}=\tilde{u}%
^{k}\partial _{k}$ is defined as
\begin{equation}
\Delta _{B}\boldsymbol{u}:=\mathrm{div}_{g}({\mathrm{grad}}_{g}\boldsymbol{u})=\mathrm{div}_{g}\left(
g^{kj}\tilde{u}_{;k}^{i}\frac{\partial }{\partial \xi ^{i}}\otimes \frac{%
\partial }{\partial \xi ^{j}}\right) =g^{kj}\tilde{u}_{;kj}^{i}\partial
_{i}. \label{eqn:Bog}
\end{equation}

\noindent 2) For the L laplacian, one can first define the symmetric
deformation rate tensor $S$ (cf. \cite{Chan2017,samavaki2020navier} or Ch
17 in \cite{taylorPDE3}),
\[
S:=\frac{1}{2}\left( {\mathrm{grad}}_{g}\boldsymbol{u}+({\mathrm{grad}}_{g}\boldsymbol{u})^{\top
}\right) =\frac{1}{2}\left( g^{mj}\tilde{u}_{;m}^{i}+g^{mi}\tilde{u}%
_{;m}^{j}\right) \partial _{i}\otimes \partial _{j}.
\]%
Then the (negative-definite) L Laplacian acting on the vector field $\boldsymbol{u}$ can
be defined and computed in the local coordinates as%
\begingroup\makeatletter\def\f@size{9}\check@mathfonts
\begin{eqnarray}
\Delta _{L}\boldsymbol{u}:=2\mathrm{div}_{g}\left( S\right) =\mathrm{div}_{g}\left( {%
\mathrm{grad}}_{g}\boldsymbol{u}+({\mathrm{grad}}_{g}\boldsymbol{u})^{\top }\right) =\left( g^{ij}%
\tilde{u}_{;ij}^{k}+g^{ki}\tilde{u}_{;ij}^{j}\right) \partial _{k}. \label{eqn:deL}
\end{eqnarray}
\endgroup

\noindent 3) The (negative-definite) Hodge Laplacian acting on the vector
field $\boldsymbol{u}=\tilde{u}^{k}\partial _{k}$ is defined as%
\[
\Delta _{H}\boldsymbol{u}:=-\sharp (\mathrm{d}^{\ast }\mathrm{d}+\mathrm{d}\mathrm{d}%
^{\ast })\flat \boldsymbol{u},
\]%
where $\mathrm{d}$ is the exterior derivative, $\mathrm{d}^{\ast }$ is the
formal adjoint of $\mathrm{d}$, and the symbols $\sharp $ and $\flat $ are
the standard musical isomorphisms. By this definition, the Hodge Laplacian
can be expressed and computed as%
\begin{eqnarray}
\Delta _{H}\boldsymbol{u} &=&\mathrm{div}_{g}({\mathrm{grad}}_{g}\boldsymbol{u}-({\mathrm{grad}}%
_{g}\boldsymbol{u})^{\top })+{\mathrm{grad}}_{g}\left( \mathrm{div}_{g}\left( \boldsymbol{u}\right)
\right)  \notag \\
&=&\left( g^{ij}\tilde{u}_{;ij}^{k}-g^{kj}\tilde{u}_{;ji}^{i}+g^{jk}\tilde{u}%
_{;ij}^{i}\right) \partial _{k}. \label{eqn:hog}
\end{eqnarray}

\subsection{Intrinsic GMLS approximation of other vector Laplacians}
In this section, we first provide the representation of $\Delta_L \boldsymbol{u}$ and $\Delta_H \boldsymbol{u}$ in the Monge coordinate system following the similar computational procedure for the Bochner $\Delta_B\boldsymbol{u}$ in Section~\ref{sec:mgc}. Recall that the notations, $\{a_{\alpha
_{1},\alpha _{2}}|2\leq \alpha _{1}+\alpha _{2}\leq l\}$, are the regression
coefficients of the parameterized surface (see equation (\ref{eqn:qhat})),
and the notations, $\beta _{\alpha _{1},\alpha _{2}}$ and $\gamma _{\alpha
_{1},\alpha _{2}}$ for $0\leq \alpha _{1}+\alpha _{2}\leq l$, are the
coefficients of respective component, $\tilde{u}^{1}$ and $\tilde{u}^{2}$,
of the vector field $\boldsymbol{u}=\tilde{u}^{k}\partial _{k}$ (see equation (\ref%
{eqn:ui})).

Then for $d=2$, the L Laplacian in (\ref{eqn:deL}) at the base point $\mathbf{x}_{0}$
can be approximated as
\begingroup\makeatletter\def\f@size{9}\check@mathfonts
\begin{eqnarray*}
&&\Delta _{L}\boldsymbol{u}(\mathbf{x}_{0}) \\
&=&(g^{ij}\tilde{u}_{;ij}^{k}+g^{ki}\tilde{u}_{;ij}^{j})\partial _{k}|_{%
\mathbf{x}_{0}}=\left[ \left( 2\tilde{u}_{;11}^{1}+\tilde{u}_{;22}^{1}+%
\tilde{u}_{;12}^{2}\right) \partial _{1}+\left( \tilde{u}_{;11}^{2}+\tilde{u}%
_{;21}^{1}+2\tilde{u}_{;22}^{2}\right) \partial _{2}\right] |_{\mathbf{x}%
_{0}}\\
&=&\left( 4\beta _{2,0}+2\beta _{0,2}+\beta _{0,0}\left[
8(a_{2,0})^{2}+(a_{1,1})^{2}+4a_{0,2}a_{2,0}\right] +\gamma _{1,1}+\gamma
_{0,0}\left[ 4a_{2,0}a_{1,1}+4a_{0,2}a_{1,1}\right] \right) \partial _{1}|_{\mathbf{x}_{0}} \\
&&+\left( \beta _{1,1}+\beta _{0,0}[4a_{2,0}a_{1,1}+4a_{0,2}a_{1,1}]+4\gamma
_{0,2}+2\gamma _{2,0}+\gamma _{0,0}\left[
8(a_{0,2})^{2}+(a_{1,1})^{2}+4a_{2,0}a_{2,0}\right] \right) \partial _{2}|_{\mathbf{x}_{0}},
\end{eqnarray*}%
\endgroup
and the Hodge Laplacian in (\ref{eqn:hog}) can be approximated as%
\begingroup\makeatletter\def\f@size{9}\check@mathfonts
\begin{eqnarray*}
&&\Delta _{H}\boldsymbol{u}(\mathbf{x}_{0})=(g^{ij}\tilde{u}_{;ij}^{k}-g^{kj}\tilde{u}%
_{;ji}^{i}+g^{jk}\tilde{u}_{;ij}^{i})\partial _{k}|_{\mathbf{x}_{0}} \\
&=&\left[ \left( \tilde{u}_{;22}^{1}-\tilde{u}_{;12}^{2}+\tilde{u}_{;11}^{1}+%
\tilde{u}_{;21}^{2}\right) \partial _{1}+\left( \tilde{u}_{;11}^{2}-\tilde{u}%
_{;21}^{1}+\tilde{u}_{;12}^{1}+\tilde{u}_{;22}^{2}\right) \partial _{2}%
\right] |_{\mathbf{x}_{0}} \\
&=&\left( 2\beta _{2,0}+2\beta _{0,2}+\beta _{0,0}\left[
4(a_{2,0})^{2}+2(a_{1,1})^{2}-4a_{0,2}a_{2,0}\right] +\gamma _{0,0}\left[
2a_{2,0}a_{1,1}+2a_{0,2}a_{1,1}\right] \right) \partial _{1}|_{\mathbf{x}_{0}} \\
&&+\left( \beta _{0,0}[2a_{2,0}a_{1,1}+2a_{0,2}a_{1,1}]+2\gamma
_{0,2}+2\gamma _{2,0}+\gamma _{0,0}\left[
4(a_{0,2})^{2}+2(a_{1,1})^{2}-4a_{2,0}a_{2,0}\right] \right) \partial _{2}|_{\mathbf{x}_{0}}.
\end{eqnarray*}
\endgroup

For the general $d$ dimensional manifold, one can obtain the L Laplacian%
\begin{equation*}
	\Delta _{L}\boldsymbol{u}(\mathbf{x}_{0})=\sum_{i=1}^d\left( \partial _{ii}%
	\tilde{u}^{k}+\partial _{ki}\tilde{u}^{i}+\sum_{s=1}^{n-d}\tilde{u}%
	^{l}\partial _{il}q^{s}\partial _{ki}q^{s}+\sum_{s=1}^{n-d}\tilde{u}%
	^{l}\partial _{kl}q^{s}\partial _{ii}q^{s}\right) \partial _{k}|_{\mathbf{x}%
		_{0}},
\end{equation*}
and the Hodge Laplacian
\begin{equation*}
	\Delta _{H}\boldsymbol{u}(\mathbf{x}_{0})=\sum_{i=1}^d\left( \partial _{ii}%
	\tilde{u}^{k}+2\sum_{s=1}^{n-d}\tilde{u}^{l}\partial _{ik}q^{s}\partial
	_{il}q^{s}-\sum_{s=1}^{n-d}\tilde{u}^{l}\partial _{ii}q^{s}\partial
	_{kl}q^{s}\right) \partial _{k}|_{\mathbf{x}_{0}}.
\end{equation*}%
Then one can follow the similar GMLS procedure in Sections~\ref{sec:ABL} and \ref{3p3highint} to approximate $\Delta_L \boldsymbol{u}$ and $\Delta_H \boldsymbol{u}$ on a point cloud.

\subsection{Proof for Proposition~\ref{prop3p2}}\label{sec:app1_3}

Here we prove Proposition~\ref{prop3p2}.

\begin{proof}
	From equation (\ref{eqn:xps}), we have%
	\begin{equation*}
		{\mathbf{x}}-{\mathbf{x}_{0}}=\sum_{s=1}^{n}p^{s}\mathbf{e}_{s},
	\end{equation*}%
	where $\left\{ \mathbf{e}_{s}\right\} _{s=1}^{n}$ are the standard
	orthonormal bases of the Euclidean space. By comparing with the Monge
	coordinate (\ref{eqn:xtq}), we arrive at
	\begin{equation*}
		\sum_{s=1}^{n}p^{s}\mathbf{e}_{s}=\sum_{i=1}^{d}\theta _{i}\boldsymbol{t}%
		_{i}({\mathbf{x}_{0}})+\sum_{s=1}^{n-d}q^{s}\boldsymbol{n}_{s}({\mathbf{x}%
			_{0}}),
	\end{equation*}%
	which implies that %
	\begin{equation}
		p^{s}=x^{s}-x_{0}^{s}=\sum_{i=1}^{d}\theta
		_{i}t_{i}^{s}+\sum_{j=1}^{n-d}q^{j}n_{j}^{s}. \label{eqn:pqs}
	\end{equation}%
	Here $\left\{ t_{i}^{s}\right\} _{i=1,\ldots ,d}^{s=1,\ldots ,n}$ is the $s$%
	th component of the $i$th tangent vector $\boldsymbol{t}_i$ and $\{n_{j}^{s}\}_{j=1,\ldots
		,n-d}^{s=1,\ldots ,n}$ is the $s$th component of the $j$th normal vector $\boldsymbol{n}_j$ \ in
	the Monge coordinate system. Then we conclude that%
	\begin{eqnarray*}
		\sum_{s=1}^{n}\partial _{ir}p^{s}\partial _{kl}p^{s}
		&=&\sum_{s=1}^{n}\partial _{ir}\left( \sum_{v=1}^{d}\theta
		_{v}t_{v}^{s}+\sum_{j=1}^{n-d}q^{j}n_{j}^{s}\right) \partial _{kl}\left(
		\sum_{h=1}^{d}\theta _{h}t_{h}^{s}+\sum_{m=1}^{n-d}q^{m}n_{m}^{s}\right) \\
		&=&\sum_{s=1}^{n}\sum_{j=1}^{n-d}\sum_{m=1}^{n-d}(\partial
		_{ir}q^{j})n_{j}^{s}(\partial
		_{kl}q^{m})n_{m}^{s}=\sum_{j=1}^{n-d}\sum_{m=1}^{n-d}(\partial
		_{ir}q^{j})(\partial _{kl}q^{m})\sum_{s=1}^{n}n_{j}^{s}n_{m}^{s} \\
		&=&\sum_{j=1}^{n-d}\sum_{m=1}^{n-d}(\partial _{ir}q^{j})(\partial
		_{kl}q^{m})\delta _{jm}=\sum_{s=1}^{n-d}\partial _{ir}q^{s}\partial
		_{kl}q^{s}.
	\end{eqnarray*}%
	Using the expresion (\ref{eqn:xps}) of $p^{s}$ and the quadratic form (\ref%
	{eqn:bthe}), we arrive at
	\begin{equation}
		\partial _{ir}p^{s}|_{{\mathbf{x}_{0}}}=2c_{ir}^{s},\text{ \ \ for all }%
		1\leq i,r\leq d\text{ and }1\leq s\leq n.  \label{eqn:cij}
	\end{equation}%
	Combining (\ref{eqn:derT}), (\ref{eqn:pqe}) and (\ref{eqn:cij}), we obtain
	that
	\begin{equation*}
		\frac{\partial \Gamma _{kl}^{i}}{\partial \theta _{r}}|_{\mathbf{x}%
			_{0}}=\sum_{s=1}^{n-d}\partial _{ir}q^{s}\partial _{kl}q^{s}|_{\mathbf{x}%
			_{0}}=\sum_{s=1}^{n}\partial _{ir}p^{s}\partial _{kl}p^{s}|_{\mathbf{x}%
			_{0}}=\sum_{s=1}^{n}4c_{ir}^{s}c_{kl}^{s}.
	\end{equation*}
\end{proof}

\section{Extrinsic formulation for other vector Laplacians \label{app:A}}
In this section, we first supplement some details for the  extrinsic formulation of the Bochner Laplacian as in \cite{harlim2023radial}. Then, we review the extrinsic formulation result for the L
Laplacian and the Hodge Laplacian acting on vector fields \cite{harlim2023radial}  and provide the extrinsic GMLS approximation of the  two vector Laplacians. In the following, we use the notations introduced in Section~\ref{sec:revext} and Section~\ref{sec:revdisH}.

\subsection{Bochner Laplacian}
\label{app:newA}
Here, we supplement a detailed derivation of the extrinsic formulation for the Bochner Laplacian  in Section~\ref{sec:revext}. Let $\boldsymbol{v}=v^{jk}\frac{\partial }{\partial \xi ^{j}}\otimes \frac{\partial }{%
\partial \xi ^{k}}\in \mathfrak{X}(M)\times \mathfrak{X}(M)$ be a $(2,0)$
tensor field, and let
\begin{equation*}
\boldsymbol{V}=\frac{\partial X^{s}}{\partial \xi ^{j}}v^{jk}\frac{\partial X^{t}}{%
\partial \xi ^{k}}\frac{\partial }{\partial X^{s}}\otimes \frac{\partial }{%
\partial X^{t}}:=V^{st}\frac{\partial }{\partial X^{s}}\otimes \frac{%
\partial }{\partial X^{t}},
\end{equation*}%
be a smooth extension of $\boldsymbol{v}$ onto an open subset in $\mathbb{R}^{3}$\ such
that $\boldsymbol{V}|_{M}=\boldsymbol{v}$. In general, the tensor field $\boldsymbol{v}$ is not symmetric, that is,
$v^{jk}\neq v^{kj}$. We now consider the extrinsic formulation for the
divergence of a $(2,0)$ tensor field $\boldsymbol{v}$. Using the definition in (\ref{eqdiv}), one can
verify the identity extension formulation for the divergence (see e.g. \cite{harlim2023radial}),
\comment{ The divergence is defined as $%
\mathrm{div}_{1}^{1}\left( v\right) =C_{1}^{1}(\nabla v)=v_{;j}^{jk}\frac{%
\partial }{\partial \xi ^{k}}$, where the super-index and sub-index in $%
\mathrm{div}_{s}^{r}\left( v\right) $\ denote the contraction for the $r$th
vector-field component (contravariant) and the $s$th $1$-form component (covariant) of $\nabla v$.}
\begingroup\makeatletter\def\f@size{9}\check@mathfonts
\begin{equation*}
\mathrm{div}_g\left( \boldsymbol{v}\right) =\mathbf{P}\mathrm{tr}^2_1\left(
\mathbf{P}\bar{\nabla}_{\mathbb{R}^{3}}\left( \boldsymbol{V}\right) \right) =\mathbf{P}%
\mathrm{tr}^2_1\left( \mathcal{G}_{r}V^{st}\frac{\partial }{\partial
X^{s}}\otimes \frac{\partial }{\partial X^{t}}\otimes \mathrm{d}X^{r}\right)
=\mathbf{P}\left( \mathcal{G}_{t}V^{st}\frac{\partial }{\partial X^{s}}%
\right) ,
\end{equation*}%
\endgroup
where $\mathcal{G}_{r}=\left( \mathbf{e}_{r}\cdot \mathbf{P}\right) \cdot
\bar{\nabla}_{\mathbb{R}^{3}}=\left( \mathbf{e}_{r}\cdot \mathbf{P}\right)
\cdot \overline{\mathrm{grad}}_{\mathbb{R}^{3}}$ is the differential
operator defined in (\ref{eqn:Gs}). In matrix form, above equation can be
written as%
\begin{eqnarray}
\mathrm{div}_g\left( \boldsymbol{v}\right) &=&\mathbf{P}\left[
\begin{array}{c}
\sum_{k=1}^{3}\mathcal{G}_{k}V^{1k} \\
\sum_{k=1}^{3}\mathcal{G}_{k}V^{2k} \\
\sum_{k=1}^{3}\mathcal{G}_{k}V^{3k}%
\end{array}%
\right] =\mathbf{P}\left[
\setlength{\arraycolsep}{0.3pt}
\begin{array}{ccc}
\mathcal{G}_{1} &  &  \\
& \mathcal{G}_{1} &  \\
&  & \mathcal{G}_{1}%
\end{array}%
\right] \left[
\setlength{\arraycolsep}{0.3pt}
\begin{array}{c}
V^{11} \\
V^{21} \\
V^{31}%
\end{array}%
\right] +\cdots +\mathbf{P}\left[
\setlength{\arraycolsep}{0.3pt}
\begin{array}{ccc}
\mathcal{G}_{3} &  &  \\
& \mathcal{G}_{3} &  \\
&  & \mathcal{G}_{3}%
\end{array}%
\right] \left[
\setlength{\arraycolsep}{0.3pt}
\begin{array}{c}
V^{13} \\
V^{23} \\
V^{33}%
\end{array}%
\right]  \notag \\
&=&\mathcal{H}_{1}\left[
\setlength{\arraycolsep}{0.3pt}
\begin{array}{c}
V^{11} \\
V^{21} \\
V^{31}%
\end{array}%
\right] +\mathcal{H}_{2}\left[
\setlength{\arraycolsep}{0.3pt}
\begin{array}{c}
V^{12} \\
V^{22} \\
V^{32}%
\end{array}%
\right] +\mathcal{H}_{3}\left[
\setlength{\arraycolsep}{0.3pt}
\begin{array}{c}
V^{13} \\
V^{23} \\
V^{33}%
\end{array}%
\right] :=\mathcal{H}_{1}\boldsymbol{V}\mathbf{|}_{1}+\mathcal{H}%
_{2}\boldsymbol{V}\mathbf{|}_{2} +\mathcal{H}_{3}\boldsymbol{V}\mathbf{|}_{3}  \label{eqn:divv}
\end{eqnarray}%
where $\mathcal{H}_{s}$ is defined as in (\ref{eqn:pV}) and $\boldsymbol{V}\mathbf{|}_s:=(V^{1s},V^{2s},V^{3s})^\top$\ is the $s$th column of tensor
field $\boldsymbol{V}$ for $s=1,2,3$.


Let $\boldsymbol{v}={\mathrm{grad}}_{g}\boldsymbol{u}\in \mathfrak{X}(M)\times \mathfrak{X}(M)$ be the
gradient of a vector field which is a $(2,0)$ tensor field. Using the
extrinsic formulation in ({\ref{eqn:gradu}}),\ one obtains the equation:
\begin{equation*}
\boldsymbol{V}=\left( \boldsymbol{V}\mathbf{|}_{1}, \boldsymbol{V}\mathbf{|}_{2}, \boldsymbol{V}\mathbf{|}_{3}\right)=\left( \mathcal{H}_{1}\boldsymbol{U},%
\mathcal{H}_{2}\boldsymbol{U},\mathcal{H}_{3}\boldsymbol{U}\right)=\mathbf{P}\left(
\overline{{\mathrm{grad}}}_{\mathbb{R}^{3}}\boldsymbol{U}\right) \mathbf{P},
\end{equation*}%
which gives us
\begin{equation}
\boldsymbol{V}\mathbf{|}_{k}=\mathcal{H}_{k}\boldsymbol{U},\text{ \ \ for }k=1,2,3%
\text{.}  \label{eqn:pV2}
\end{equation}%
Combining (\ref{eqn:divv}) and (\ref{eqn:pV2}), one can obtain the extrinsic
formulation for the Bochner Laplacian of a vector field, $\Delta
_{B}\boldsymbol{u}=g^{ij}u_{;ij}^{k}\frac{\partial }{\partial \xi ^{k}}=\mathrm{div}_g({\mathrm{grad}}_{g}{\boldsymbol{u})}$, as
\begin{eqnarray*}
\Delta _{B}\boldsymbol{u} &=&\mathrm{div}_g({\mathrm{grad}}_{g}\boldsymbol{u})=\mathrm{div}_g\left( \boldsymbol{v}\right) =\mathcal{H}_{1}\boldsymbol{V}\mathbf{|}_{1}+%
\mathcal{H}_{2}\boldsymbol{V}\mathbf{|}_{2}+\mathcal{H}_{3}\boldsymbol{V}\mathbf{|}_{3} \\
&=&\mathcal{H}_{1}\mathcal{H}_{1}\boldsymbol{U}+\mathcal{H}_{2}\mathcal{H}_{2}\boldsymbol{U}+\mathcal{H%
}_{3}\mathcal{H}_{3}\boldsymbol{U}:=\bar{\Delta}_{B}\boldsymbol{U}.
\end{eqnarray*}%

\subsection{L Laplacian}\label{sec:inLlap}

We first review the extrinsic formulation for the L Laplacian (see e.g. \cite%
{suchde2021meshfree,harlim2023radial}). Let $\boldsymbol{v}={\mathrm{grad}}_{g}\boldsymbol{u}$\ be the
gradient of a vector field $\boldsymbol{u}$ and let $\boldsymbol{V}$ be a smooth extension of $\boldsymbol{v}$ onto
an open subset in $\mathbb{R}^{3}$\ such that $\boldsymbol{V}|_{M}=\boldsymbol{v}$. Based on the
definition of the L Laplacian in (\ref{eqn:deL}), one can obtain the
extrinsic formulation for the L Laplacian:
\begin{eqnarray}
\Delta _{L}\boldsymbol{u} &=&\mathrm{div}_g\left( {\boldsymbol{v}+\boldsymbol{v}}^{\top }\right) =\mathbf{P}%
\mathrm{tr}^2_1\left( \mathbf{P}\bar{\nabla}_{\mathbb{R}^{3}}(\boldsymbol{V}+\boldsymbol{V}^{\top
})\right) =\mathbf{P}\left[
\setlength{\arraycolsep}{0.3pt}
\begin{array}{c}
\sum_{k=1}^{3}\mathcal{G}_{k}\left( V^{1k}+V^{k1}\right) \\
\sum_{k=1}^{3}\mathcal{G}_{k}\left( V^{2k}+V^{k2}\right) \\
\sum_{k=1}^{3}\mathcal{G}_{k}\left( V^{3k}+V^{k3}\right)%
\end{array}%
\right]  \notag \\
&=&\mathcal{H}_{1}\left[\boldsymbol{V}\mathbf{|}_{1}+ \left(\bar{\boldsymbol{V}}_{1}\right)^\top
\right] +\mathcal{H}_{2}\left[ \boldsymbol{V}\mathbf{|}_2 +\left(\bar{\boldsymbol{V}}_{2}\right) ^{\top }
\right] +\mathcal{H}_{3}\left[\boldsymbol{V}
\mathbf{|}_{3} +\left(\bar{\boldsymbol{V}}_{3}\right) ^{\top }\right] ,  \label{eqn:divL}
\end{eqnarray}%
where $\bar{\boldsymbol{V}}_{i}=\left[ V_{i1},V_{i2},V_{i3}\right] $ is
the $i$th row of the $(2,0)$ tensor field $\boldsymbol{V}$ for $i=1,2,3$. Here, we
have used the notations $\mathcal{G}_{j}$ in (\ref{eqn:Gs}) and $\mathcal{H}%
_{j}$ in (\ref{eqn:pV}) for $j=1,2,3$, and we can further calculate $\bar{\boldsymbol{V}}^{i}$ as%
\begin{eqnarray}
(\bar{\boldsymbol{V}}^{i})^\top &=&\left[
V_{i1},V_{i2},V_{i3}\right]^\top
=\left(\left[P_{i1},P_{i2},P_{i3}\right]\left[\begin{matrix}
\mathcal{G}_{1}U^{1} & \mathcal{G}_{2}U^{1} & \mathcal{G}_{3}U^{1} \\
\mathcal{G}_{1}U^{2} & \mathcal{G}_{2}U^{2} & \mathcal{G}_{3}U^{2} \\
\mathcal{G}_{1}U^{3} & \mathcal{G}_{2}U^{3} & \mathcal{G}_{3}U^{3}%
\end{matrix}\right]\right)^\top  \notag \\
&=&\left[
\setlength{\arraycolsep}{2.3pt}
\begin{array}{ccc}
{P}_{i1}\mathcal{G}_{1} & {P}_{i2}\mathcal{G}_{1} & {P}_{i3}\mathcal{G}_{1}
\\
{P}_{i1}\mathcal{G}_{2} & {P}_{i2}\mathcal{G}_{2} & {P}_{i3}\mathcal{G}_{2}
\\
{P}_{i1}\mathcal{G}_{3} & {P}_{i2}\mathcal{G}_{3} & {P}_{i3}\mathcal{G}_{3}%
\end{array}%
\right] \left[
\setlength{\arraycolsep}{0.3pt}
\begin{array}{c}
U^{1} \\
U^{2} \\
U^{3}%
\end{array}%
\right] :=\mathcal{S}_{i}\boldsymbol{U}.  \label{eqn:Vcol}
\end{eqnarray}%
Substituting (\ref{eqn:Vcol}) into (\ref{eqn:divL}), we obtain%
\begin{equation*}
\Delta _{L}\boldsymbol{u}=\mathcal{H}_{1}\left[ \mathcal{H}_{1}+\mathcal{S}_{1}\right]
\boldsymbol{U}+\cdots +\mathcal{H}_{3}\left[ \mathcal{H}_{3}+\mathcal{S}_{3}\right] \boldsymbol{U}:=%
\bar{\Delta}_{L}\boldsymbol{U}.
\end{equation*}%
For the discretization, one can approximate the differential operator $%
\mathcal{S}_{i}$ in the stencil as,%
\begin{equation}
\mathbf{S}_{i}=\left[
\begin{array}{ccc}
\mathrm{diag}\left( \mathbf{p}_{i1}\right) \mathbf{G}_{1} & \mathrm{diag}%
\left( \mathbf{p}_{i2}\right) \mathbf{G}_{1} & \mathrm{diag}\left( \mathbf{p}%
_{i3}\right) \mathbf{G}_{1} \\
\mathrm{diag}\left( \mathbf{p}_{i1}\right) \mathbf{G}_{2} & \mathrm{diag}%
\left( \mathbf{p}_{i2}\right) \mathbf{G}_{2} & \mathrm{diag}\left( \mathbf{p}%
_{i3}\right) \mathbf{G}_{2} \\
\mathrm{diag}\left( \mathbf{p}_{i1}\right) \mathbf{G}_{3} & \mathrm{diag}%
\left( \mathbf{p}_{i2}\right) \mathbf{G}_{3} & \mathrm{diag}\left( \mathbf{p}%
_{i3}\right) \mathbf{G}_{3}%
\end{array}%
\right] _{3K\times 3K},  \label{eqn:Si}
\end{equation}%
where $\mathbf{p}_{ij}=\left( {P}_{ij}(\mathbf{x}_{0,1}),\ldots ,{P}_{ij}(%
\mathbf{x}_{0,K})\right) \in \mathbb{R}^{K\times 1}$ for $1\leq i,j\leq 3$.
Then, the L Laplacian can be subsequently approximated in the stencil $S_{{%
\mathbf{x}_{0}}}$ as,
\begin{equation}
(\Delta _{L}\boldsymbol{u})|_{S_{{\mathbf{x}_{0}}}}=(\bar{\Delta}_{L}\boldsymbol{U})|_{S_{{\mathbf{x}_{0}}%
}}=\sum_{\ell =1}^{3}(\mathcal{H}_{\ell }\left( \mathcal{H}_{\ell }+\mathcal{S%
}_{\ell }\right) \boldsymbol{U})|_{S_{{\mathbf{x}_{0}}}}\approx \sum_{\ell =1}^{3}\mathbf{H%
}_{\ell }\left( \mathbf{H}_{\ell }+\mathbf{S}_{\ell }\right) \mathbf{U}_{{%
\mathbf{x}_{0}}},  \label{eqn:HSl}
\end{equation}%
where $\mathbf{U}_{{\mathbf{x}_{0}}}$ has been defined in (\ref{eqn:Sx0}).


We now consider the reduction for the L Laplacian following the similar
procedure to the Bochner Laplacian as in Section \ref{sec:dimred}. Using the
definition of $\mathbf{H}_{i}$ in (\ref{eqn:Hi}) and $\mathbf{S}_{i}$ in (%
\ref{eqn:Si}), one can write the L Laplacian in the stencil $S_{{\mathbf{x}%
_{0}}}$ as
\begin{equation}
(\Delta _{L}\boldsymbol{u})|_{S_{{\mathbf{x}_{0}}}}\approx \sum_{\ell =1}^{3}\mathbf{H}%
_{\ell }\left( \mathbf{H}_{\ell }+\mathbf{S}_{\ell }\right) \mathbf{U}_{{%
\mathbf{x}_{0}}}=\sum_{\ell =1}^{3}\mathbf{P}^{\otimes }\mathbf{H}_{\ell }%
\mathbf{P}^{\otimes }\left( \mathbf{H}_{\ell }+\mathbf{S}_{\ell }\right)
\mathbf{P}^{\otimes }\mathbf{U}_{{\mathbf{x}_{0}}},  \notag
\end{equation}%
where made use of $\mathbf{P}^{\otimes }\mathbf{S}_{\ell }=\mathbf{S}_{\ell
} $ since each column of $\mathbf{S}_{\ell }$ always lives in the tangent
bundle of $M$. Using the decomposition $\mathbf{P}^{\otimes }=\mathbf{T}%
^{\otimes }\mathbf{T}^{\otimes \top }$, one arrives at%
\begin{eqnarray}
\Delta _{L}\boldsymbol{u}|_{S_{{\mathbf{x}_{0}}}} &\approx &\sum_{\ell =1}^{3}\mathbf{T}%
^{\otimes }(\mathbf{T}^{\otimes \top }\mathbf{H}_{\ell }\mathbf{T}^{\otimes
})[\mathbf{T}^{\otimes \top }\left( \mathbf{H}_{\ell }+\mathbf{S}_{\ell
}\right) \mathbf{T}^{\otimes }](\mathbf{T}^{\otimes \top }\mathbf{U}_{{%
\mathbf{x}_{0}}})  \notag \\
&:=& \mathbf{T}^{\otimes }[\sum_{\ell =1}^{3}\mathbf{R}_{\ell }\left(
\mathbf{R}_{\ell }+\mathbf{M}_{\ell }\right) ](\mathbf{T}^{\otimes \top }%
\mathbf{U}_{{\mathbf{x}_{0}}}),  \label{eqn:rmL}
\end{eqnarray}%
where $\mathbf{R}_{\ell }:=\mathbf{T}^{\otimes \top }\mathbf{H}_{\ell }%
\mathbf{T}^{\otimes }=\mathbf{T}^{\otimes \top }\left( \mathbf{I}_{3}\mathbf{%
\otimes G}_{\ell }\right) \mathbf{T}^{\otimes }$ has been defined in (\ref%
{eqn:TTB}), $\mathbf{M}_{\ell }:=\mathbf{T}^{\otimes \top }\mathbf{S}_{\ell }%
\mathbf{T}^{\otimes }$ is the reduction of $\mathbf{S}_{\ell }$, and $%
\sum_{\ell =1}^{3}\mathbf{R}_{\ell }\left( \mathbf{R}_{\ell }+\mathbf{M}%
_{\ell }\right) \in \mathbb{R}^{2K\times 2K}$\ is the L Laplacian matrix
that maps the coordinate representation of $\boldsymbol{u}|_{S_{{\mathbf{x}_{0}}}}$\ to
the coordinate representation of $(\Delta _{L}\boldsymbol{u})|_{S_{{\mathbf{x}_{0}}}}$
w.r.t. the global basis $\{\boldsymbol{t}_{j}(\mathbf{x}_{0,k})%
\}_{j=1,2}^{k=1,\ldots ,K}$.

Next, we compute the differential matrix $\mathbf{M}_{\ell }$
componentwisely. Fix some $\ell $, one can compute a $2\times 2$\ block
matrix component of $\mathbf{M}_{\ell }$ whose rows corresponding to $%
\mathbf{x}_{0,s}$ and columns corresponding to $\mathbf{x}_{0,t}$,
\begingroup\makeatletter\def\f@size{9}\check@mathfonts
\begin{eqnarray*}
&&\left[
\setlength{\arraycolsep}{0.3pt}
\begin{array}{c}
\boldsymbol{t}_{1}(\mathbf{x}_{0,s})^{\top } \\
\boldsymbol{t}_{2}(\mathbf{x}_{0,s})^{\top }%
\end{array}%
\right] _{2\times 3}\left[
\setlength{\arraycolsep}{2.5pt}
\begin{array}{ccc}
P_{\ell1 }(\mathbf{x}_{0,s})\left( \mathbf{G}_{1}\right) _{st} & P_{\ell 2 }(%
\mathbf{x}_{0,s})\left( \mathbf{G}_{1}\right) _{st} & P_{\ell 3 }(\mathbf{x}%
_{0,s})\left( \mathbf{G}_{1}\right) _{st} \\
P_{\ell 1 }(\mathbf{x}_{0,s})\left( \mathbf{G}_{2}\right) _{st} & P_{\ell 2 }(%
\mathbf{x}_{0,s})\left( \mathbf{G}_{2}\right) _{st} & P_{\ell 3 }(\mathbf{x}%
_{0,s})\left( \mathbf{G}_{2}\right) _{st} \\
P_{\ell1 }(\mathbf{x}_{0,s})\left( \mathbf{G}_{3}\right) _{st} & P_{\ell 2 }(%
\mathbf{x}_{0,s})\left( \mathbf{G}_{3}\right) _{st} & P_{\ell 3 }(\mathbf{x}%
_{0,s})\left( \mathbf{G}_{3}\right) _{st}%
\end{array}%
\right] 
\left[
\setlength{\arraycolsep}{0.5pt}
\begin{array}{cc}
\boldsymbol{t}_{1}(\mathbf{x}_{0,t}) & \boldsymbol{t}_{2}(\mathbf{x}_{0,t})%
\end{array}%
\right] _{3\times 2} \\
&=&\left( \left[
\setlength{\arraycolsep}{0.3pt}
\begin{array}{c}
\boldsymbol{t}_{1}(\mathbf{x}_{0,s})^{\top } \\
\boldsymbol{t}_{2}(\mathbf{x}_{0,s})^{\top }%
\end{array}%
\right] \left[
\setlength{\arraycolsep}{0.3pt}
\begin{array}{c}
\left( \mathbf{G}_{1}\right) _{st} \\
\left( \mathbf{G}_{2}\right) _{st} \\
\left( \mathbf{G}_{3}\right) _{st}%
\end{array}%
\right] \right) 
\Big( \left[
\setlength{\arraycolsep}{2.5pt}
\begin{array}{ccc}
P_{\ell  1}(\mathbf{x}_{0,s}) & P_{\ell 2 }(\mathbf{x}_{0,s}) & P_{\ell 3}(%
\mathbf{x}_{0,s})%
\end{array}%
\right] 
\left[
\setlength{\arraycolsep}{2.5pt}
\begin{array}{cc}
\boldsymbol{t}_{1}(\mathbf{x}_{0,t}) & \boldsymbol{t}_{2}(\mathbf{x}_{0,t})%
\end{array}%
\right] \Big) :=\mathbf{\hat{M}}_{st}^{\ell }.
\end{eqnarray*}%
\endgroup
Stacking all these block matrices $\{\mathbf{\hat{M}}_{st}^{\ell
}\}_{s,t=1}^{K}$\ for each $\ell =1,2,3$,\ one can obtain a reshapement of $%
\mathbf{M}_{\ell }$ via some appropriate row and column interchanges.
Therefore, the weights $\{\boldsymbol{w}^{k}\in \mathbb{R}^{2\times
2}\}_{k=1}^{K}$\ for approximating the L Laplacian at the base $\mathbf{x}%
_{0}$ corresponding to the $2$ rows of the matrix $\sum_{\ell =1}^{3}\mathbf{%
R}_{\ell }\left( \mathbf{R}_{\ell }+\mathbf{M}_{\ell }\right) $ in (\ref%
{eqn:rmL}) become%
\begin{equation*}
(\boldsymbol{w}^{1},\ldots ,\boldsymbol{w}^{K})=\sum_{\ell =1}^{3}\left(
\mathbf{\hat{R}}_{11}^{\ell },\ldots ,\mathbf{\hat{R}}_{1K}^{\ell }\right) %
\left[
\setlength{\arraycolsep}{2.0pt}
\begin{array}{ccc}
\mathbf{\hat{R}}_{11}^{\ell }+\mathbf{\hat{M}}_{11}^{\ell } & \cdots &
\mathbf{\hat{R}}_{1K}^{\ell }+\mathbf{\hat{M}}_{1K}^{\ell } \\
\vdots & \ddots & \vdots \\
\mathbf{\hat{R}}_{K1}^{\ell }+\mathbf{\hat{M}}_{K1}^{\ell } & \cdots &
\mathbf{\hat{R}}_{KK}^{\ell }+\mathbf{\hat{M}}_{KK}^{\ell }%
\end{array}%
\right] \in \mathbb{R}^{2\times 2K}.
\end{equation*}%
Here, these weights satisfy $\Delta _{L}\boldsymbol{u}({\mathbf{x}_{0}})\approx
\sum_{k=1}^{K}\boldsymbol{w}^{k}\boldsymbol{u}|_{{\mathbf{x}_{0,k}}}$ at the base point $%
\mathbf{x}_{0}$.


\subsection{Hodge Laplacian}\label{sec:apphodge}

We now review the extrinsic formulation for the Hodge Laplacian \cite%
{harlim2023radial}, which can be written as%
\begin{equation}
\Delta _{H}\boldsymbol{u}=\mathbf{P}\mathrm{tr}^2_1\left( \mathbf{P}\bar{\nabla}_{%
\mathbb{R}^{3}}(\boldsymbol{V}-\boldsymbol{V}^{\top })\right) +\mathbf{P}\overline{{\mathrm{grad}}}_{%
\mathbb{R}^{3}}(\mathrm{tr}^1_1[\mathbf{P}\bar{\nabla}_{\mathbb{R}%
^{3}}\boldsymbol{U}]),  \label{eqn:Hodu}
\end{equation}%
where $\boldsymbol{V}=\mathbf{P}\left( \overline{{\mathrm{grad}}}_{\mathbb{R}^{3}}\boldsymbol{U}\right) \mathbf{P}$ is an extension of $\boldsymbol{v}={\mathrm{grad}}_{g}\boldsymbol{u}$ such that $%
\boldsymbol{V}|_{M}=\boldsymbol{v}$. We note that the first term in above (\ref{eqn:Hodu}) is only
different from the extrinsic formulation for the L Laplacian in (\ref%
{eqn:divL}) by a negative sign. Thus, it only remains to compute the last
term in (\ref{eqn:Hodu}),%
\begin{equation}
\mathbf{P}\overline{{\mathrm{grad}}}_{\mathbb{R}^{3}}(\mathrm{tr}^1_1[%
\mathbf{P}\bar{\nabla}_{\mathbb{R}^{3}}\boldsymbol{U}])=\mathbf{P}\overline{{\mathrm{grad}%
}}_{\mathbb{R}^{3}}(\sum_{k=1}^{3}\mathcal{G}_{k}U^{k})=\left[
\setlength{\arraycolsep}{0.3pt}
\begin{array}{c}
\mathcal{G}_{1} \\
\mathcal{G}_{2} \\
\mathcal{G}_{3}%
\end{array}%
\right] \left[ \mathcal{G}_{1},\mathcal{G}_{2},\mathcal{G}_{3}\right] \left[
\setlength{\arraycolsep}{0.3pt}
\begin{array}{c}
U^{1} \\
U^{2} \\
U^{3}%
\end{array}%
\right] .  \label{eqn:ggu2}
\end{equation}%
For the discretization, one only needs to substitute each $\mathcal{G}_{i}$\
above with the differential matrix $\mathbf{G}_{i}$ for $i=1,2,3$.

Next, we only need to consider the dimension reduction for the
discretization of equation (\ref{eqn:ggu2}). Using the decomposition $%
\mathbf{P}^{\otimes }=\mathbf{T}^{\otimes }\mathbf{T}^{\otimes \top }$, we
arrive at%
\begin{eqnarray*}
\mathbf{P}\overline{{\mathrm{grad}}}_{\mathbb{R}^{3}}\left( \mathrm{tr}^1_1\left[ \mathbf{P}\bar{\nabla}_{\mathbb{R}^{3}}\boldsymbol{U}\right] \right) |_{S_{%
\mathbf{x}_{0}}} &\approx &\mathbf{T}^{\otimes }\left( \mathbf{T}^{\otimes
\top }\left[
\setlength{\arraycolsep}{0.3pt}
\begin{array}{c}
\mathbf{G}_{1} \\
\mathbf{G}_{2} \\
\mathbf{G}_{3}%
\end{array}%
\right] \right) \Big( \big[\mathbf{G}_{1},\mathbf{G}_{2},\mathbf{G}_{3}\big]%
\mathbf{T}^{\otimes }\Big) \Big( \mathbf{T}^{\otimes \top }\mathbf{U}_{{%
\mathbf{x}_{0}}}\Big)  \\
&:= &\mathbf{T}^{\otimes }\mathbf{J}\left( \mathbf{T}^{\otimes \top }\mathbf{%
U}_{{\mathbf{x}_{0}}}\right) .
\end{eqnarray*}%
Thus, the Hodge Laplacian can be approximated as
\begin{equation*}
(\Delta _{H}\boldsymbol{u})|_{S_{{\mathbf{x}_{0}}}}\approx \mathbf{T}^{\otimes }\left(
\Big[\sum_{\ell =1}^{3}\mathbf{R}_{\ell }\left( \mathbf{R}_{\ell }-\mathbf{M}%
_{\ell }\right) \Big]+\mathbf{J}\right) (\mathbf{T}^{\otimes \top }\mathbf{U}_{{%
\mathbf{x}_{0}}}).
\end{equation*}%
Let $\left( \mathbf{\hat{J}}_{11},\ldots ,\mathbf{\hat{J}}_{1K}\right) \in
\mathbb{R}^{2\times 2K}$ be the $2$ rows of $\mathbf{J}$ corresponding to
the approximation of $\Delta _{H}u$\ at the base point $\mathbf{x}_{0}$.
Then, the weights $\{\boldsymbol{w}^{k}\in \mathbb{R}^{2\times 2}\}_{k=1}^{K}
$\ for the Hodge Laplacian become%
\begingroup\makeatletter\def\f@size{9}\check@mathfonts
\begin{equation*}
(\boldsymbol{w}^{1},\ldots ,\boldsymbol{w}^{K})=\sum_{\ell =1}^{3}\left(
\mathbf{\hat{R}}_{11}^{\ell },\ldots ,\mathbf{\hat{R}}_{1K}^{\ell }\right) %
\left[
\setlength{\arraycolsep}{0.3pt}
\begin{array}{ccc}
\mathbf{\hat{R}}_{11}^{\ell }-\mathbf{\hat{M}}_{11}^{\ell } & \cdots  &
\mathbf{\hat{R}}_{1K}^{\ell }-\mathbf{\hat{M}}_{1K}^{\ell } \\
\vdots  & \ddots  & \vdots  \\
\mathbf{\hat{R}}_{K1}^{\ell }-\mathbf{\hat{M}}_{K1}^{\ell } & \cdots  &
\mathbf{\hat{R}}_{KK}^{\ell }-\mathbf{\hat{M}}_{KK}^{\ell }%
\end{array}%
\right] +\left( \mathbf{\hat{J}}_{11},\ldots ,\mathbf{\hat{J}}_{1K}\right),
\end{equation*}%
\endgroup
which are the $2$ rows of the matrix $\sum_{\ell =1}^{3}\mathbf{R}_{\ell
}\left( \mathbf{R}_{\ell }-\mathbf{M}_{\ell }\right) +\mathbf{J}$
corresponding to the base point $\mathbf{x}_{0}$ such that $\Delta _{H}\boldsymbol{u}({%
\mathbf{x}_{0}})\approx \sum_{k=1}^{K}\boldsymbol{w}^{k}\boldsymbol{u}|_{{\mathbf{x}_{0,k}}}
$.

\bibliographystyle{abbrv}
\bibliography{references}

\begin{thebibliography}{10}

\bibitem{alvarez2021local}
D.~{\'A}lvarez, P.~Gonz{\'a}lez-Rodr{\'\i}guez, and M.~Kindelan.
\newblock A local radial basis function method for the laplace--beltrami
  operator.
\newblock {\em J. Sci. Comput.}, 86, 2021.

\bibitem{arnold2006finite}
D.~N. Arnold, R.~S. Falk, and R.~Winther.
\newblock Finite element exterior calculus, homological techniques, and
  applications.
\newblock {\em Acta Numer.}, 15:1--155, 2006.

\bibitem{2015Discrete}
O.~Azencot, M.~Ovsjanikov, F.~Chazal, and M.~Ben-Chen.
\newblock Discrete derivatives of vector fields on surfaces - an operator
  approach.
\newblock {\em ACM Trans. Graph.}, 34(3):1--13, 2015.

\bibitem{barrett2015numerical}
J.~W. Barrett, H.~Garcke, and R.~N{\"u}rnberg.
\newblock Numerical computations of the dynamics of fluidic membranes and
  vesicles.
\newblock {\em Phys. Rev. E}, 92(5):052704, 2015.

\bibitem{berry2020spectral}
T.~Berry and D.~Giannakis.
\newblock Spectral exterior calculus.
\newblock {\em Comm. Pure Appl. Math.}, 73(4):689--770, 2020.

\bibitem{bertalmio2001variational}
M.~Bertalm{\i}o, L.-T. Cheng, S.~Osher, and G.~Sapiro.
\newblock Variational problems and partial differential equations on implicit
  surfaces.
\newblock {\em J. Comput. Phys.}, 174(2):759--780, 2001.

\bibitem{Chan2017}
C.~H. Chan, M.~Czubak, and M.~M. Disconzi.
\newblock The formulation of the navier-stokes equations on riemannian
  manifolds.
\newblock {\em J. Geom. Phys.}, 121:335--346, 2017.

\bibitem{crane2010trivial}
K.~Crane, M.~Desbrun, and P.~Schr{\"o}der.
\newblock Trivial connections on discrete surfaces.
\newblock {\em Comput. Graph. Forum}, 29(5):1525--1533, 2010.

\bibitem{curtiss1952integration}
C.~F. Curtiss and J.~O. Hirschfelder.
\newblock Integration of stiff equations.
\newblock {\em Proc. Natl. Acad. Sci. USA}, 38(3):235--243, 1952.

\bibitem{desbrun2006discrete}
M.~Desbrun, E.~Kanso, and Y.~Tong.
\newblock Discrete differential forms for computational modeling.
\newblock {\em Int. Conf. Comput. Graph. Interactive Tech., ACM SIGGRAPH 2006
  Courses, 7}, pages 39--54, 2006.

\bibitem{do1992riemannian}
M.~P. Do~Carmo.
\newblock {\em Riemannian geometry}, volume~6.
\newblock Springer, 1992.

\bibitem{donoho2003hessian}
D.~L. Donoho and C.~Grimes.
\newblock Hessian eigenmaps: Locally linear embedding techniques for
  high-dimensional data.
\newblock {\em Proc. Natl. Acad. Sci. USA}, 100(10):5591--5596, 2003.

\bibitem{fisher2007design}
M.~Fisher, P.~Schr\"{o}der, M.~Desbrun, and H.~Hoppe.
\newblock Design of tangent vector fields.
\newblock {\em ACM Trans. Graph.}, 26(3):56–es, July 2007.

\bibitem{Natasha2015Solving}
N.~Flyer and B.~Fornberg.
\newblock Solving pdes with radial basis functions.
\newblock {\em Acta Numer.}, 24:215--258, 2015.

\bibitem{flyer2016role}
N.~Flyer, B.~Fornberg, V.~Bayona, and G.~A. Barnett.
\newblock On the role of polynomials in rbf-fd approximations: I. interpolation
  and accuracy.
\newblock {\em J. Comput. Phys.}, 321:21--38, 2016.

\bibitem{flyer2009radial}
N.~Flyer and G.~B. Wright.
\newblock A radial basis function method for the shallow water equations on a
  sphere.
\newblock {\em Proc. R. Soc. A.}, 465(2106):1949--1976, 2009.

\bibitem{Fuselier2009Stability}
E.~J. Fuselier and G.~B. Wright.
\newblock Stability and error estimates for vector field interpolation and
  decomposition on the sphere with rbfs.
\newblock {\em SIAM J. Numer. Anal.}, 47:3213--3239, 2009.

\bibitem{fuselier2013high}
E.~J. Fuselier and G.~B. Wright.
\newblock A high-order kernel method for diffusion and reaction-diffusion
  equations on surfaces.
\newblock {\em J. Sci. Comput.}, 56(3):535--565, 2013.

\bibitem{gillette2017finite}
A.~Gillette, M.~Holst, and Y.~Zhu.
\newblock Finite element exterior calculus for evolution problems.
\newblock {\em J. Comp. Math.}, 35(2):187--212, 2017.

\bibitem{greer2006improvement}
J.~B. Greer.
\newblock An improvement of a recent eulerian method for solving pdes on
  general geometries.
\newblock {\em J. Sci. Comput.}, 29(3):321--352, 2006.

\bibitem{gross2018hydrodynamic}
B.~J. Gross and P.~J. Atzberger.
\newblock Hydrodynamic flows on curved surfaces: Spectral numerical methods for
  radial manifold shapes.
\newblock {\em J. Comput. Phys.}, 371:663--689, 2018.

\bibitem{gross2020meshfree}
B.~J. Gross, N.~Trask, P.~Kuberry, and P.~J. Atzberger.
\newblock {Meshfree methods on manifolds for hydrodynamic flows on curved
  surfaces: A Generalized Moving Least-Squares (GMLS) approach}.
\newblock {\em J. Comput. Phys.}, 409:109340, 2020.

\bibitem{gross2018trace}
S.~Gross, T.~Jankuhn, M.~A. Olshanskii, and A.~Reusken.
\newblock A trace finite element method for vector-laplacians on surfaces.
\newblock {\em SIAM J. Numer. Anal.}, 56(4):2406--2429, 2018.

\bibitem{harlim2023radial}
J.~Harlim, S.~W. Jiang, and J.~W. Peoples.
\newblock Radial basis approximation of tensor fields on manifolds: from
  operator estimation to manifold learning.
\newblock {\em J. Mach. Learn. Res.}, 24(345):1--85, 2023.

\bibitem{hirani2003discrete}
A.~N. Hirani.
\newblock {\em Discrete exterior calculus, Ph.D. thesis}.
\newblock California Institute of Technology, 2003.

\bibitem{holme2007large}
R.~Holme.
\newblock 8.04 - large-scale flow in the core.
\newblock In G.~Schubert, editor, {\em Treatise on Geophysics}, pages 107--130.
  Elsevier, Amsterdam, 2007.

\bibitem{holst2012geometric}
M.~Holst and A.~Stern.
\newblock Geometric variational crimes: Hilbert complexes, finite element
  exterior calculus, and problems on hypersurfaces.
\newblock {\em Found. Comput. Math.}, 12:263--293, 2012.

\bibitem{jiang2024generalized}
S.~W. Jiang, R.~Li, Q.~Yan, and J.~Harlim.
\newblock Generalized finite difference method on unknown manifolds.
\newblock {\em J. Comput. Phys.}, 502:112812, 2024.

\bibitem{jones2023generalized}
A.~M. Jones, P.~A. Bosler, P.~A. Kuberry, and G.~B. Wright.
\newblock Generalized moving least squares vs. radial basis function finite
  difference methods for approximating surface derivatives.
\newblock {\em Comput. Math. Appl.}, 147:1--13, 2023.

\bibitem{kim1985application}
J.~Kim and P.~Moin.
\newblock Application of a fractional-step method to incompressible
  navier-stokes equations.
\newblock {\em J. Comput. Phys.}, 59(2):308--323, 1985.

\bibitem{lai2013local}
R.~Lai, J.~Liang, and H.~Zhao.
\newblock A local mesh method for solving pdes on point clouds.
\newblock {\em Inverse Probl. Imaging}, 7(3), 2013.

\bibitem{lee2018introduction}
J.~M. Lee.
\newblock {\em Introduction to Riemannian manifolds}.
\newblock Springer, 2018.

\bibitem{lehto2017radial}
E.~Lehto, V.~Shankar, and G.~B. Wright.
\newblock A radial basis function (rbf) compact finite difference (fd) scheme
  for reaction-diffusion equations on surfaces.
\newblock {\em SIAM J. Sci. Comput.}, 39(5):A2129--A2151, 2017.

\bibitem{levin1998approximation}
D.~Levin.
\newblock The approximation power of moving least-squares.
\newblock {\em Math. Comput.}, 67(224):1517--1531, 1998.

\bibitem{liang2012geometric}
J.~Liang, R.~Lai, T.~W. Wong, and H.~Zhao.
\newblock Geometric understanding of point clouds using laplace-beltrami
  operator.
\newblock In {\em 2012 IEEE conference on computer vision and pattern
  recognition}, pages 214--221. IEEE, 2012.

\bibitem{liang2013solving}
J.~Liang and H.-K. Zhao.
\newblock Solving partial differential equations on point clouds.
\newblock {\em SIAM J. Sci. Comput.}, 35(3):A1461--A1486, 2013.

\bibitem{liu1997moving}
W.-K. Liu, S.~Li, and T.~Belytschko.
\newblock Moving least-square reproducing kernel methods (i) methodology and
  convergence.
\newblock {\em Comput. Methods Appl. Mech. Engrg.}, 143(1-2):113--154, 1997.

\bibitem{mirzaei2012generalized}
D.~Mirzaei, R.~Schaback, and M.~Dehghan.
\newblock On generalized moving least squares and diffuse derivatives.
\newblock {\em IMA J. Numer. Anal.}, 32(3):983--1000, 2012.

\bibitem{monge1809application}
G.~Monge.
\newblock {\em Application de l'analyse {\`a} la g{\'e}om{\'e}trie {\`a}
  l'usage de l'Ecole imp{\'e}riale polytechnique}.
\newblock Veuve Bernard, 1809.

\bibitem{morita2001geometry}
S.~Morita.
\newblock {\em Geometry of differential forms}.
\newblock Number 201. American Mathematical Soc., 2001.

\bibitem{nayroles1992generalizing}
B.~Nayroles, G.~Touzot, and P.~Villon.
\newblock Generalizing the finite element method: diffuse approximation and
  diffuse elements.
\newblock {\em Comput. Mech.}, 10(5):307--318, 1992.

\bibitem{nitschke2019hydrodynamic}
I.~Nitschke, S.~Reuther, and A.~Voigt.
\newblock Hydrodynamic interactions in polar liquid crystals on evolving
  surfaces.
\newblock {\em Phys. Rev. Fluids}, 4(4):044002, 2019.

\bibitem{peoples2024higher}
J.~W. Peoples and J.~Harlim.
\newblock A higher order local mesh method for approximating laplacians on
  unknown manifolds.
\newblock {\em arXiv preprint arXiv:2405.15735}, 2024.

\bibitem{petras2018rbf}
A.~Petras, L.~Ling, and S.~J. Ruuth.
\newblock An rbf-fd closest point method for solving pdes on surfaces.
\newblock {\em J. Comput. Phys.}, 370:43--57, 2018.

\bibitem{piret2012orthogonal}
C.~Piret.
\newblock The orthogonal gradients method: A radial basis functions method for
  solving partial differential equations on arbitrary surfaces.
\newblock {\em J. Comput. Phys.}, 231(14):4662--4675, 2012.

\bibitem{pressley2010elementary}
A.~N. Pressley.
\newblock {\em Elementary differential geometry}.
\newblock Springer Science \& Business Media, 2010.

\bibitem{prycel1993SL}
J.~D. Pryce.
\newblock {\em Numerical Solution of Sturm-Liouville Problems}.
\newblock Oxford University Press, 1993.

\bibitem{ruuth2008simple}
S.~J. Ruuth and B.~Merriman.
\newblock A simple embedding method for solving partial differential equations
  on surfaces.
\newblock {\em J. Comput. Phys.}, 227(3):1943--1961, 2008.

\bibitem{samavaki2020navier}
M.~Samavaki and J.~Tuomela.
\newblock Navier--stokes equations on riemannian manifolds.
\newblock {\em J. Geom. Phys.}, 148:103543, 2020.

\bibitem{Schoberl2014C11IO}
J.~Sch{\"o}berl.
\newblock C++11 implementation of finite elements in ngsolve.
\newblock {\em Institute for analysis and scientific computing, Vienna
  University of Technology}, 2014.

\bibitem{shankar2015radial}
V.~Shankar, G.~B. Wright, R.~M. Kirby, and A.~L. Fogelson.
\newblock A radial basis function (rbf)-finite difference (fd) method for
  diffusion and reaction--diffusion equations on surfaces.
\newblock {\em J. Sci. Comput.}, 63(3):745--768, 2015.

\bibitem{sharp2019vector}
N.~Sharp, Y.~Soliman, and K.~Crane.
\newblock The vector heat method.
\newblock {\em ACM Trans. Graph.}, 38(3):19, 2019.

\bibitem{singer2012vector}
A.~Singer and H.-T. Wu.
\newblock Vector diffusion maps and the connection laplacian.
\newblock {\em Comm. Pure Appl. Math.}, 65(8):1067--1144, 2012.

\bibitem{sober2020manifold}
B.~Sober and D.~Levin.
\newblock Manifold approximation by moving least-squares projection (mmls).
\newblock {\em Constr. Approx.}, 52(3):433--478, 2020.

\bibitem{stam2003flows}
J.~Stam.
\newblock Flows on surfaces of arbitrary topology.
\newblock {\em ACM Trans. Graph.}, 22(3):724--731, 2003.

\bibitem{suchde2021meshfree}
P.~Suchde.
\newblock A meshfree lagrangian method for flow on manifolds.
\newblock {\em Internat. J. Numer. Methods Fluids}, 93(6):1871--1894, 2021.

\bibitem{suchde2019meshfree}
P.~Suchde and J.~Kuhnert.
\newblock {A meshfree generalized finite difference method for surface PDEs}.
\newblock {\em Comput. Math. Appl.}, 78(8):2789--2805, 2019.

\bibitem{suchde2019fully}
P.~Suchde and J.~Kuhnert.
\newblock A fully lagrangian meshfree framework for pdes on evolving surfaces.
\newblock {\em J. Comput. Phys.}, 395:38--59, 2019.

\bibitem{taylorPDE1}
M.~Taylor.
\newblock {\em Partial Differential Equations I: Basic Theory.}, volume 115 of
  {\em Applied Mathematical Sciences}.
\newblock Springer, New York, 2nd edition, 2011.

\bibitem{taylorPDE3}
M.~Taylor.
\newblock {\em Partial Differential Equations III: Nonlinear Equations.},
  volume 117 of {\em Applied Mathematical Sciences}.
\newblock Springer, New York, 2nd edition, 2011.

\bibitem{tyagi2013tangent}
H.~Tyagi, E.~Vural, and P.~Frossard.
\newblock Tangent space estimation for smooth embeddings of riemannian
  manifolds.
\newblock {\em Inf. Inference}, 2(1):69--114, 2013.

\bibitem{vaxman2016directional}
A.~Vaxman, M.~Campen, O.~Diamanti, D.~Panozzo, D.~Bommes, K.~Hildebrandt, and
  M.~Ben-Chen.
\newblock Directional field synthesis, design, and processing.
\newblock In {\em Computer graphics forum}, volume~35, pages 545--572. Wiley
  Online Library, 2016.

\bibitem{Wendland2005Scat}
H.~Wendland.
\newblock {\em Scattered Data Approximation}.
\newblock Cambridge University Press, 2005.

\bibitem{xu2003eulerian}
J.-J. Xu and H.-K. Zhao.
\newblock An eulerian formulation for solving partial differential equations
  along a moving interface.
\newblock {\em J. Sci. Comput.}, 19(1):573--594, 2003.

\bibitem{zhang2004principal}
Z.~Zhang and H.~Zha.
\newblock Principal manifolds and nonlinear dimensionality reduction via
  tangent space alignment.
\newblock {\em SIAM J. Sci. Comput.}, 26(1):313--338, 2004.

\end{thebibliography}

\end{document}